\documentclass[11pt,oneside]{amsart}

\usepackage{fullpage}

\usepackage{amsthm,amsfonts,amssymb,amsmath,amsxtra,mathtools}
\usepackage[dvipsnames]{xcolor}
\usepackage[all]{xy}
\SelectTips{cm}{}
\usepackage{tikz}
\usepackage{verbatim}
\usepackage{todonotes}
\usepackage{mathrsfs}
\usepackage{booktabs,makecell}
\usepackage{diagbox}
\usepackage{epigraph}
\usepackage{bm}
\usepackage{tikz-cd}
\usepackage{relsize}
\usepackage{braket}
\usepackage{lmodern}
\usepackage{array}
\usepackage[pagebackref,breaklinks]{hyperref}
\usepackage{appendix} 
\usepackage[lite,abbrev,msc-links,alphabetic]{amsrefs}
\usepackage{indentfirst}
\usepackage{yhmath}

\usepackage{mathdots}

\usepackage{dsfont}
\allowdisplaybreaks

\RequirePackage{xspace}
\RequirePackage{etoolbox}
\RequirePackage{varwidth}
\RequirePackage{enumitem}
\RequirePackage{tensor}
\RequirePackage{longtable}
\RequirePackage{multirow}
\RequirePackage{bigints}

\hypersetup{  
	colorlinks=true,
	linkcolor=blue,
	citecolor=MidnightBlue,
	urlcolor=magenta,
}

\DeclareMathOperator{\Hom}{Hom}
\DeclareMathOperator{\End}{End}

\DeclareMathOperator{\Lie}{Lie}

\newcommand{\val}{{\mathrm{val}}}

\renewcommand{\O}{\mathrm{O}}

\newcommand{\isoarrow}{%
	\ifbool{@display}{\overset{\sim}{\longrightarrow}}{\xrightarrow\sim}%
}

\theoremstyle{plain}
\newtheorem{theorem}{Theorem}[section]
\newtheorem{lemma}[theorem]{Lemma}
\newtheorem{proposition}[theorem]{Proposition}

\newtheorem{conjecture}[theorem]{Conjecture}

\theoremstyle{definition}
\newtheorem{definition}[theorem]{Definition}
\newtheorem{remark}[theorem]{Remark}

\theoremstyle{remark}

\newcommand{\tx}[1]{\text{#1}}

\newcommand{\BC}{{\mathbb C}}

\newcommand{\BE}{{\mathbb E}}
\newcommand{\BF}{{\mathbb F}}

\newcommand{\BL}{{\mathbb L}}

\newcommand{\BQ}{{\mathbb Q}}

\newcommand{\BV}{{\mathbb V}}

\newcommand{\BX}{{\mathbb X}}

\newcommand{\BZ}{{\mathbb Z}}

\newcommand{\CA}{{\mathcal A}}

\newcommand{\CE}{{\mathcal E}}
\newcommand{\CF}{{\mathcal F}}
\newcommand{\CG}{{\mathcal G}}

\newcommand{\CI}{{\mathcal I}}

\newcommand{\CN}{{\mathcal N}}

\newcommand{\CR}{{\mathcal R}}
\newcommand{\CS}{{\mathcal S}}
\newcommand{\CT}{{\mathcal T}}

\newcommand{\CY}{{\mathcal Y}}
\newcommand{\CZ}{{\mathcal Z}}

\newcommand{\ScS}{{\mathscr S}}

\newcommand{\alphad}{\mathrm{Den}}

\DeclareFontFamily{U}{matha}{\hyphenchar\font45}
\DeclareFontShape{U}{matha}{m}{n}{
	<5> <6> <7> <8> <9> <10> gen * matha
	<10.95> matha10 <12> <14.4> <17.28> <20.74> <24.88> matha12
}{}
\DeclareSymbolFont{matha}{U}{matha}{m}{n}
\DeclareFontFamily{U}{mathx}{\hyphenchar\font45}
\DeclareFontShape{U}{mathx}{m}{n}{
	<5> <6> <7> <8> <9> <10>
	<10.95> <12> <14.4> <17.28> <20.74> <24.88>
	mathx10
}{}
\DeclareSymbolFont{mathx}{U}{mathx}{m}{n}

\DeclareMathSymbol{\obot}         {2}{matha}{"6B}

\newcommand{\C}{\mathbb{C}}
\newcommand{\h}{\mathbb{H}}

\newcommand{\cZ}{\mathcal{Z}}

\newcommand{\cY}{\mathcal{Y}}

\newcommand{\cN}{\mathcal{N}}

\newcommand{\cM}{\mathcal{M}}

\newcommand{\Spf}{\mathrm{Spf}\, }

\newcommand{\diag}{\mathrm{Diag}}

\numberwithin{equation}{section}

\title{Arithmetic Siegel-Weil Formulas for Unitary Shimura Varieties with Iwahori Levels}

\author[Sungyoon Cho]{Sungyoon Cho}
\address{Department of Mathematics, POSTECH, 77, Cheongam-ro, Nam-gu, Pohang-si, Gyeongsangbuk-do, Korea}
\email{sungyooncho@postech.ac.kr}

\author[Qiao He]{Qiao He}
\address{Zhizhen School of Interdisciplinary Mathematical Sciences, Chinese University of Hong Kong}
\email{qiaohe@cuhk.edu.hk}

\author[Chul-hee Lee]{Chul-hee Lee}
\address{Center for Artificial Intelligence and Natural Sciences, Korea Institute for Advanced Study, 85 Hoegiro, Dongdaemun-gu, Seoul, Korea}
\email{chlee@kias.re.kr}

\author[Zhiyu Zhang]{Zhiyu Zhang}
\address{Zhizhen School of Interdisciplinary Mathematical Sciences, Chinese University of Hong Kong}
\email{zhiyuzhang@cuhk.edu.hk}

\date{\today}
\begin{document}
	\begin{abstract}
	In this paper, we formulate conjectural formulas for the arithmetic intersection numbers of special cycles on unitary Shimura varieties with Iwahori level structure. We first introduce a framework for formulating Kudla-Rapoport conjectures: we develop the theory of weighted representation densities and special cycles for general parahoric level structures, prove dimension-reduction formulas on both the analytic and geometric sides and duality relations for special cycles. We then formulate the Kudla-Rapoport conjecture for Iwahori level structures based on these properties. Finally, we prove the conjecture when $n=2$.
	\end{abstract}

	\maketitle

	\tableofcontents

	\section{Introduction}
	 
	 In \cite{KR1} and \cite{KR2}, Kudla and Rapoport formulated a conjectural formula for the arithmetic intersection numbers of special cycles on unitary Shimura varieties with hyperspecial level structure. This conjecture, the so--called Kudla-Rapoport conjecture, was proved by Li and Zhang in \cite{LZ}. After this, the Kudla--Rapoport conjecture was formulated and proved for various level structures (\cite{Cho}, \cite{LL2}, \cite{HSY3}, \cite{HLSY}, \cite{chz}, \cite{Yao}, \cite{Luo}).
	 
	 In this paper, we formulate the Kudla--Rapoport conjecture for unitary Shimura varieties with Iwahori level structure at an unramified prime. More precisely, we formulate a conjectural formula for the arithmetic intersection numbers of special cycles on unitary Rapoport--Zink spaces associated with unitary Shimura varieties with Iwahori level structure.
	 
	 We now describe our results in more detail. We fix an odd prime $p$. Let $F_0$ be a finite unramified extension of $\BQ_p$ with ring of integers $O_{F_0}$, and residue field $\BF_q$. Let $F$ be the unramified quadratic extension of $F_0$ with ring of integers $O_F$, and residue field $\BF_{q^2}$. We fix a uniformizer $\pi$ of $O_F$ and $O_{F_0}$. Let $\breve F$ be the completion of the maximal unramified extension of $F$ with ring of integers $O_{\breve{F}}$, and residue field $\BF$. We denote by $^*$ the nontrivial Galois automorphism of $F$ over $F_0$.
	 
	 Let $n$ be a positive integer and let $\mathds{h}=(h_1,\dots,h_m)$ a tuple of integers such that $0\leq h_1<\dots<h_m \leq n$ and $h_1 \equiv \dots \equiv h_m (\tx{mod } 2)$. Then, we can consider the unitary Rapoport-Zink space $\CN^{[\mathds{h}]}_n$ over $\Spf O_{\breve{F}}$ (Definition \ref{definition3.3.2}). Here, $\mathds{h}$ is related to the general parahoric level structure of the corresponding unitary Shimura variety.  We note that there is a dual isomorphism $$\theta:\CN_n^{[\mathds{h}]} \rightarrow \CN_{n}^{[\mathds{h}^{\vee}]},$$ between the two Rapoport-Zink spaces, where $\mathds{h}^{\vee}=(n-h_m,\dots,n-h_1)$ (Definition \ref{definition3.3.3}). Let $\BV$ be the space of special homomorphisms (Definition \ref{definition3.3.4}). For each nonzero special homomorphism $x \in \BV$, and $h_i$, $1 \leq i \leq m$, we can attach Cartier divisors $\CZ_{h_i}(x)$ and $\CY_{h_i}(x)$ (Definition \ref{definition3.3.6}). We note that there are dual relations between $\CZ_{h_i}$ and $\CY_{n-h_i}$ via the dual isomorphism $\theta:\CN_n^{[\mathds{h}]} \rightarrow \CN_{n}^{[\mathds{h}^{\vee}]}$.
	 
	  We denote by $\CE_{n}^{\mathds{h}}$ the set of all $n$-tuples of integers $\alpha=(\alpha_1,\dots,\alpha_n)$ such that $-n \leq \alpha_1 \leq  \dots \leq \alpha_n \leq n$, $\alpha_1 \equiv \dots \equiv \alpha_n \pmod{2}$, and for all $i$, $\alpha_i =-h_j \tx{ or } h_j$ for some $j$. Fix $\alpha \in \CE_n^{\mathds{h}}$ and assume that $\alpha_1 \leq \dots \leq \alpha_{h-1} < 0 \leq \alpha_h \leq \dots \leq \alpha_n$. Let $x_1,\dots,x_{h-1}, y_1,\dots,y_{n-h+1}$ be special homomorphisms in $\BV$. Then, we consider special cycles
	  \begin{equation*}
	  	\CZ_{-\alpha_1}(x_1),\dots,\CZ_{-\alpha_{h-1}}(x_{h-1}),\CY_{\alpha_h}(y_1),\dots,\CY_{\alpha_n}(y_{n-h+1}),
	  \end{equation*}
	  in $\CN^{[\mathds{h}]}_n$. We denote by $$\CI_{\alpha}^{\mathds{h}}(x_1,\dots,x_{h-1},y_1,\dots,y_{n-h+1})$$ the arithmetic intersection number of $\CZ_{-\alpha_1}(x_1)$,$\CZ_{-\alpha_2}(x_2)$,$\dots$,$\CZ_{-\alpha_{h-1}}(x_{h-1})$,$\CY_{\alpha_h}(y_1)$,$\dots$,$\CY_{\alpha_n}(y_{n-h+1})$ in $\CN^{[\mathds{h}]}_n$ (Definition \ref{definition3.3.10}).
	  
	  The Kudla--Rapoport conjecture relates the arithmetic intersection number $\CI_{\alpha}^{\mathds{h}}(x_1,$ $\dots,$ $x_{h-1}$, $y_1,$ $\dots$, $y_{n-h+1})$ to weighted representation densities, which we now define.
	  
	 Denote by $\CE_{n}$ the set of all $n$-tuples of integers $\alpha=(\alpha_1,\dots,\alpha_n)$ such that $-n \leq \alpha_1 \leq  \dots \leq \alpha_n \leq n$, $\alpha_1 \equiv \dots \equiv \alpha_n \pmod{2}$. For each $\alpha \in \CE_{n}$, we define a weighted representation density $W_{\alpha}(B,r)$, where $B$ is a hermitian matrix and $r$ is a nonnegative integer (Definition \ref{definition2.1}).
	 
	 We prove dimension reduction formulas for weighted representation densities (Proposition \ref{proposition2.10} and Proposition \ref{proposition2.25}) and for special cycles (Proposition \ref{proposition3.9} and Proposition \ref{proposition3.11}).
     
     We expect the Kudla--Rapoport conjecture for general parahoric levels to form a system of conjectures compatible with the dimension reduction formulas and to take the following form: for $\alpha \in \CE_{n}^{\mathds{h}}$,
	 \begin{equation}\label{eq:conj}
	 	\CI_{\alpha}^{\mathds{h}}(x_1,\dots,x_{h-1},y_1,\dots,y_{n-h+1})=C_{\mathds{h}}\lbrace W_{\alpha}'(B)+\tx{correction terms}\rbrace,
	 \end{equation}
	 where $B$ is the hermitian matrix
	 \begin{equation*}
	 	B=	\left(\begin{array}{cc} 
	 		h(x_i,x_j) & h(x_i,y_l)\\
	 		h(y_k,x_j)& h(y_k,y_l)
	 	\end{array} 
	 	\right)_{\substack{1\leq i,j \leq h-1\\ 1 \leq k,l \leq n-h+1}},
	 \end{equation*}
	 and $C_{\mathds{h}}$ is a constant related to the volume, with respect to our chosen measure, of the parahoric group corresponding to $\mathds{h}$. We expect that the correction terms are a linear combination of the $W_{\beta}(B)$, i.e., the non-derivative weighted representation densities.
	 
	Now, assume that $\mathds{h}=(1,3,\dots,2n-1)$, $\alpha=(-2n+1,-2n+3,\dots,-1,1,3,\dots,2n-1) \in \CE_{2n}^{\mathds{h}}$. Then, $\CN_{2n}^{[\mathds{h}]}$ is the Rapoport--Zink space associated with the unitary Shimura variety with Iwahori level structure. Also, $\CI_{\alpha}^{\mathds{h}}(x_1,\dots,x_{n},y_1,\dots,y_{n})$ is the arithmetic intersection number of special cycle $\CZ_{2n-1}(x_1),\dots,\CZ_{1}(x_{n}),\CY_{1}(y_1),\dots,\CY_{2n-1}(y_{n})$.
	 In particular, the types of the special cycles are all distinct. In this case, we conjecture the following precise identity.
	 
	 \begin{conjecture}[Theorem \ref{theorem3.1}, Conjecture \ref{conjecture4.3}]\label{conjectureintro} 	Let $\mathds{h}=(1,3,\dots,2n-1)$ and $\alpha=(-2n+1,-2n+3,\dots,-1,1,\dots,2n-1) \in \CE_{2n}^{\mathds{h}}$. Let $\lbrace x_1$,$\dots$,$x_{n}$,$y_1$,$\dots$,$y_{n} \rbrace$ be a basis of the space of special homomorphisms $\BV$. Then, we have
	 		\begin{equation*}\begin{array}{ll}
	 			\CI_{\alpha}^{[\mathds{h}]}(x_1,\dots,x_{n},y_1,\dots,y_{n})
	 			&=\dfrac{q^{6n^2-n}}{(q+1)^2(q^2-1)^{n-1}}\lbrace W_{\alpha}'(B)-q^{-4n^2+2n}\mathlarger{\sum}_{\beta}C_{\beta,2n-1-\beta_{2n}}W_{\beta}(B,0) \rbrace,
	 		\end{array}
	 	\end{equation*}
	 	where $C_{\beta,i}$ is defined by
	 \begin{equation*}
	 	C_{\beta,i}:=\left\lbrace \begin{array}{ll}
	 		q^{4n^2-2n\delta_{\beta}-\sum_{i=1}^{2n} \beta_i}(-1)^k\lbrace \dfrac{1}{q^{2k-1}+1}+\dfrac{1}{q^{2k}-1} \rbrace & \tx{if }i=2k-1,\\
	 		q^{4n^2-2n\delta_{\beta}-\sum_{i=1}^{2n} \beta_i}(-1)^k\lbrace \dfrac{1}{q^{2k-2}-1}+\dfrac{1}{q^{2k-1}+1} \rbrace & \tx{if }i=2k-2,\\0 &\tx{if } i=0.
	 	\end{array}\right.
	 \end{equation*}
	 Here, $\beta$ runs over the subset of $\CE_{2n}$ satisfying the following conditions:
	 \begin{enumerate}
	 	\item $\beta_i=-\beta_{2n+1-i}$ for all $1 \leq i \leq 2n$.
	 	\item $\alpha-\beta$ is increasing, i.e., $\alpha_i-\beta_i \leq \alpha_{i+1}-\beta_{i+1}$ for all $1 \leq i \leq 2n-1$.
	 \end{enumerate}
	 \end{conjecture}
	 
	 We note that, by the dimension reduction formulas, Conjecture \ref{conjectureintro} gives conjectural formulas for the arithmetic intersection numbers of special cycles at all Iwahori levels and in arbitrary dimension, when the types of the special cycles are all distinct. We refer to Remark \ref{remark4.5} for details.  We also prove that the conjecture holds when the Rapoport-Zink space is $\CN_2^{[\mathds{h}]}$, i.e., in dimension $2$. We refer to Section \ref{section5} for details.
	 
	 We now explain how we formulate Conjecture \ref{conjectureintro}. To determine the correction terms, we use the duality of special cycles (Lemma \ref{lemma3.12}). More precisely, we use duality relations between $\CZ_{h_i}$ and $\CY_{n-h_i}$ via the dual isomorphism $\theta:\CN_n^{[\mathds{h}]} \rightarrow \CN_n^{[\mathds{h}^{\vee}]}$. Since $\theta$ is an isomoprhism, it preserves the intersection numbers of special cycles. Therefore, the analytic side $C_{\mathds{h}}\lbrace W_{\alpha}'(B)+\tx{correction terms}\rbrace$ of \eqref{eq:conj} should also be invariant under this change (Remark \ref{remark4.1}). We prove that this condition is sufficient to determine the correction terms when $\mathds{h}=(1,3,\dots,2n-1)$ and $\alpha=(-2n+1,-2n+3,\dots,-1,1,\dots,2n-1) \in \CE_{2n}^{\mathds{h}}$ (Theorem \ref{theorem3.1}). This is how we formulate Conjecture \ref{conjectureintro}.
	 
	 Since Lemma \ref{lemma3.12} holds for arbitrary $\mathds{h}$ and $\alpha$, it is natural to expect that our strategy of formulating a system of conjectures via the duality of special cycles can also be applied at any parahoric level and for any choice of special cycles. For this, we refer the reader to Remark \ref{remark4.6} and Remark \ref{remark4.7}. Although this strategy works in many, if not all, cases and the formulas have many interesting properties, it remains difficult to give a uniform formulation because of the large number of cases involved (e.g., $1386$ cases for $n=6$, $19305$ cases for $n=8$). We hope to return to this problem in the future.

	\subsection*{Notation}
	\begin{enumerate}
		\item We fix a prime $p>2$. Let $F_0$ be a finite unramified extension of $\BQ_p$ with ring of integers $O_{F_0}$, and residue field $\BF_q$. Let $F$ be the unramified quadratic extension of $F_0$ with ring of integers $O_F$, and residue field $\BF_{q^2}$. We fix a uniformizer $\pi$ of $O_F$ and $O_{F_0}$. Let $\breve F$ be the completion of the maximal unramified extension of $F$ with ring of integers $O_{\breve{F}}$, and residue field $\BF$. We denote by $^*$ the nontrivial Galois automorphism of $F$ over $F_0$.
		\item For a lattice $L$ with hermitian form $(\text{ },\text{ })$, we denote by $L^{\vee}$ the dual lattice with respect to $(\text{ },\text{ })$.
	\end{enumerate}

	\subsection*{Acknowledgements} We would like to thank Chao Li, Yifeng Liu, Michael Rapoport, and Wei Zhang for helpful discussions and comments. S.C. was supported by the National Research Foundation of Korea (NRF) grant funded by the Korea government (MSIT) (No. RS-2026-25508638), (No. RS-2026-25485874), and (No. RS -2021-NR060139).
    C.-h. L. was supported by a KIAS Individual Grant (AP067301) via the Center for Artificial Intelligence and Natural Sciences at Korea Institute for Advanced Study. Q. H. and  Z. Z. are supported by a research start up fund of the Chinese University of Hong Kong. Q.H. was supported by the AMS-Simons Travel Grant Program during the preparation of this work.
 	
		\subsection*{Declaration of AI use}
	The authors used AI solely to check for typographical errors and improve the language of the manuscript.
    
	\section{Weighted representation densities}\label{sec: weighted rep}
	\subsection{Weighted representation densities}\label{subsection2.1}
	In this subsection, we recall the definition of weighted representation densities and prove some properties. We write $^*$ for the nontrivial Galois automorphism of $F$ over $F_0$. We fix the standard additive character $\psi:F_0 \rightarrow \BC^{\times}$ which is trivial on $O_{F_0}$. Let $V^+$ (resp. $V^{-}$) be a split (resp. non-split) $n$-dimensional hermitian vector space over $F$. We write $\CS((V^{\pm})^{n})$ for the space of Schwartz functions on $(V^{\pm})^{n}$. We write $V_{r,r}$ for the split hermitian space of signature $(r,r)$ and $L_{r,r}$ for a self-dual lattice in $V_{r,r}$. Let $\phi_{r,r}$ be the characteristic function of $(L_{r,r})^{n}$. Let $(V^{\pm})^{[r]}$ be the space $V^{\pm} \otimes V_{r,r}$. For any function $\phi \in \CS((V^{\pm})^{n})$, we define a function $\phi^{[r]}$ by $\phi \otimes \phi_{r,r} \in \CS(((V^{\pm})^{[r]})^{n})$.
	
	We write $\Gamma_{n}$ for the Iwahori subgroup
	\begin{equation*}
		\Gamma_{n}\coloneqq \lbrace \gamma=(\gamma_{ij}) \in GL_{n}(O_{F}) \mid \gamma_{ij} \in \pi O_{F} \text{ if }i>j \rbrace.
	\end{equation*}
	
	We define the set $V_{n}(F)$ and $X_{n}(F)$ as follows:
	\begin{equation*}
		V_{n}(F)=\lbrace Y \in M_{n,n}(F) \mid \text{}^tY^*=Y \rbrace,
	\end{equation*}
	\begin{equation*}
		X_{n}(F)=\lbrace X \in GL_{n}(F) \mid \text{}^tX^*=X \rbrace.
	\end{equation*}
	
	For $g \in GL_{n}(F)$ and $X \in X_{n}(F)$, we define the group action of $GL_{n}(F)$ on $X_{n}(F)$ by $g \cdot X=gX^tg^*$. For $X, Y \in V_{n}(F)$, we denote by $\langle X,Y \rangle=\mathrm{Tr}(XY)$. For $X \in M_{m,n}(F)$ and $A \in V_{m}(F)$, we denote by $A[X]=^tX^*AX.$ For a hermitian matrix $A \in X_m(F)$, we define
	\begin{equation*} A^{[r]}=\left(\begin{array}{cc} A & \\ & I_{2r} \end{array}\right).\end{equation*}

	Now, let us recall the definition of weighted representation densities in \cite{Cho}.
	\begin{definition}[cf. {\cite[Definition 3.1]{Cho}}]\label{definition2.1} Let $t,n$ be integers such that $-n \leq t \leq n$. \begin{enumerate}
			\item We denote by $\CE_n$ the set of all $n$-tuples of integers $\alpha=(\alpha_1,\dots,\alpha_n)$ such that  $-n \leq \alpha_1 \leq \dots \leq \alpha_n \leq n$ and $\alpha_1 \equiv \dots \equiv \alpha_n \pmod{2}$. We define $\alpha \in \CE_n$ to be even (resp. odd) if all $\alpha_i$ are even (resp. odd). Also, we define
			\begin{equation*}
				\delta_{\alpha}:=\left\lbrace \begin{array}{ll} 0 & \text{ if } n-\alpha \text{ is even},\\  1 & \text{ if } n-\alpha \text{ is odd}. \end{array}\right.
			\end{equation*}

			\item Assume that $n$ is even. Let $\Lambda_{n,n}$ (resp. $\Lambda_{n-1,n}$) be a lattice of rank $n$ in $V^+$ (resp. in $V^{-}$) with hermitian form 
			\begin{equation*}
				A_{n,0}\coloneqq \scalebox{0.8}{$\left(\begin{array}{ccc}  && \pi^{-1}\\
						& \iddots &\\
						\pi^{-1} &&
					\end{array}\right)$} \mathlarger{(}\text{resp. }	A_{n,1}\coloneqq \scalebox{0.8}{$\left(\begin{array}{cccc} 1  &&&\\ &&&\pi^{-1}\\
						&& \iddots &\\
						&\pi^{-1} &&
					\end{array}\right)$} \mathlarger{)},
			\end{equation*}
			with respect to a basis $\lbrace b_1, b_2,\dots, b_n \rbrace$.
			
			\item Assume that $n$ is odd. Let $\Lambda_{n,n}$ (resp. $\Lambda_{n-1,n}$) be a lattice of rank $n$ in $V^-$ (resp. in $V^{+}$) with hermitian form 
			\begin{equation*}
				A_{n,0}\coloneqq \scalebox{0.8}{$\left(\begin{array}{ccc}  && \pi^{-1}\\
						& \iddots &\\
						\pi^{-1} &&
					\end{array}\right)$} \mathlarger{(}\text{resp. }	A_{n,1}\coloneqq \scalebox{0.8}{$\left(\begin{array}{cccc} 1  &&&\\ &&&\pi^{-1}\\
						&& \iddots &\\
						&\pi^{-1} &&
					\end{array}\right)$} \mathlarger{)},
			\end{equation*}
			with respect to a basis $\lbrace b_1, b_2,\dots, b_n \rbrace$.
			
			\item For $t \equiv n \pmod{2}$, we write $\Lambda_{t,n}$ for the sublattice of $\Lambda_{n,n}$ with basis \begin{equation*}\lbrace b_1,\dots,b_{(n+t)/2},\pi b_{(n+t+2)/2},\dots,\pi b_n \rbrace.\end{equation*}
			
			\item For $t \not\equiv n \pmod{2}$, we write $\Lambda_{t,n}$ for the sublattice of $\Lambda_{n-1,n}$ with basis \begin{equation*} \lbrace b_1,\dots,b_{(n+t+1)/2},\pi b_{(n+t+3)/2},\dots,\pi b_n \rbrace\end{equation*}.

			\item For $\alpha=(\alpha_1,\dots,\alpha_n) \in \CE_n$, we denote by $1_{\alpha} \in \CS((V^{\pm})^n)$ the characteristic function of
			\begin{equation*}
				\Lambda_{\alpha_1,n}\times \Lambda_{\alpha_2,n}\times \dots \times \Lambda_{\alpha_n,n}.
			\end{equation*}
			
			\item 	For $B \in X_{n}(F)$ and $\alpha=(\alpha_1,\dots,\alpha_n)$, we define
			\begin{equation*}
				W_{\alpha}(B,r)\coloneqq \int_{V_{n}(F)}\int_{M_{n+2r,n}(F)}\psi(\langle Y,A_{n,\delta_{\alpha}}^{[r]}[X]-B\rangle)1_{\alpha}^{[r]}(X)dXdY.
			\end{equation*}
			
			Here, $dY$ (resp. $dX$) is the Haar measure on $V_{n}(F)$ (resp. $M_{n+2r,n}(F)$) such that
			\begin{equation*}
				\int_{V_{n}(O_F)}dY=1 \text{ }(\text{resp. }\int_{M_{n+2r,n}(O_F)}dX=1 ).
			\end{equation*}
			
			\item For $A \in X_m(O_F)$ and $B \in X_n(O_F)$, we define $\alphad(A,B)$ by
			\begin{equation*}
				\alphad(A,B)=\lim_{d \rightarrow \infty} (q^{-d})^{n(2m-n)}\vert \lbrace x \in M_{m,n}(O_F/\pi^dO_F) \mid A[x]\equiv B (\text{mod} \pi^d)\rbrace\vert.
			\end{equation*}
		\end{enumerate}

	\end{definition}
	
	\begin{remark}
		Note that for $0 \leq t \leq n$, $\Lambda_{t,n}$ is a vertex lattice satisfying $\Lambda_{-t,n}=\Lambda_{t,n}^{\vee} \overset{t}{\subset} \Lambda_{t,n} \overset{n-t}{\subset} \pi^{-1}\Lambda_{t,n}^{\vee}$.
	\end{remark}

	\begin{proposition}[cf. {\cite{Hi}}, {\cite[Lemma 3.5]{Cho}}]\label{proposition2.3}  For $A \in X_{m}(F)$, $B \in X_{n}(F)$ and $\alpha=(\alpha_1,\dots, \alpha_n)$, we have
		\begin{equation*}
			\alphad(A,B)=\mathlarger{\sum}_{Y \in \Gamma_{n} \setminus X_{n}(F)}\dfrac{\CG(Y,B)\CF(Y,A)}{\alpha(Y;\Gamma_{n})},
		\end{equation*}
		
		\begin{equation*}
			W_{\alpha}(B,r)=\mathlarger{\sum}_{Y \in \Gamma_{n} \setminus X_{n}(F)}\dfrac{\CG(Y,B)\CF_{\alpha}(Y,r)}{\alpha(Y;n)}.
		\end{equation*}
		
		Here, we define $\CF_{\alpha}(Y,r)$ and $\CG(Y,B)$ as follows:
		\begin{equation*}\begin{array}{l}
				\CF(Y,A)\coloneqq \mathlarger{\int}_{M_{m,n}(F)}\psi(\langle Y,A[X] \rangle)dX,\\
				\CF_{\alpha}(Y,r)\coloneqq \mathlarger{\int}_{M_{n+2r,n}(F)}\psi(\langle Y,A_{n,\delta_{\alpha}}^{[r]}[X] \rangle)1_{\alpha}^{[r]}(X)dX,\\
				\CG(Y,B)\coloneqq \mathlarger{\int}_{\Gamma_{n}}\psi(\langle Y, -B[\gamma]\rangle)d\gamma,
			\end{array}
		\end{equation*}
		where $d\gamma$ is the Haar measure on $M_{n,n}(O_F)$ such that $\int_{M_{n,n}(O_F)}d\gamma=1$.
		
		Also, we define $\alpha(Y;n)$ by $\alpha(Y;n)\coloneqq \lim_{d\rightarrow \infty} q^{-dn^2}N_d(Y;\Gamma_{n}),$ where $N_d(Y;\Gamma_{n})=\vert \lbrace \gamma \in \Gamma_{n}( \text{mod } \pi^d) \vert \gamma \cdot Y \equiv Y ( \text{mod } \pi^{d})\rbrace \vert$.

	\end{proposition}
	
	\begin{definition} For $r \geq 0$, $\alpha \in \CE_n$, and $B \in X_n(O_F)$, we can regard $\CF_{\alpha}(Y,r)$ and $W_{\alpha}(B,r)$ as functions of $X=(-q)^{-2r}$. We define
		\begin{equation*}
			\CF'_{\alpha}(Y)\coloneqq -\dfrac{d}{dX}\CF_{\alpha}(Y,r)\vert_{X=1}, \text{ and }W'_{\alpha}(B)\coloneqq -\dfrac{d}{dX}W_{\alpha}(B,r)\vert_{X=1}.
		\end{equation*}
	\end{definition}
	
	\begin{definition}\label{definition2.5}
		Assume that $0 \leq h \leq n$.
		\begin{enumerate}
			\item We define the set $\CR_n$ by
			\begin{equation*}
				\CR_n:=\lbrace Y_{\sigma,e} \mid (\sigma, e)\in \CS_n \times \BZ^n, \sigma^2=1, e_i=e_{\sigma(i)}, \text{ }\forall i\rbrace,
			\end{equation*}
			where $\CS_n$ is the symmetric group of degree $n$ and $Y_{\sigma,e}:=\sigma \diag(\pi^{e_1},\pi^{e_2},\dots,\pi^{e_n})$.
			
			By \cite[Theorem 1]{Hi}, the set $\CR_n$ forms a complete set of representatives of $\Gamma_n \backslash X_n$.
			
			\item For $Y=Y_{\sigma,e} \in \CR_n$ and $k \in \BZ$, we define the following sets.
			\begin{itemize}
				\item $I^{?}_k(Y):=\lbrace i  \mid e_i?k \rbrace$, where $? \in \lbrace >,<,=,\leq, \geq \rbrace$.
				\item $A_1^{h}(Y)=A_1^{h}(\sigma):=\lbrace 1 \leq j \leq n-h \mid \sigma(j)=j \rbrace$.
				\item $A_2^{h}(Y)=A_2^{h}(\sigma):=\lbrace 1 \leq j \leq n-h \mid \sigma(j) \neq j, 1 \leq \sigma(j) \leq n-h \rbrace$.
				\item $B_1^{h}(Y)=B_1^{h}(\sigma):=\lbrace 1 \leq j \leq n-h \mid n-h+1 \leq \sigma(j) \leq n \rbrace$.
				\item $B_2^{h}(Y)=B_2^{h}(\sigma):=\lbrace n-h+1 \leq j \leq n \mid 1 \leq \sigma(j) \leq n-h \rbrace$.
				\item $C_1^{h}(Y)=C_1^{h}(\sigma):=\lbrace n-h+1 \leq j \leq n \mid \sigma(j)=j \rbrace$.
				\item $C_2^{h}(Y)=C_2^{h}(\sigma):=\lbrace n-h+1 \leq j \leq n \mid \sigma(j) \neq j, n-h+1 \leq \sigma(j) \leq n \rbrace$.
				\item $A^{h}(Y)=A^h(\sigma):=A_1^{h}(Y) \cup A_2^{h}(Y)$
				\item $B^{h}(Y)=A^h(\sigma):=B_1^{h}(Y) \cup B_2^{h}(Y)$
				\item $C^{h}(Y)=A^h(\sigma):=C_1^{h}(Y) \cup C_2^{h}(Y)$
			\end{itemize}
			
			\item For $Y=Y_{\sigma,e} \in \CR_n$, we define the following numbers.
			\begin{itemize}
				\item $d_1(Y)=d_1(\sigma)=\vert \lbrace i \mid i=\sigma(i) \rbrace \vert $.
				\item $d_2(Y)=d_2(\sigma)=\frac{1}{2}\vert \lbrace i \mid i \neq \sigma(i) \rbrace \vert $.
				\item For $k \in \BZ$ and $I \subset \lbrace 1,2,\dots,n \rbrace$, $t(\sigma,I):=\vert \lbrace (i,j) \in I\times I \mid i < j <\sigma(i), \sigma(j) < \sigma(i) \rbrace\vert.$
				\item $t(Y):=\sum_{k \in \BZ}\vert \lbrace (i,j) \in I_k^{=}(Y)\times I_k^{=}(Y) \mid i < j <\sigma(i), \sigma(j) < \sigma(i) \rbrace\vert.$
			\end{itemize}
			
			\item For $Y=Y_{\sigma,e} \in \CR_n$, we define the following elements.
			\begin{itemize}
				\item For $k \in \BZ$ and $? \in \lbrace >,<,=,\geq, \leq \rbrace$, we define $Y^?_k$ by
				\begin{equation*}
					Y^{?}_k=Y_{\overline{\sigma},\overline{e}} \in \CR_{\overline{n}},
				\end{equation*}
				where $\overline{n}=\vert I^?_{k}(Y)\vert $, $\overline{\sigma}=\sigma\vert_{I^?_{k}(Y)} \in \CS_{\overline{n}}$, and $\overline{e} \in \BZ^{\overline{n}}$ such that $\overline{e}_i=(e_i)_{i \in I^?_{k}(Y)}$.
			\end{itemize}
			
		\end{enumerate}
	\end{definition}
	
	\begin{proposition}\label{proposition2.6}
		For $Y=Y_{\sigma,e} \in \CR_n$ and $\alpha \in \CE_n$, we have
		\begin{equation*}\begin{array}{ll}
				\CF_{\alpha}(Y,r)=q^{-n^2+n\delta_{\alpha}+\sum_{j=1}^n\alpha_j}f_{n+2r}(Y)&\prod_{j=1}^n (-q)^{-\min(0,\alpha_j+\alpha_{\sigma(j)}) (\min(0,e_j+1)-\min(0,e_j))/2}\\
				&\times(-q)^{-\max(0,\alpha_j+\alpha_{\sigma(j)}) (\min(0,e_j)-\min(0,e_j-1))/2},
			\end{array}
		\end{equation*}
		where $f_k(Y)=\prod_{j=1}^{n}(-q)^{k\min(0,e_j)}$.
	\end{proposition}
	
	\begin{proof}
		The proof of this proposition is similar to the proof of \cite[Lemma 3.15]{Cho}. For $X \in M_{n+2r,n}(F)$ and $\tau \in \CS_{n+2r}, a\in \CE_{n+2r}$, where $A_{n,\delta_{\alpha}}^{[r]}=Y_{\tau,a}$, we have
		\begin{equation*}
			\langle Y,A_{n,\delta_{\alpha}}^{[r]}[X] \rangle=\sum_{i=1}^{n}\sum_{k=1}^{n+2r}\pi^{e_i+a_k}X_{ki}X^*_{\tau(k)\sigma(i)}.
		\end{equation*}
		Here, when $n$ is even and $\alpha$ is odd, we have $a=(0,-1,\dots,-1)$, $\tau(1)=1$, $\tau(i)=n+2-i$ for all $2 \leq i \leq n$, $\tau(i)=i$ for all $n+1 \leq i \leq n+2r$. Also, by the definitions of $1_{\alpha}^{[r]}(X)$, $\Lambda_{\alpha_i,n}$, we have that in the integral
		\begin{equation*}
			\CF_{\alpha}(Y,r)\coloneqq \mathlarger{\int}_{M_{n+2r,n}(F)}\psi(\langle Y,A_{n,\delta_{\alpha}}^{[r]}[X] \rangle)1_{\alpha}^{[r]}(X)dX,
		\end{equation*}
		$X_{kj}$ runs over $O_F$ if $1 \leq k \leq (n+\alpha_j+1)/2$ or $n+1 \leq k$, and over $\pi O_F$ if $(n+\alpha_j+3)/2 \leq k \leq n$.
		
		Now, we can follow the proof of \cite[Lemma 3.15]{Cho} to prove the proposition when $n$ is even and $\alpha$ is odd. The other cases are similar.
	\end{proof}
	
	\begin{proposition}[{\cite[Proposition 2.2]{Hi}}]\label{proposition2.7} For $Y=\sigma (\pi^e1_n) \in \CR_n$, we have
		\begin{equation*}
			\alpha(Y,n)=(q+1)^{d_1(\sigma)}(q(q^2-1))^{d_2(\sigma)}q^{2t(Y)+(e-1)n^2}.
		\end{equation*}
		
	\end{proposition}
	\begin{proposition}\label{proposition2.8} For $Y \in \CR_n$ and $k \in \BZ$, we have
		\begin{equation*}
			\alpha(Y,n)=\alpha(Y^{<}_k,n_0)\alpha(Y^{\geq}_k,n_1)q^{2n_1(\sum_{i \in I_k^{<}(Y)}e_i)-2(\sum_{i \in I_k^{<}(Y)} i)+n_0(n_0+1)}.
		\end{equation*}
		
	\end{proposition}
	\begin{proof}
		This follows from \cite[Theorem 2.5]{Hi}.
	\end{proof}

	\subsection{Reduction formulas for weighted representation densities I}\label{subsection2.2}
	In this subsection, we prove some reduction formulas for weighted representation densities. Assume that $B=Y_{\tau, \lambda}$ with $\tau(h)=h$, $\lambda_h=-1$. Also, assume that $\alpha \in \CE_n$ with $\alpha_{h-1}<0<\alpha_h$. We first introduce the following notation.
	
	\begin{definition}\label{definition2.9}
		For $1 \leq h \leq n$, assume that $\alpha=(\alpha_1,\dots,\alpha_n) \in \CE_n$, $\alpha_{h-1} < 0 \leq \alpha_h$, and $0<\alpha_{h+1}$. We define $\overline{\alpha}^{1} \in \CE_{n-1}$ by
		\begin{equation*}
			\overline{\alpha}^{1}:=(\alpha_1+1,\dots,\alpha_{h-1}+1,\alpha_{h+1}-1,\dots,\alpha_n-1).
		\end{equation*}
		Assume further that $1 \leq l \leq h-1 \leq n$, $-\alpha_{h-1} \leq \dots \leq -\alpha_{l}<\alpha_h \leq -\alpha_{l-1}$. Then, we define $\overline{\alpha}^{i}$ by
		\begin{equation*}
			\overline{\alpha}^{i}:=(\alpha_1+1,\dots,\alpha_{h-i}+1,\alpha_{h-i+1}-1,\dots,\alpha_{h-1}-1,\alpha_{h+1}-1,\dots,\alpha_n-1), \tx{ }2 \leq i \leq h-l+1.
		\end{equation*}
	\end{definition}

	\begin{proposition}\label{proposition2.10}
		For $1 \leq h \leq n$, assume that $\alpha=(\alpha_1,\dots,\alpha_n) \in \CE_n$, and $0 \leq \alpha_h$. Let $B=Y_{\tau, \lambda}$, where $\tau(h)=h$ and $\lambda_h=-1$. Let $\overline{B}$ be the $n-1 \times n-1$ matrix obtained by removing the $h$-th row and the $h$-th column from $B$. Then, we have the following.
		
		\begin{enumerate}
			\item If $\alpha_h=0$, then $W_{\alpha}(B,r)=0$.
			
			\item If $0< \alpha_h \leq -\alpha_{h-1}$, then we have
			\begin{equation*}
				W_{\alpha}(B,r)=q^{n\delta_{\alpha}-(n-1)\delta_{\overline{\alpha}^1}+\alpha_h-3n+1}(1-(-q)^{-\alpha_h})W_{\overline{\alpha}^{1}}(\overline{B},r)
			\end{equation*}
			
			\item For $1 \leq l \leq h-1 \leq n$, assume that $-\alpha_{h-1} \leq \dots \leq -\alpha_{l}<\alpha_h \leq -\alpha_{l-1}$. Then, we have
			\begin{equation*}\begin{array}{ll}
					W_{\alpha}(B,r)&=q^{n\delta_{\alpha}-(n-1)\delta_{\overline{\alpha}^1}+\alpha_h-3n+1}
					\\&\times \lbrace ((-q)^{-\alpha_h-\alpha_{h-1}}-(-q)^{-\alpha_h}))W_{\overline{\alpha}^1}(\overline{B},r)\\
					& +\sum_{j=2}^{h-l} ((-q)^{-\alpha_h-\alpha_{h-j}}-(-q)^{-\alpha_h-\alpha_{h-j+1}})W_{\overline{\alpha}^j}(\overline{B},r)\\
					&+(1-(-q)^{-\alpha_h-\alpha_l})W_{\overline{\alpha}^{h-l+1}}(\overline{B},r)\rbrace.
				\end{array}
			\end{equation*}
			Here, we note that the coefficients of $W_{\overline{\alpha}^i}(\overline{B},r)$ are $0$ if $\overline{\alpha}^i \notin \CE_{n-1}$.
		\end{enumerate} 
		
	\end{proposition}
	
	We will prove this proposition at the end of the subsection.\\

	Note that by Proposition \ref{proposition2.3}, we have $W_{\alpha}(B,r)=\mathlarger{\sum}_{Y \in \CR_n}\dfrac{\CG(Y,B)\CF_{\alpha}(Y,r)}{\alpha(Y;n)}.$ Also, by using the formula for $\CG(Y,B)$ in the proof of \cite[Lemma 3.8]{Cho}, we can show that for $Y=Y_{\sigma,e} \in \CR_n$, $\CG(Y,B) \neq 0$ only when $e_h= 0$ and $\sigma(h) \leq h$, or when $e_h \geq 1$. Therefore, we only need to consider the set $\lbrace Y=Y_{\sigma,e}\in \CR_n \mid e_h=0, \sigma(h)\leq h \text{ or } e_h \geq 1 \rbrace$. We further decompose this set as follows.
	
	\begin{definition}\label{definition2.11}Assume that $\widetilde{Y}=Y_{\widetilde{\sigma},\widetilde{e}} \in \lbrace Y=Y_{\sigma,e}\in \CR_n \mid e_h=0, \sigma(h)\leq h \text{ or } e_h \geq 1 \rbrace$ 
		\begin{enumerate}
			\item We define the set $S_h(\widetilde{Y})$ by setting
			\begin{equation*}
				S_h(\widetilde{Y}):=\lbrace Y=Y_{\sigma,e} \mid \begin{array}{c}
					e_h=0, \sigma(h) \leq h \text{ or } e_h \geq 1, \text{ and  } \sigma(i)=\widetilde{\sigma}(i), e_i=\widetilde{e}_i \\
					\forall i \in I_0^{<}(\widetilde{Y}) \cup (I_0^{=}(\widetilde{Y}) \cap (B^{n-h+1}(\widetilde{Y}) \cup C^{n-h+1}(\widetilde{Y})), i \neq h, \sigma(h) \end{array}\rbrace
			\end{equation*}
			
			\item For $\lbrace i_1,i_2,\dots,i_t \rbrace= \lbrace 1,2,\dots,h-1\rbrace \setminus (I_0^<(\widetilde{Y}) \cup (I_0^{=}(\widetilde{Y}) \cap B_1^{n-h+1}(\widetilde{Y})))$, $i_1 \leq i_2 \leq \dots \leq i_t$, we define the following sets:
			\begin{equation*}\begin{array}{l}
					S_{h}(\widetilde{Y})_0:=\lbrace Y=Y_{\sigma,e} \in S_{h}(\widetilde{Y}) \mid e_h \geq 1\rbrace;\\
					S_{h}(\widetilde{Y})_{i_k}:=\lbrace Y=Y_{\sigma,e} \in S_{h}(\widetilde{Y}) \mid e_h=e_{i_k}=0, \sigma(h)=i_k\rbrace, 1 \leq k \leq t;\\
					S_{h}(\widetilde{Y})_h:=\lbrace Y=Y_{\sigma,e} \in S_{h}(\widetilde{Y}) \mid e_h=0, \sigma(h)=h\rbrace.
				\end{array}
			\end{equation*}
			
		\end{enumerate}
	\end{definition}
	
	\begin{proposition}\label{proposition2.12} For $1 \leq n$, we have
		$
		\mathlarger{\sum}_{\substack{Y=Y_{\sigma,e} \in \CR_n,\\ e_i \geq 1 \text{ }\forall 1 \leq i \leq n}} \dfrac{1}{\alpha(Y,n)}=\dfrac{q^n}{(q^2-1)^n}.
		$
	\end{proposition}
	\begin{proof}
		For $0 \leq k \leq n$, let $A_k=\left( \begin{array}{cc} I_{n-k} & \\ & \pi^{-1}I_k \end{array} \right)$, and $B=\pi^{-1}I_n$. By \cite[Proposition 3.1]{Cho3}, we have
		\begin{equation*}
			\alphad(A_k,B)=\mathlarger{\sum}_{Y \in \CR_n}\frac{\CG(Y,B)\CF(Y,A_k)}{\alpha(Y,n)},
		\end{equation*}
		where $\CF(Y,A_k)=\mathlarger{\prod}_{1 \leq i \leq n} (-q)^{(n-k)\min(0,e_i)+k\min(0,e_i-1)}$. Also, we have that 
		\begin{equation}\label{eq2.1}\begin{array}{l}
				\alphad(A_i,B)=0, \tx{ } 0 \leq i \leq n-1;\\
				\alphad(A_n,B)=\alphad(\pi^{-1}I_n,\pi^{-1}I_n)=q^{-n^2}\prod_{i=1}^n(1-(-q)^{-i}).
			\end{array}
		\end{equation}
		Now, consider the $n+1 \times n+1$ Vandermonde matrix $\CA=((-q)^{-(i-1)(j-1)})_{1 \leq i,j \leq n+1}$ and let
		\begin{equation*}
			\left( \begin{array}{l} k_0 \\ k_1 \\ \vdots \\ k_n \end{array}\middle)=\CA^{-1} \middle(\begin{array}{l} 1\\0  \\ \vdots \\ 0 \end{array}\right).
		\end{equation*}
		Then, we have that $\sum_{i=0}^{n} k_i \CF(Y,A_i)=\left\lbrace \begin{array}{ll} 1 & \text{if } e_i \geq 1, \text{ for all }i\\
			0 & \tx{otherwise} \end{array}\right.$.
		
		Therefore, we have that $\sum_{i=0}^{n} k_i \alphad(A_i,B)=\mathlarger{\sum}_{\substack{Y=Y_{\sigma,e} \in \CR_n,\\ e_i \geq 1 \text{ }\forall 1 \leq i \leq n}}\dfrac{\CG(Y,B)}{\alpha(Y,n)}.$ Now, combining \eqref{eq2.1} and the facts that $k_n=\prod_{i=1}^n(1-(-q)^{-i})$, and $\CG(Y,B)=q^{-n(n-1}(1-q^{-2})^n$ for all $Y=Y_{\sigma,e_i} \in \CR_n$ with $e_i \geq 1$ for all $i$ (see \cite[Proposition 4.2]{Hi}), we have the proposition.
	\end{proof}
	
	\begin{proposition}\label{proposition2.13} For $2 \leq k \leq t$, we have
		$
		\mathlarger{\sum}_{Y \in S_h(\widetilde{Y})_{i_{k-1}}} \dfrac{1}{\alpha(Y,n)}=\dfrac{1}{q^2}\mathlarger{\sum}_{Y \in S_h(\widetilde{Y})_{i_{k}}} \dfrac{1}{\alpha(Y,n)}.
		$
	\end{proposition}
	\begin{proof}
		We consider a bijection $\CT_{\widetilde{Y},k}:S_h(\widetilde{Y})_{i_k} \rightarrow S_h(\widetilde{Y})_{i_k-1}$ sending $Y=Y_{\sigma,e} \in S_h(\widetilde{Y})_{i_k}$ to the following element $Y' \in S_h(\widetilde{Y})_{i_{k-1}}$:
		
		\begin{enumerate}
			\item If $\sigma(i_{k-1})=i_l \neq i_{k-1}$, we define $Y'=Y_{\sigma', e'}$, where
			\begin{equation*}\begin{array}{l}
					\sigma'(h)=i_{k-1}, e_h'=e_{i_{k-1}}'=0,\tx{ }
					\sigma'(i_k)=i_l, e_{i_k}'=e_{i_l}'=e_{i_l}=e_{i_{k-1}},;\\
					\sigma(i)=\sigma'(i), e_i=e_i' \text{ for all } i \notin \lbrace i_l,i_{k-1},i_k,h\rbrace.
				\end{array}
			\end{equation*}
			
			\item If $\sigma(i_{k-1})=i_{k-1}$, we define $Y'=Y_{\sigma', e'}$, where
			\begin{equation*}\begin{array}{l}
					\sigma'(h)=i_{k-1}, e_h'=e_{i_{k-1}}'=0,\tx{ }
					\sigma'(i_k)=i_k, e_{i_k}=e_{i_{k-1}},;\\
					\sigma(i)=\sigma'(i), e_i=e_i' \text{ for all } i \notin \lbrace i_{k-1},i_k,h\rbrace.
				\end{array}
			\end{equation*}
		\end{enumerate}
		
		Now, note that $Y_0^<$ and $I_0^<(Y)$ are fixed by the definitions of $S_h(\widetilde{Y})_{i_k}$ and $S_h(\widetilde{Y})_{i_{k-1}}$.
		
		Also, by Proposition \ref{proposition2.8}, we have
		\begin{equation*}
			\alpha(Y,n)=\alpha(Y^{<}_0,n_0)\alpha(Y^{\geq}_0,n_1)q^{2n_1(\sum_{i \in I_0^{<}(Y)}e_i)-2(\sum_{i \in I_0^{<}} i)+n_0(n_0+1)},
		\end{equation*}
		where $n_0=\vert I_0^<(Y) \vert$, $n_1=n-n_0$. Therefore, we have that the proposition is equivalent to
		\begin{equation*}
			\mathlarger{\sum}_{Y \in S_h(\widetilde{Y})_{i_{k-1}}} \dfrac{1}{\alpha(Y_0^{\geq},n)}=\dfrac{1}{q^2}\mathlarger{\sum}_{Y \in S_h(\widetilde{Y})_{i_{k}}} \dfrac{1}{\alpha(Y_0^{\geq},n)}.
		\end{equation*}
		This shows that it suffices to assume that $Y=Y_0^{\geq}$, i.e., $I_0^{<}(\widetilde{Y})=\emptyset$.
		
		By Proposition \ref{proposition2.8} again, we have that
		\begin{equation*}
			\alpha(Y,n)=\alpha(Y_0^{=},n_0)\alpha(Y_1^{\geq},n_1)q^{-2(\sum_{i \in I_0^{=}(Y)} i)+n_0(n_0+1)},
		\end{equation*}
		where $n_0=\vert I_0^=(Y) \vert$, $n_1=n-n_0$.
		
		Now, we consider the set $S_h(\widetilde{Y})_{i_k}(Y_0^=,I_0^=(Y))$ of all $\CY=Y_{\sigma,e} \in \CR_n$ such that $e_i \geq 0$ for all $i$, $\CY_0^==Y_0^{=}$, and $I_0^=(\CY)=I_0^=(Y)$. In other words, for a fixed $Y_0^=, I_0^=(Y)$, we vary $\CY_1^{\geq}$ over $\CR_{n_1}$.
		
		Then, by Proposition \ref{proposition2.12}, we have
		\begin{equation*}
			\mathlarger{\sum}_{\CY \in S_h(\widetilde{Y})_{i_k}(Y_0^=,I_0^=(Y))} \dfrac{1}{\alpha(\CY,n)}=\dfrac{q^{2(\sum_{i \in I_0^{=}(Y)} i)-n_0(n_0+1)}}{\alpha(Y_0^=,n_0)}\times \dfrac{q^{n_1}}{(q^2-1)^{n_1}}.
		\end{equation*}
		
		Similarly, for $Y'=\CT_{\widetilde{Y},k}(Y)$, $Y_0'^{=}$, and $I_0^{=}(Y')$, we have
		\begin{equation*}
			\mathlarger{\sum}_{\CY' \in \CT_{\widetilde{Y},k}(Y)( S_h(\widetilde{Y})_{i_k}(Y_0^=,I_0^=(Y)))} \dfrac{1}{\alpha(\CY',n)}=\dfrac{q^{2(\sum_{i \in I_0^{=}(Y')} i)-n_0(n_0+1)}}{\alpha(Y_0'^=,n_0)}\times \dfrac{q^{n_1}}{(q^2-1)^{n_1}}.
		\end{equation*}
		Therefore, it suffices to show that
		\begin{equation}\label{eq2.2}
			\dfrac{q^{2(\sum_{i \in I_0^{=}(Y)} i}}{\alpha(Y_0^=,n_0)}=q^2 \times \dfrac{q^{2(\sum_{i \in I_0^{=}(Y')} i)}}{\alpha(Y_0'^=,n_0)}.
		\end{equation}
		
		By Proposition \ref{proposition2.7}, we have
		\begin{equation*}
			\begin{array}{l}
				\alpha(Y_0^=,n_0)=(q+1)^{d_1(\overline{\sigma})}((q(q^2-1))^{d_2(\overline{\sigma})}q^{2t(\sigma,I_0^{=}(Y))-n_0^2},\\
				\alpha(Y_0'^=,n_0)=(q+1)^{d_1(\overline{\sigma'})}((q(q^2-1))^{d_2(\overline{\sigma'})}q^{2t(\sigma',I_0^{=}(Y'))-n_0^2},
			\end{array}
		\end{equation*}
		where $\overline{\sigma}=\sigma\vert_{I_0^{=}(Y)}$ and  $\overline{\sigma'}=\sigma'\vert_{I_0^{=}(Y')}$. Note that $d_1(\overline{\sigma})=d_1(\overline{\sigma'})$ and $d_2(\overline{\sigma})=d_2(\overline{\sigma'})$. Now, we have the following two cases.
		\begin{enumerate}
			\item (Case 1) If $e_{i_{k-1}}=0$, then $I_0^{=}(Y)=I_0^{=}(Y')$. Also, we have $t(\sigma,I_0^=(Y))-t(\sigma',I_0^=(Y'))=-1$, and hence $\alpha(Y_0^=,n_0)=q^{-2}\alpha(Y_0'^=,n_0)$. This implies that \eqref{eq2.2} holds.
			
			\item (Case 2) If $e_{i_{k-1}} \geq 1$, then $I_0^=(Y) \setminus \lbrace i_k \rbrace=I_0^=(Y') \setminus \lbrace i_{k-1} \rbrace$. Also we have $t(\sigma,I_0^=(Y))-t(\sigma',I_0^=(Y'))=i_k-i_{k-1}-1$. These imply that \eqref{eq2.2} holds.
		\end{enumerate}
		This finishes the proof of the proposition.
	\end{proof}
	
	\begin{proposition}\label{proposition2.14} We have
		$
		\mathlarger{\sum}_{Y \in S_h(\widetilde{Y})_{0}} \dfrac{1}{\alpha(Y,n)}=\dfrac{1}{q^2-1}\mathlarger{\sum}_{Y \in S_h(\widetilde{Y})_{i_{1}}} \dfrac{1}{\alpha(Y,n)}.
		$
	\end{proposition}
	\begin{proof}
		As in the proof of Proposition \ref{proposition2.13}, it suffices to assume that $\widetilde{Y}=\widetilde{Y}_0^{\geq}$. Also, as in that proof, for a fixed $Y=Y_{\sigma,e} \in S_h(\widetilde{Y})_{i_{1}}$, we consider the set $S_h(\widetilde{Y})_{i_{1}}(Y_0^=,I_0^=(Y))$ of all $\CY=Y_{\sigma',e'} \in \CR_n$ such that $e'_i \geq 0$ for all $i$, $\CY_0^==Y_0^{=}$, and $I_0^=(\CY)=I_0^=(Y)$. We assume that
		\begin{equation*}
			\begin{array}{l}
				\lbrace i_2, \dots ,i_t \rbrace=\lbrace l_1,\dots,l_u,j_1,\dots,j_v \rbrace,\tx{ } l_1 \leq \dots \leq l_u, j_1 \leq \dots \leq j_v,\\
				e_{l_k}=0,\tx{ }1 \leq k \leq u, \tx{ } e_{j_k} \geq 1, \tx{ }1 \leq k \leq v.
			\end{array}
		\end{equation*}
		
		We now consider the following subsets of $S_h(\widetilde{Y})_{0}$:
		
		$
		S_h(\widetilde{Y})_{0}(Y_0^=,I_0^=(Y))_h=\lbrace \CY'=Y_{\sigma',e'}\in \CR_n \mid \begin{array}{l}e_h'\geq1, \sigma'(i_1)=i_1, e_{i_1}'=0,  \tx{ }e_j' \geq 1 \tx{ for }j \notin I_0^{=}(Y)\\ e_i'=e_i, \sigma'(i)=\sigma(i) \tx{ for }i \in I_0^{=}(Y)\setminus \lbrace i_1,h \rbrace\end{array} \rbrace,$
		
		$	S_h(\widetilde{Y})_{0}(Y_0^=,I_0^=(Y))_0=\lbrace \CY'=Y_{\sigma',e'}\in \CR_n \mid \begin{array}{l}e_h'\geq1, e_{i_1}'\geq 1,  \tx{ }e_j' \geq 1 \tx{ for }j \notin I_0^{=}(Y)\\ e_i'=e_i, \sigma'(i)=\sigma(i) \tx{ for }i \in I_0^{=}(Y)\setminus \lbrace i_1,h \rbrace\end{array} \rbrace,$

		$		S_h(\widetilde{Y})_{0}(Y_0^=,I_0^=(Y))_{j_k}=\lbrace \CY'=Y_{\sigma',e'}\in \CR_n \mid \begin{array}{l}e_h'\geq1, \sigma(i_1)=j_k, e_{i_1}=e_{j_k}=0,\\
			e_j' \geq 1 \tx{ for }j \notin I_0^{=}(Y) \cup \lbrace j_k \rbrace, \\ e_i'=e_i, \sigma'(i)=\sigma(i) \tx{ for }i \in I_0^{=}(Y)\setminus \lbrace i_1,h \rbrace\end{array} \rbrace.
		$
		
		First, we claim that \begin{equation}\label{eq2.3}\mathlarger{\sum}_{\CY' \in 	S_h(\widetilde{Y})_{0}(Y_0^=,I_0^=(Y))_h} \dfrac{1}{\alpha(\CY',n)}=\dfrac{1}{q^{2v+1}(q+1)}\mathlarger{\sum}_{\CY \in S_h(\widetilde{Y})_{i_{1}}(Y_0^=,I_0^=(Y))} \dfrac{1}{\alpha(\CY,n)}.\end{equation}
		
		By Proposition \ref{proposition2.8}, for $\CY \in S_h(\widetilde{Y})_{i_{1}}(Y_0^=,I_0^=(Y))$, we have
		\begin{equation*}\begin{array}{ll}
				\alpha(\CY,n)&=\alpha(\CY_0^{=},n_0)\alpha(\CY_1^{\geq},n_1)q^{-2(\sum_{i \in I_0^{=}(\CY)} i)+n_0(n_0+1)}\\
				&=\alpha(Y_0^{=},n_0)q^{-2(\sum_{i \in I_0^{=}(Y)} i)+n_0(n_0+1)}\alpha(\CY_1^{\geq},n_1)
			\end{array}
		\end{equation*}
		where $n_0=\vert I_0^=(Y) \vert$, $n_1=n-n_0$.
		
		Therefore, by Proposition \ref{proposition2.12}, we have \begin{equation}\label{eq2.4}
			\mathlarger{\sum}_{\CY \in S_h(\widetilde{Y})_{i_{1}}(Y_0^=,I_0^=(Y))} \dfrac{1}{\alpha(\CY,n)}=\dfrac{q^{2(\sum_{i \in I_0^{=}(Y)} i)-n_0(n_0+1)}}{\alpha(Y_0^{=},n_0)}\dfrac{q^{n_1}}{(q^2-1)^{n_1}}
		\end{equation}
		
		Similarly, for any $Y' \in S_h(\widetilde{Y})_{0}(Y_0^=,I_0^=(Y))_h$, we have
		\begin{equation}\label{eq2.5}
			\mathlarger{\sum}_{\CY' \in 	S_h(\widetilde{Y})_{0}(Y_0^=,I_0^=(Y))_h} \dfrac{1}{\alpha(\CY',n)}=\dfrac{q^{2(\sum_{i \in I_0^{=}(Y) \setminus \lbrace h\rbrace} i)-(n_0-1)n_0}}{\alpha(Y_0'^{=},n_0-1)}\dfrac{q^{n_1+1}}{(q^2-1)^{n_1+1}}.
		\end{equation}
		Now, for $Y_0^==Y_{\overline{\sigma},\overline{e}} \in \CR_{n_0}$ and $Y_0'^==Y_{\overline{\sigma'},\overline{e'}} \in \CR_{n_0-1}$, we have that
		\begin{equation*}
			d_1(\overline{\sigma})=d_1({\overline{\sigma'}})-1, d_2({\overline{\sigma}})=d_2({\overline{\sigma'}})+1, t(Y_0^=)-t(Y_0'^=)=h-2-v.
		\end{equation*}
		Combining this with Proposition \ref{proposition2.7}, \eqref{eq2.4}, \eqref{eq2.5}, we have that \eqref{eq2.3} holds.

		Next, we claim that \begin{equation}\label{eq2.6}\mathlarger{\sum}_{\CY' \in 	S_h(\widetilde{Y})_{0}(Y_0^=,I_0^=(Y))_{j_k}} \dfrac{1}{\alpha(\CY',n)}=\dfrac{1}{q^{2(v-k+1)}}\mathlarger{\sum}_{\CY \in S_h(\widetilde{Y})_{i_{1}}(Y_0^=,I_0^=(Y))} \dfrac{1}{\alpha(\CY,n)}.\end{equation}
		
		As above, for any $Y' \in S_h(\widetilde{Y})_{0}(Y_0^=,I_0^=(Y))_{j_k}$, we have
		\begin{equation}\label{eq2.7}
			\mathlarger{\sum}_{\CY' \in 	S_h(\widetilde{Y})_{0}(Y_0^=,I_0^=(Y))_{j_k}} \dfrac{1}{\alpha(\CY',n)}=\dfrac{q^{2(\sum_{i \in I_0^{=}(Y)\cup\lbrace j_k \rbrace \setminus \lbrace h\rbrace} i)-n_0(n_0+1)}}{\alpha(Y_0'^{=},n_0)}\dfrac{q^{n_1}}{(q^2-1)^{n_1}}.
		\end{equation}
		Now, for $Y_0^==Y_{\overline{\sigma},\overline{e}} \in \CR_{n_0}$ and $Y_0'^==Y_{\overline{\sigma'},\overline{e'}} \in \CR_{n_0}$, we have that
		\begin{equation*}
			d_1(\overline{\sigma})=d_1({\overline{\sigma'}}), d_2({\overline{\sigma}})=d_2({\overline{\sigma'}}), t(Y_0^=)-t(Y_0'^=)=h-v-j_k+k-1.
		\end{equation*}
		Combining this with Proposition \ref{proposition2.7}, \eqref{eq2.4}, \eqref{eq2.7}, we have that \eqref{eq2.6} holds.
		
		Finally, we claim that
		\begin{equation}\label{eq2.8}\mathlarger{\sum}_{\CY' \in 	S_h(\widetilde{Y})_{0}(Y_0^=,I_0^=(Y))_{0}} \dfrac{1}{\alpha(\CY',n)}=\dfrac{1}{q^{2v+1}(q^2-1)}\mathlarger{\sum}_{\CY \in S_h(\widetilde{Y})_{i_{1}}(Y_0^=,I_0^=(Y))} \dfrac{1}{\alpha(\CY,n)}.\end{equation}
		
		As above, for any $Y' \in S_h(\widetilde{Y})_{0}(Y_0^=,I_0^=(Y))_{0}$, we have
		\begin{equation}\label{eq2.9}
			\mathlarger{\sum}_{\CY' \in 	S_h(\widetilde{Y})_{0}(Y_0^=,I_0^=(Y))_{0}} \dfrac{1}{\alpha(\CY',n)}=\dfrac{q^{2(\sum_{i \in I_0^{=}(Y) \setminus \lbrace i_1, h\rbrace} i)-(n_0-2)(n_0-1)}}{\alpha(Y_0'^{=},n_0-2)}\dfrac{q^{n_1+2}}{(q^2-1)^{n_1+2}}.
		\end{equation}
		Now, for $Y_0^==Y_{\overline{\sigma},\overline{e}} \in \CR_{n_0}$ and $Y_0'^==Y_{\overline{\sigma'},\overline{e'}} \in \CR_{n_0-2}$, we have that
		\begin{equation*}
			d_1(\overline{\sigma})=d_1({\overline{\sigma'}}), d_2({\overline{\sigma}})=d_2({\overline{\sigma'}})+1, t(Y_0^=)-t(Y_0'^=)=h-v+i_1-3.
		\end{equation*}
		Combining this with Proposition \ref{proposition2.7}, \eqref{eq2.4}, \eqref{eq2.9}, we have that \eqref{eq2.8} holds.
		
		Now, combining \eqref{eq2.3}, \eqref{eq2.6}, and \eqref{eq2.8}, we have the proof of the proposition.
	\end{proof}

	\begin{proposition}\label{proposition2.15}
		We have
		$
		\mathlarger{\sum}_{Y \in S_h(\widetilde{Y})_{h}} \dfrac{1}{\alpha(Y,n)}=\dfrac{q^2}{q+1}\mathlarger{\sum}_{Y \in S_h(\widetilde{Y})_{i_{t}}} \dfrac{1}{\alpha(Y,n)}.
		$
	\end{proposition}
	\begin{proof} As in the proof of Proposition \ref{proposition2.13}, it suffices to assume that $\widetilde{Y}=\widetilde{Y}_0^{\geq}$. Also, by Proposition \ref{proposition2.12} and Proposition \ref{proposition2.13}, it suffices to show that
		\begin{equation*}
			\mathlarger{\sum}_{Y \in S_h(\widetilde{Y})_{h}} \dfrac{1}{\alpha(Y,n)}=q^{2t}(q-1)\mathlarger{\sum}_{Y \in S_h(\widetilde{Y})_{0}} \dfrac{1}{\alpha(Y,n)}
		\end{equation*}
		
		For a fixed $Y=Y_{\sigma,e} \in S_h(\widetilde{Y})_{h}$, we consider the set $S_h(\widetilde{Y})_{h}(Y_0^=,I_0^=(Y))$ of all $\CY=Y_{\sigma',e'} \in \CR_n$ such that $e'_i \geq 0$ for all $i$, $\CY_0^==Y_0^{=}$, and $I_0^=(\CY)=I_0^=(Y)$. Also, we consider the set $S_h(\widetilde{Y})_{0}(Y_0^=,I_0^=(Y))$ of all $\CY'=Y_{\sigma',e'} \in \CR_n$ such that $e'_i \geq 0$ for all $i$, $I_0^=(\CY')=I_0^=(Y)\setminus \lbrace h\rbrace$, and $\CY_0'^=$ is the matrix obtained by removing the $h$-th row and the $h$-th column from $Y_0^{=}$.
		
		As in the proof of Proposition \ref{proposition2.14}, for any $Y' \in S_h(\widetilde{Y})_{0}(Y_0^=,I_0^=(Y))$, $n_0=\vert I_0^=(Y) \vert$, $n_1=n-n_0$, we have
		\begin{equation*}
			\mathlarger{\sum}_{\CY \in 	S_h(\widetilde{Y})_{h}(Y_0^=,I_0^=(Y))} \dfrac{1}{\alpha(\CY,n)}=\dfrac{q^{2(\sum_{i \in I_0^{=}(Y)} i)-n_0(n_0+1)}}{\alpha(Y_0'^{=},n_0)}\dfrac{q^{n_1}}{(q^2-1)^{n_1}},
		\end{equation*}
		\begin{equation*}
			\mathlarger{\sum}_{\CY' \in 	S_h(\widetilde{Y})_{0}(Y_0^=,I_0^=(Y))} \dfrac{1}{\alpha(\CY',n)}=\dfrac{q^{2(\sum_{i \in I_0^{=}(Y)\setminus \lbrace h\rbrace} i)-(n_0-1)n_0}}{\alpha(Y_0'^{=},n_0-1)}\dfrac{q^{n_1+1}}{(q^2-1)^{n_1+1}}.
		\end{equation*}
		Now, for $Y_0^==Y_{\overline{\sigma},\overline{e}} \in \CR_{n_0}$ and $Y_0'^==Y_{\overline{\sigma'},\overline{e'}} \in \CR_{n_0-2}$, we have that
		\begin{equation*}
			d_1(\overline{\sigma})=d_1({\overline{\sigma'}})+1, d_2({\overline{\sigma}})=d_2({\overline{\sigma'}}), t(Y_0^=)-t(Y_0'^=)=h-1-t.
		\end{equation*}
		Combining this with Proposition \ref{proposition2.7}, we have the proof of the proposition.
		
	\end{proof}
	
	\begin{proposition}\label{proposition2.16}
		Let $B=Y_{\tau, \lambda}$, $\tau(h)=h$, $\lambda_h=-1$, $\lambda_1,\dots \lambda_{h-1} \geq 0$, $\lambda_{h+1},\dots,\lambda_n \geq -1$. Then,
		\begin{enumerate}
			\item For $Y, Y' \in \cup_{k=1}^{t} S_h(\widetilde{Y})_{i_k} \cup S_h(\widetilde{Y})_0$, we have $\CG(Y,B)=\CG(Y',B)$.
			
			\item For $Y \in S_h(\widetilde{Y})_h$, $Y' \in \cup_{k=1}^{t} S_h(\widetilde{Y})_{i_k} \cup S_h(\widetilde{Y})_0$, we have $\CG(Y,B)=-(q-1)^{-1}\CG(Y',B)$.
		\end{enumerate}
	\end{proposition}
	\begin{proof}
		By the formula for $\CG(Y,B)$ in the proof of \cite[Proposition 4.2]{Hi} or \cite[Lemma 3.8]{Cho}, we have that for $Y=Y_{\sigma,e}$,
		\begin{equation*}
			\begin{array}{ll}
				\CG(Y,B)&=(1-q^{-2})^nq^{-n(n-1)}\mathlarger{\prod}_{k<j,\tau(k)<\sigma(j)}(-q)^{\min(0,e_j+\lambda_k)}\mathlarger{\prod}_{\substack{k>j,\tau(k)<\sigma(j) \text{ or }\\k<j, \tau(k)>\sigma(j)}}(-q)^{\min(0,e_j+\lambda_k+1)}\\
				& \mathlarger{\prod}_{k>j,\tau(k)>\sigma(j)}(-q)^{\min(0,e_j+\lambda_k+2)}  \mathlarger{\prod}_{k=j=\tau(k)=\sigma(j)}I_1^*(e_j+\lambda_k)  \mathlarger{\prod}_{\substack{k=j,\tau(k)<\sigma(j)\tx{ or }\\k<j, \tau(k)=\sigma(j)}}J_1(e_j+\lambda_k)^{1/2}\\
				& \mathlarger{\prod}_{\substack{k=j,\tau(k)>\sigma(j)\tx{ or }\\k>j, \tau(k)=\sigma(j)}}J_1(e_j+\lambda_k+1)^{1/2} \mathlarger{\prod}_{\substack{k=j,\tau(k)=\sigma(j)\tx{ and }\\k \neq \tau(k) \tx{or } j \neq \sigma(j)}}(\dfrac{q^2J_1(e_j+\lambda_k)-J_1(e_j+\lambda_k+1)}{q^2-1})^{1/2},
			\end{array}
		\end{equation*}
		where $J_1(e)=\left\lbrace \begin{array}{ll} 1 & \tx{ if }e \geq 0\\ 0 &\tx{ if } e \leq -1\end{array}\right.$, $I^*_1(e)=\left\lbrace \begin{array}{ll} 1 & \tx{ if }e \geq 0\\ -(q-1)^{-1} &\tx{ if } e=-1 \\ 0& \tx{ if }e \leq -2\end{array}\right.$
		
		For $Y=Y_{\sigma,e}, Y'=Y_{\sigma',e'} \in \cup_{k=1}^{t} S_h(\widetilde{Y})_{i_k} \cup S_h(\widetilde{Y})_0$, we have that $e_i=e_i'$, $\sigma(i)=\sigma'(i)$ for all $i \in \lbrace i_1,\dots,i_t,h\rbrace$. Therefore, we only need to compare the terms containing $e_i, e_i'$, $i \in \lbrace i_1,\dots,i_t,h\rbrace$. By definition of the sets $S_h(\widetilde{Y})_?$, we have that $e_i \geq 0$ for all $i \in \lbrace i_1,\dots,i_t,h\rbrace$. Also, we have $\sigma(h)\leq h-1$ or $e_h \geq 1$, and $\sigma(h) \leq h-1$ or $e_h' \geq 1$. Furthermore, $\lambda_1,\dots, \lambda_{h-1} \geq 0$, $\lambda_h, \dots \lambda_{n} \geq -1$. Then, we can see that the terms containing $e_i$ (resp. $e_i'$), $i \in \lbrace i_1,\dots,i_t,h \rbrace$ in $\CG(Y,B)$ (resp. $\CG(Y',B)$) are all $1$. This finishes the proof of (1).
		
		Similarly, $Y=Y_{\sigma,e} \in S_h(\widetilde{Y})_h$ and $Y' \in \cup_{k=1}^{t} S_h(\widetilde{Y})_{i_k} \cup S_h(\widetilde{Y})_0$, the only term in $\CG(Y,B)$ that is not equal to $1$ and that contains $e_i$, $i \in \lbrace i_1,\dots,i_t,h \rbrace$ is $I_1^*(e_h+\lambda_h)=I_1^*(-1)=-(q-1)^{-1}$. This shows that (2) holds.
	\end{proof}

	\begin{proposition}\label{proposition2.17}
		For $1 \leq h \leq n$, assume that $\alpha=(\alpha_1,\dots,\alpha_n) \in \CE_n$, $\alpha_{h-1} < 0 \leq \alpha_h$, and $0<\alpha_{h+1}$. Also, Let $B=Y_{\tau, \lambda}$ with $\tau(h)=h$, $\lambda_h=-1$,  $\lambda_1,\dots \lambda_{h-1} \geq 0$, $\lambda_{h+1},\dots,\lambda_n \geq -1$. Let $Y=Y_{\sigma,e}$ with $\sigma(h)=h$, $e_h=0$. Let $\overline{B}$ (resp. $\overline{Y}=Y_{\overline{\sigma},\overline{e}}$) be the $n-1 \times n-1$ matrix obtained by removing the $h$-th row and the $h$-th column from $B$ (resp. $Y$). Then, we have
		\begin{equation*}
			\dfrac{\CF_{\alpha}(Y,r)\CG(Y,B)}{\alpha(Y,n)}=q^{n\delta_{\alpha}-(n-1)\delta_{\overline{\alpha}^1}+\alpha_h-3n+1}(-(-q)^{-\alpha_h})\dfrac{\CF_{\overline{\alpha}^1}(\overline{Y},r)\CG(\overline{Y},\overline{B})}{\alpha(\overline{Y},n-1)}
		\end{equation*}
		
	\end{proposition}
	\begin{proof}
		First, we denote by $n_0:=\vert I_0^<(Y) \vert$, $n_1:=\vert I_0^=(Y) \vert$, 
		
		$k_1:=\vert I_0^<(Y) \cap \lbrace 1, \dots,h-1 \rbrace \vert$, $k_2:=\vert I_0^<(Y) \cap \lbrace h+1, \dots,n \rbrace \vert$, $k_3:=\vert I_0^=(Y) \cap \lbrace h+1, \dots,n \rbrace \vert$.
		
		By the formula for $\CG(Y,B)$ in the proof of Proposition \ref{proposition2.16}, we can see that
		\begin{equation}\label{eq2.10}
			\begin{array}{l}
				\dfrac{\CG(Y,B)}{\CG(\overline{Y},\overline{B})}=(1-q^{-2})q^{-2(n-1)}\times (\tx{terms containing } e_h \tx{ or } \lambda_h)\\
				=(1-q^{-2})q^{-2(n-1)}I_1^*(-1)(-q)^{\sum_{i\in I_0^<(Y)}e_i+\sum_{i \in I_0^{<}(Y) \cap A^{n-h+1}(Y)}1-\sum_{i \in I_0^{\leq}(Y) \cap C^{n-h+1}(Y),  i \neq h}1}\\
				=(1-q^{-2})q^{-2(n-1)}(-(q-1)^{-1})(-q)^{\sum_{i\in I_0^<(Y)}e_i+(k_1-k_2)-\vert I_0^=(Y) \cap C^{n-h+1}(Y)\vert+1}.
			\end{array}
		\end{equation}
		
		By Proposition \ref{proposition2.6}, we have
		\begin{equation}\label{eq2.11}
			\begin{array}{l}
				\dfrac{\CF_{\alpha}(Y,r)}{\CF_{\overline{\alpha}^1}(\overline{Y},r)}\\
				=q^{n\delta_{\alpha}-(n-1)\delta_{\overline{\alpha}^1}+\alpha_h-n-2h+2}
				(-q)^{-\alpha_h+\sum_{i\in I_0^<(Y)}e_i+\vert I_0^{<}(Y) \cap A^{n-h+1}(Y)\vert-
					\vert (I_0^{\leq}(Y) \cap C^{n-h+1}(Y)) \setminus \lbrace h \rbrace\vert} \\
				=q^{n\delta_{\alpha}-(n-1)\delta_{\overline{\alpha}^1}+\alpha_h-n-2h+2}
				(-q)^{-\alpha_h+\sum_{i\in I_0^<(Y)}e_i+(k_1-k_2)-
					\vert (I_0^{=}(Y) \cap C^{n-h+1}(Y))\vert+1} .
			\end{array}
		\end{equation}
		
		Now, by Proposition \ref{proposition2.8} and the fact that $\sum_{i \in I_0^<(\overline{Y})} i=\sum_{i I_0^{<}(Y)}i-k_2$, we have
		\begin{equation}\label{eq2.12}
			\dfrac{\alpha(Y,n)}{\alpha(\overline{Y},n-1)}=\dfrac{\alpha(Y_0^{\geq},n-n_0)}{\alpha(\overline{Y}_0^{\geq},n-n_0-1)}*q^{2\sum_{i\in I_0^{<}(Y)}e_i-2k_2}.
		\end{equation}
		
		Also, by Proposition \ref{proposition2.7} and Proposition \ref{proposition2.8}, we have that for $\sigma_0:=\sigma\vert_{I_0^=(Y)}$,
		\begin{equation}\label{eq2.13}
			\begin{array}{l}
				\alpha(Y_0^{\geq},n-n_0)=\alpha(Y_0^{=},n_1)\alpha(Y_0^{>},n-n_0-n_1)q^{-2\sum_{i \in I_0^{=}(Y^{\geq}_0)} i +n_1(n_1+1)},\\
				\alpha(Y_0^{\geq},n-n_0-1)=\alpha(\overline{Y}_0^{=},n_1-1)\alpha(Y_0^{>},n-n_0-n_1)q^{-2\sum_{i \in I_0^{=}(Y^{\geq}_0)} i +2(h-k_1)+2k_3+(n_1-1)n_1},	\\
				\alpha(Y_0^{=},n_1)=(q+1)^{d_1(\sigma_0)}(q(q^2-1))^{d_2(\sigma_0)}q^{2t(\sigma,I_0^{=}(Y))-n_1^2},\\
				\alpha(\overline{Y}_0^{=},n_1-1)=(q+1)^{d_1(\sigma_0)-1}(q(q^2-1))^{d_2(\sigma_0)}q^{2t(\overline{\sigma},I_0^=(\overline{Y}))-(n_1-1)^2}.
			\end{array}
		\end{equation}
		
		We can compute that $t(\sigma,I_0^{=}(Y))-t(\overline{\sigma},I_0^=(\overline{Y}))=\vert B_1^{n-h+1}(Y) \cap I_0^=(Y)\vert=\vert B_2^{n-h+1}(Y) \cap I_0^=(Y)\vert $. Combining this with \eqref{eq2.12}, \eqref{eq2.13}, and the fact that $\vert (C^{n-h+1}\cup B_2^{n-h+1}(Y)) \cap I_0^=(Y) \vert= k_3+1$, we have
		\begin{equation}\label{eq2.14}
			\dfrac{\alpha(Y,n)}{\alpha(\overline{Y},n-1)}=(q+1)q^{-2h-2(k_2-k_1)-2\vert C^{n-h+1}(Y) \cap I_0^=(Y) \vert+2 \sum_{i \in I_0^{<}(Y)} e_i +3}.
		\end{equation}
		
		Now, combining \eqref{eq2.10}, \eqref{eq2.11}, and \eqref{eq2.14}, we have the proof of the proposition.
	\end{proof}
	
	\begin{definition}\label{definition2.18}
		For $\widetilde{Y} \in \lbrace Y=Y_{\sigma,e} \in \CR_n \mid e_h=0, \sigma(h)\leq h \tx{ or } e_h \geq 1 \rbrace$ and $S_h(\widetilde{Y})_h$, we define $\overline{S}_h(\widetilde{Y}) \subset \CR_{n-1}$ by
		\begin{equation*}
			\overline{S}_h(\widetilde{Y}):=\left\lbrace \overline{Y} \in \CR_{n-1} \middle\vert \overline{Y}=\begin{array}{l} \tx{the matrix obtained by removing} \\ \tx{the h-th row and the h-th column from Y} \in S_h(\widetilde{Y})_h.
			\end{array}\right \rbrace.
		\end{equation*}
	\end{definition}
	\begin{proposition}\label{proposition2.19}
		For $1 \leq h \leq n$, assume that $\alpha=(\alpha_1,\dots,\alpha_n) \in \CE_n$, $0 \leq \alpha_h$. Let $B=Y_{\tau, \lambda}$ with $\tau(h)=h$, $\lambda_h=-1$. Let $\overline{B}$ be the $n-1 \times n-1$ matrix obtained by removing the $h$-th row and the $h$-th column from $B$. The following statements hold.
		\begin{enumerate}
			\item If $\alpha_h=0$, then 
			\begin{equation*}
				\mathlarger{\sum}_{j \in \lbrace 0,i_1,\dots,i_t,h \rbrace}\mathlarger{\sum}_{Y \in S_h(\widetilde{Y})_{j}}	\dfrac{\CF_{\alpha}(Y,r)\CG(Y,B)}{\alpha(Y,n)}=0.
			\end{equation*}
			\item If $\alpha_{h-1}< 0 <\alpha_h$, then we have that
			\begin{equation}\label{eq2.15}\mathlarger{\sum}_{Y \in S_h(\widetilde{Y})_{h}}	\dfrac{\CF_{\alpha}(Y,r)\CG(Y,B)}{\alpha(Y,n)}=q^{n\delta_{\alpha}-(n-1)\delta_{\overline{\alpha}^1}+\alpha_h-3n+1}(-(-q)^{-\alpha_h})\mathlarger{\sum}_{\overline{Y} \in \overline{S}_h(\widetilde{Y})}\dfrac{\CF_{\overline{\alpha}^1}(\overline{Y},r)\CG(\overline{Y},\overline{B})}{\alpha(\overline{Y},n-1)}\end{equation}.
			
			\item Assume that $0<\alpha_h \leq -\alpha_{h-1}$. Then,
			\begin{equation}\label{eq2.16}
				\mathlarger{\sum}_{j \in \lbrace 0,i_1,\dots,i_t \rbrace}\mathlarger{\sum}_{Y \in S_h(\widetilde{Y})_{j}}	\dfrac{\CF_{\alpha}(Y,r)\CG(Y,B)}{\alpha(Y,n)}=q^{n\delta_{\alpha}-(n-1)\delta_{\overline{\alpha}^1}+\alpha_h-3n+1}\mathlarger{\sum}_{\overline{Y} \in \overline{S}_h(\widetilde{Y})}\dfrac{\CF_{\overline{\alpha}^1}(\overline{Y},r)\CG(\overline{Y},\overline{B})}{\alpha(\overline{Y},n-1)}.
			\end{equation}
		\end{enumerate}
		
		\begin{enumerate}[resume]
			\item 	Assume that for $1 \leq l \leq h-1$, $-\alpha_l<\alpha_h \leq -\alpha_{l-1}$. Then,
			\begin{equation*}\begin{array}{ll}
					\mathlarger{\sum}_{j \in \lbrace 0,i_1,\dots,i_t,h \rbrace}\mathlarger{\sum}_{Y \in S_h(\widetilde{Y})_{j}}	\dfrac{\CF_{\alpha}(Y,r)\CG(Y,B)}{\alpha(Y,n)}&=q^{n\delta_{\alpha}-(n-1)\delta_{\overline{\alpha}^1}+\alpha_h-3n+1}\mathlarger{\sum}_{\overline{Y} \in \overline{S}_h(\widetilde{Y})}\dfrac{\CG(\overline{Y},\overline{B})}{\alpha(\overline{Y},n-1)}
					\\ &\times \lbrace \CF_{\overline{\alpha}^{1}}(\overline{Y},r)((-q)^{-\alpha_h-\alpha_{h-1}}-(-q)^{-\alpha_h})\\
					&+\sum_{j=2}^{h-l} \CF_{\overline{\alpha}^{j}}(\overline{Y},r)((-q)^{-\alpha_h-\alpha_{h-j}}-(-q)^{-\alpha_h-\alpha_{h-j+1}})\\
					&+\CF_{\overline{\alpha}^{h-l+1}}(\overline{Y},r)(1-(-q)^{-\alpha_h-\alpha_{l}}) \rbrace.
				\end{array}
			\end{equation*}	
		\end{enumerate}
	\end{proposition}
	\begin{proof} By Proposition \ref{proposition2.6}, we have that $\CF_{\alpha}(Y,r)=\CF_{\alpha}(Y',r)$ for all $Y, Y' \in \cup_{j \in \lbrace 0,i_1,\dots,i_t,h\rbrace } S_{h}(\widetilde{Y})_j$. Therefore, by Proposition \ref{proposition2.13}, \ref{proposition2.14}, \ref{proposition2.15}, and \ref{proposition2.16}, we have that
		\begin{equation*}
			\mathlarger{\sum}_{j \in \lbrace 0,i_1,\dots,i_t,h \rbrace}\mathlarger{\sum}_{Y \in S_h(\widetilde{Y})_{j}}	\dfrac{\CF_{\alpha}(Y,r)\CG(Y,B)}{\alpha(Y,n)}=\mathlarger{\sum}_{Y \in S_h(\widetilde{Y})_{h}}	\dfrac{\CF_{\alpha}(Y,r)\CG(Y,B)}{\alpha(Y,n)}\lbrace 1-\frac{q^2-1}{q^2}-\dots-\frac{q^2-1}{q^{2t}}-\frac{1}{q^{2t}}\rbrace=0.
		\end{equation*}
		This completes the proof of (1).
		
		(2) follows directly from Proposition \ref{proposition2.17}.
		
		By Proposition \ref{proposition2.6}, for $Y \in S_{h}(\widetilde{Y})_h$, $Y' \in S_{h}(\widetilde{Y})_{j}, j \in \lbrace 0,i_1,\dots,i_t \rbrace$, we have
		\begin{equation}\label{eq2.17}
			\CF_{\alpha}(Y',r)=\left\lbrace \begin{array}{ll}(-q)^{\min(-\alpha_i,\alpha_h)}\CF_{\alpha}(Y,r) & \tx{if }i \in \lbrace i_1,\dots,i_t\rbrace,\\
				(-q)^{\alpha_h}\CF_{\alpha}(Y,r) & \tx{if }i=0.
			\end{array}\right.
		\end{equation}
		
		Also, for $\overline{Y} \in \overline{S}_h(\widetilde{Y})$ and $1 \leq k \leq h-l+1$, we have
		\begin{equation}\label{eq2.18}
			\CF_{\overline{\alpha}^{k}}(\overline{Y},r)=(-q)^{-2\vert \lbrace i_j \mid 1 \leq j\leq t, h-k+1 \leq i_j \leq h-1 \rbrace \vert}\CF_{\overline{\alpha}^1}(\overline{Y},r).
		\end{equation}
		
		Combining \eqref{eq2.17}, Proposition \ref{proposition2.13}, \ref{proposition2.14}, \ref{proposition2.15}, and \ref{proposition2.16}, we have
		\begin{equation}\label{eq2.19}
			\mathlarger{\sum}_{Y \in S_h(\widetilde{Y})_{i_k}}	\dfrac{\CF_{\alpha}(Y,r)\CG(Y,B)}{\alpha(Y,n)}=-\dfrac{q^2-1}{q^{2(t-k+1)}}(-q)^{\min(-\alpha_{i_k},\alpha_h)}\mathlarger{\sum}_{Y \in S_h(\widetilde{Y})_{h}}	\dfrac{\CF_{\alpha}(Y,r)\CG(Y,B)}{\alpha(Y,n)}.
		\end{equation}
		
		Now, assume that $\alpha_h \leq -\alpha_{h-1}$. By \eqref{eq2.17}, Proposition \ref{proposition2.13}, \ref{proposition2.14}, \ref{proposition2.15}, and \ref{proposition2.16}, we have
		
		\begin{equation*}	\mathlarger{\sum}_{j \in \lbrace 0,i_1,\dots,i_t \rbrace}\mathlarger{\sum}_{Y \in S_h(\widetilde{Y})_{j}}	\dfrac{\CF_{\alpha}(Y,r)\CG(Y,B)}{\alpha(Y,n)}=-(-q)^{\alpha_h}\mathlarger{\sum}_{Y \in S_h(\widetilde{Y})_{h}}	\dfrac{\CF_{\alpha}(Y,r)\CG(Y,B)}{\alpha(Y,n)}.
		\end{equation*}
		Combining this with \eqref{eq2.15}, we have the proof of \eqref{eq2.16}.

		Assume that for $1 \leq l \leq h-1$, $-\alpha_l<\alpha_h \leq -\alpha_{l-1}$, and $\alpha_{h-1}<0$.
		Assume further that $l \leq  i_k$. 
		
		We claim that
		\begin{equation}\label{eq2.20}\begin{array}{ll}
				\mathlarger{\sum}_{Y \in S_h(\widetilde{Y})_{i_k}}	\dfrac{\CF_{\alpha}(Y,r)\CG(Y,B)}{\alpha(Y,n)}&=q^{n\delta_{\alpha}-(n-1)\delta_{\overline{\alpha}^1}+\alpha_h-3n+1}\mathlarger{\sum}_{\overline{Y} \in \overline{S}_h(\widetilde{Y})}\dfrac{\CG(\overline{Y},\overline{B})}{\alpha(\overline{Y},n-1)}\\
				&\times \lbrace \CF_{\overline{\alpha}^{h-i_{k+1}+1}}(\overline{Y},r)((-q)^{-\alpha_h-\alpha_{i_{k+1}-1}})\\
				&+\sum_{j=2}^{i_{k+1}-i_k} \CF_{\overline{\alpha}^{h-i_{k+1}+j}}(\overline{Y},r)((-q)^{-\alpha_h-\alpha_{i_{k+1}-j}}-(-q)^{-\alpha_h-\alpha_{i_{k+1}-j+1}})\\
				&+\CF_{\overline{\alpha}^{h-i_{k}+1}}(\overline{Y},r)(-(-q)^{-\alpha_h-\alpha_{i_{k}}}) \rbrace.
			\end{array}
		\end{equation}
		
		By \eqref{eq2.18}, we have
		\begin{equation*}\begin{array}{ll}
				\CF_{\overline{\alpha}^{h-i_{k+1}+j}}(\overline{Y},r)=(-q)^{-2(t-k)}\CF_{\overline{\alpha}^1}(\overline{Y},r), & 1 \leq j \leq i_{k+1}-i_k;\\
				\CF_{\overline{\alpha}^{h-i_{k}+1}}(\overline{Y},r)=(-q)^{-2(t-k+1)}\CF_{\overline{\alpha}^1}(\overline{Y},r).&
			\end{array}
		\end{equation*}
		
		Now, the right hand side of \eqref{eq2.20} is
		\begin{equation*}\begin{array}{l}
				q^{n\delta_{\alpha}-(n-1)\delta_{\overline{\alpha}^1}+\alpha_h-3n+1}\mathlarger{\sum}_{\overline{Y} \in \overline{S}_h(\widetilde{Y})}\dfrac{\CG(\overline{Y},\overline{B})}{\alpha(\overline{Y},n-1)}
				\times\CF_{\overline{\alpha}^{1}}(\overline{Y},r) \dfrac{q^2-1}{q^{2(t-k+1)}}(-q)^{-\alpha_{i_k}+\alpha_h}.
			\end{array}
		\end{equation*}
		
		Combining this with \eqref{eq2.15} and \eqref{eq2.19}, we have the proof of \eqref{eq2.20}.
		
		Now, assume that $i_k < l \leq i_{k+1}$.  We claim that
		\begin{equation}\label{eq2.21}\begin{array}{ll}
				\mathlarger{\sum}_{j \in \lbrace 0,i_1,\dots,i_k \rbrace}\mathlarger{\sum}_{Y \in S_h(\widetilde{Y})_{j}}	\dfrac{\CF_{\alpha}(Y,r)\CG(Y,B)}{\alpha(Y,n)}&=q^{n\delta_{\alpha}-(n-1)\delta_{\overline{\alpha}^1}+\alpha_h-3n+1}\mathlarger{\sum}_{\overline{Y} \in \overline{S}_h(\widetilde{Y})}\dfrac{\CG(\overline{Y},\overline{B})}{\alpha(\overline{Y},n-1)}
				\\ &\times \lbrace \CF_{\overline{\alpha}^{h-i_{k+1}+1}}(\overline{Y},r)((-q)^{-\alpha_h-\alpha_{i_{k+1}-1}})\\
				&+\sum_{j=2}^{i_{k+1}-i_k} \CF_{\overline{\alpha}^{h-i_{k+1}+j}}(\overline{Y},r)((-q)^{-\alpha_h-\alpha_{i_{k+1}-j}}-(-q)^{-\alpha_h-\alpha_{i_{k+1}-j+1}})\\
				&+\CF_{\overline{\alpha}^{h-l+1}}(\overline{Y},r)(1-(-q)^{-\alpha_h-\alpha_{l}}) \rbrace.
			\end{array}
		\end{equation}
		
		By \eqref{eq2.17}, \eqref{eq2.18}, Proposition \ref{proposition2.13}, \ref{proposition2.14}, \ref{proposition2.15}, and \ref{proposition2.16}, we have
		\begin{equation}\label{eq2.22}
			\mathlarger{\sum}_{j \in \lbrace 0,i_1,\dots,i_k \rbrace}\mathlarger{\sum}_{Y \in S_h(\widetilde{Y})_{j}}	\dfrac{\CF_{\alpha}(Y,r)\CG(Y,B)}{\alpha(Y,n)}=-\dfrac{1}{q^{2(t-k)}}(-q)^{\alpha_h}\mathlarger{\sum}_{Y \in S_h(\widetilde{Y})_{h}}	\dfrac{\CF_{\alpha}(Y,r)\CG(Y,B)}{\alpha(Y,n)}
		\end{equation}
		On the other hand, by \eqref{eq2.18}, we have
		\begin{equation*}\begin{array}{ll}
				\CF_{\overline{\alpha}^{h-i_{k+1}+j}}(\overline{Y},r)=(-q)^{-2(t-k)}\CF_{\overline{\alpha}^1}(\overline{Y},r), & 1 \leq j \leq i_{k+1}-l+1.
			\end{array}
		\end{equation*}
		Therefore, the right hand side of \eqref{eq2.21} is
		\begin{equation*}\begin{array}{l}
				q^{n\delta_{\alpha}-(n-1)\delta_{\overline{\alpha}^1}+\alpha_h-3n+1}\mathlarger{\sum}_{\overline{Y} \in \overline{S}_h(\widetilde{Y})}\dfrac{\CG(\overline{Y},\overline{B})}{\alpha(\overline{Y},n-1)}
				\times\CF_{\overline{\alpha}^{1}}(\overline{Y},r) \dfrac{1}{q^{2(t-k)}}.
			\end{array}
		\end{equation*}
		Combining this with \eqref{eq2.15} and \eqref{eq2.22}, we have the proof of \eqref{eq2.21}.
		Combining \eqref{eq2.20} and \eqref{eq2.21}, we have the proof of (4) when $\alpha_{h-1}<0$.
		
		Finally, assume that $\alpha_{h-1}=0$. Assume further that $\alpha_{m-1}<0, \alpha_{m}=\dots=\alpha_{h-1}=0$ for some $l \leq m \leq h-1$. In this case, we cannot apply Proposition \ref{proposition2.17}. Therefore, we consider $\alpha'=(\alpha_1,\dots,\alpha_{m-1},-2,-2,\dots,-2,\alpha_h,\dots,\alpha_n) \in \CE_n$ instead of $\alpha$. We note that $\overline{\alpha'}^1=\overline{\alpha}^{h-m+1}$. Assume that $i_k <m \leq i_{k+1}$. By Proposition \ref{proposition2.6}, we have that for $Y \in S_h(\widetilde{Y})_j$, $j \in \lbrace 0, i_1,\dots,i_t,h \rbrace$,
		\begin{equation*}
			\CF_{\alpha}(Y,r)=\CF_{\alpha'}(Y,r)(-q)^{2 \vert \lbrace i_j \mid 1 \leq j\leq t, m \leq i_j \leq h-1 \rbrace \vert}=\CF_{\alpha'}(Y,r)(-q)^{2(t-k)}.
		\end{equation*}
		Now, by Proposition \ref{proposition2.17}, we have
		\begin{equation}\label{eq2.23}\begin{array}{l}\mathlarger{\sum}_{Y \in S_h(\widetilde{Y})_{h}}	\dfrac{\CF_{\alpha}(Y,r)\CG(Y,B)}{\alpha(Y,n)}=\mathlarger{\sum}_{Y \in S_h(\widetilde{Y})_{h}}	\dfrac{\CF_{\alpha'}(Y,r)\CG(Y,B)}{\alpha(Y,n)}(-q)^{2(t-k)}
				\\=q^{n\delta_{\alpha}-(n-1)\delta_{\overline{\alpha}^1}+\alpha_h-3n+1}(-(-q)^{-\alpha_h})(-q)^{2(t-k)}\mathlarger{\sum}_{\overline{Y} \in \overline{S}_h(\widetilde{Y})}\dfrac{\CF_{\overline{\alpha}^{h-m+1}}(\overline{Y},r)\CG(\overline{Y},\overline{B})}{\alpha(\overline{Y},n-1)}\end{array}.\end{equation}
		Combining \eqref{eq2.17}, Proposition \ref{proposition2.13}, \ref{proposition2.14}, \ref{proposition2.15}, and \ref{proposition2.16}, we have
		\begin{equation*}\begin{array}{ll}\mathlarger{\sum}_{j \in \lbrace i_{k+1},\dots,i_t,h \rbrace}\mathlarger{\sum}_{Y \in S_h(\widetilde{Y})_{h}}	\dfrac{\CF_{\alpha}(Y,r)\CG(Y,B)}{\alpha(Y,n)}&=q^{n\delta_{\alpha}-(n-1)\delta_{\overline{\alpha}^1}+\alpha_h-3n+1}(-(-q)^{-\alpha_h})\mathlarger{\sum}_{\overline{Y} \in \overline{S}_h(\widetilde{Y})}\dfrac{\CF_{\overline{\alpha}^{h-m+1}}(\overline{Y},r)\CG(\overline{Y},\overline{B})}{\alpha(\overline{Y},n-1)}.
		\end{array}\end{equation*}
		Now, we can follow the above arguments to prove (4) when $\alpha_{h-1}=0$.
		This finishes the proof of the proposition.
	\end{proof}
	
	\begin{proof}[Proof of Proposition \ref{proposition2.10}]
		This follows from Proposition \ref{proposition2.3} and Proposition \ref{proposition2.19} where $\widetilde{Y}$ runs over the set $\lbrace Y=Y_{\sigma,e}\in \CR_n \mid e_h=0, \sigma(h)\leq h \text{ or } e_h \geq 1 \rbrace$.
	\end{proof}
	\subsection{Reduction formulas for weighted representation densities II}\label{subsection2.3}
	In this subsection, we prove some reduction formulas for weighted representation densities. Assume that $B=Y_{\tau, \lambda}$ with $\tau(1)=1$, $\lambda_1=0$. We first introduce the following notation.
	
	\begin{definition}\label{definition2.20} Let $0 \leq k \leq n$.
		\begin{enumerate}
			\item For $X=\left( \begin{array}{cc} A & B \\C & D \end{array}\right) \in X_n(F)$, where $A \in M_{n-k,n-k}(F), B \in M_{n-k,k}(F), C \in M_{k,n-k}(F), D \in M_{k,k}(F)$, we denote by $X^{\vee_k}$ the element $X^{\vee_k}:=\left( \begin{array}{cc} \pi D & C \\B & \pi^{-1}A \end{array}\right)$. Also, we denote by $X^{\wedge_k}$ the element $X^{\wedge_k}:=\left( \begin{array}{cc} \pi^{-1} D & C \\B & \pi A \end{array}\right)$
			\item For $\alpha=(\alpha_1,\dots,\alpha_n) \in \CE_n$ with $\alpha_{n-k} \leq 0 \leq \alpha_{n-k+1}$, we denote by $\alpha^{\vee_k}$ the element $\alpha^{\vee_k}=(\alpha_{n-k+1}-n,\dots,\alpha_n-n,\alpha_1+n,\dots,\alpha_{n-k}+n)\in \CE_n$.
		\end{enumerate}
	\end{definition}
	
	\begin{definition}\label{definition2.21}
		For $1 \leq h \leq n$, assume that $\alpha=(\alpha_1,\dots,\alpha_n) \in \CE_n$,$-n<\alpha_1$, and $\alpha_{n-h+1}\leq 0<\alpha_{n-h+2}$. We define $^1\overline{\alpha} \in \CE_{n-1}$ by
		\begin{equation*}
			^1\overline{\alpha}:=(\alpha_2,\dots,\alpha_n).
		\end{equation*}
		Assume further that $n-h+2 \leq l \leq n$, $-\alpha_{l}<\alpha_1 \leq -\alpha_{l-1}$. Then, we define $^i\overline{\alpha}$ by
		\begin{equation*}
			^i\overline{\alpha}:=(\alpha_2,\dots,\alpha_{n-i+1},\alpha_{n-i+2}-2,\dots,\alpha_n-2), \tx{ }2 \leq i \leq n-l+2.
		\end{equation*}
	\end{definition}

	\begin{proposition}\label{proposition2.22}
		For $1 \leq k \leq n$ and $\alpha \in \CE_n$ with $\alpha_{n-k} \leq 0 \leq \alpha_{n-k+1}$ and $Y \in \CR_n$, we have
		\begin{equation*}
			\CF_{\alpha}(Y,0)=\CF_{\alpha^{\vee_k}}(Y^{\wedge_k},0)q^{n\delta_{\alpha}-n\delta_{\alpha^{\vee_k}}-n(n-2k)}.
		\end{equation*}
	\end{proposition}
	
	\begin{proposition}[cf. {\cite[Lemma 3.8]{Cho}}]\label{proposition2.23} For $1 \leq k \leq n$, and $Y,B \in \CR_n$, we have
		\begin{equation*}
			\CG(Y,B)=\CG(Y^{\wedge_k},B^{\vee_k}).
		\end{equation*}
		
	\end{proposition}
	
	\begin{proposition}[{\cite[Lemma 3.9]{Cho}}]\label{proposition2.24} For $1 \leq k \leq n$ and $Y \in \CR_n$, we have
		\begin{equation*}
			\alpha(Y^{\wedge_k},n)=q^{n(n-2k)}\alpha(Y,n).
		\end{equation*}
	\end{proposition}

	\begin{proposition}\label{proposition2.25}
		For $1 \leq h \leq n$, assume that $\alpha=(\alpha_1,\dots,\alpha_n) \in \CE_n$ and $\alpha_{n-h+1} \leq 0 \leq\alpha_{n-h+2}$. Let $B=Y_{\tau, \lambda}$ with $\tau(1)=1$, $\lambda_1=0$. Let $\overline{B}$ be the $n-1 \times n-1$ matrix obtained by removing the first row and the first column from $B$. For $r \geq 0$, we further assume that $-n< \alpha_1, \alpha_n <n$ if $r \geq 1$.  Then, we have the following.
		
		\begin{enumerate}
			\item If $\alpha_1 \leq -\alpha_{n}$, then we have
			\begin{equation*}
				W_{\alpha}(B,r)=q^{n\delta_{\alpha}-(n-1)\delta_{^1\overline{\alpha}}+\alpha_1-2n+1}(1-(-q)^{-\alpha_1-n-2r})W_{^1\overline{\alpha}}(\overline{B},r)
			\end{equation*}
			
			\item For $n-h+2 \leq l \leq n$, assume that $-\alpha_{l}<\alpha_1 \leq -\alpha_{l-1}$. Then, we have
			\begin{equation*}\begin{array}{ll}
					W_{\alpha}(B,r)&=q^{n\delta_{\alpha}-(n-1)\delta_{^1\overline{\alpha}}+\alpha_1-2n+1}
					\\&\times \lbrace ((-q)^{-\alpha_1-\alpha_{n}}-(-q)^{-\alpha_1-n-2r}))W_{\tx{}^1\overline{\alpha}}(\overline{B},r)\\
					& +\sum_{j=2}^{n-l+1} ((-q)^{-\alpha_1-\alpha_{n+1-j}}-(-q)^{-\alpha_1-\alpha_{n+2-j}})W_{\tx{}^j\overline{\alpha}}(\overline{B},r)\\
					&+(1-(-q)^{-\alpha_1-\alpha_l})W_{\tx{}^{n-l+2}\overline{\alpha}}(\overline{B},r)\rbrace.
				\end{array}
			\end{equation*}
		\end{enumerate} 
		
	\end{proposition}
	\begin{proof}
		For simplicity, we write $^{\vee}$ for $^{\vee_{h-1}}$. Note that
		\begin{equation*}
			\alpha^{\vee}=(\alpha_{n-h+2}-n,\dots,\alpha_n-n,\alpha_1+n,\dots,\alpha_{n-h+1}+n).
		\end{equation*}
		By Proposition \ref{proposition2.22}, Proposition \ref{proposition2.23}, and Proposition \ref{proposition2.24}, we have
		\begin{equation}\label{eq2.24}
			\begin{array}{ll}
				W_{\alpha}(B,r)&=\mathlarger{\sum}_{Y \in \CR_n}\dfrac{\CF_{\alpha}(Y,0)f_{2r}(Y)\CG(Y,B)}{\alpha(Y,n)}\\
				&=q^{n\delta_{\alpha}-n\delta_{\alpha^{\vee}}}\mathlarger{\sum}_{Y \in \CR_n}\dfrac{\CF_{\alpha^{\vee}}(Y^{\wedge_{h-1}},0)f_{2r}(Y)\CG(Y^{\wedge_{h-1}},B^{\vee})}{\alpha(Y^{\wedge_{h-1}},n)}\\
				&=q^{n\delta_{\alpha}-n\delta_{\alpha^{\vee}}}\mathlarger{\sum}_{Y \in \CR_n}\dfrac{\CF_{\alpha^{\vee}}(Y,0)\CG(Y,B^{\vee})f_{2r}(Y^{\wedge_{n-h+1}})}{\alpha(Y,n)}.
			\end{array}
		\end{equation}
		
		Similarly, we have
		\begin{equation}\label{eq2.25}
			\begin{array}{ll}
				W_{\tx{}^{k}\overline{\alpha}}(\overline{B},r)
				&=q^{(n-1)\delta_{\tx{}^k\overline{\alpha}}-(n-1)\delta_{\overline{\alpha^{\vee}}^{k}}}\mathlarger{\sum}_{\overline{Y} \in \CR_{n-1}}\dfrac{\CF_{\overline{\alpha^{\vee}}^k}(\overline{Y},0)\CG(\overline{Y},\overline{B}^{\vee})f_{2r}(\overline{Y}^{\wedge_{n-h}})}{\alpha(\overline{Y},n-1)}.
			\end{array}
		\end{equation}
		When $r=0$, the proposition follows from Proposition \ref{proposition2.10} together with \eqref{eq2.24} and \eqref{eq2.25}.
		
		Now, assume that $r \geq 1$. By the same arguments as in Section \ref{subsection2.2}, for $Y=Y_{\sigma,e}\in\CR_n$, we know that $\CG(Y,B^{\vee})\neq 0$ only when $e_h=0$ and $\sigma(h)\leq h$, or when $e_h \geq 1$. Therefore, for $\widetilde{Y} \in \lbrace Y=Y_{\sigma,e}\in \CR_n \mid e_h=0, \sigma(h)\leq h \tx{ or }e_h \geq 1 \rbrace$, we consider the sets $S_{h}(\widetilde{Y})_j$, $j\in \lbrace 0,i_1,\dots,i_t,h\rbrace$ in Definition \ref{definition2.11} and $\overline{S}_h(\widetilde{Y})$ in Definition \ref{definition2.18}. Note that for $Y \in S_h(\widetilde{Y})_h, Y' \in S_h(\widetilde{Y})_j$, $j\in \lbrace 0,i_1,\dots,i_t \rbrace$, and $\overline{Y} \in \overline{S}_h(\widetilde{Y})$, we have
		\begin{equation}\label{eq2.26}
			\begin{array}{ll}
				f_{2r}(Y^{\wedge_{n-h+1}})=(-q)^{-2r}f_{2r}(\overline{Y}^{\wedge_{n-h}}),\\
				f_{2r}(Y'^{\wedge_{n-h+1}})=f_{2r}(\overline{Y}^{\wedge_{n-h}}).
			\end{array}
		\end{equation}
		
		By Proposition \ref{proposition2.17} together with \eqref{eq2.24}, \eqref{eq2.25} and \eqref{eq2.26}, we have
		\begin{equation}\label{eq2.27}
			\begin{array}{l}
				\mathlarger{\sum}_{Y \in S_h(\widetilde{Y})_{h}}	\dfrac{\CF_{\alpha^{\vee}}(Y,0)\CG(Y,B^{\vee})}{\alpha(Y,n)}f_{2r}(Y^{\wedge_{n-h+1}})\\
				=q^{n\delta_{\alpha^{\vee}}-(n-1)\delta_{\overline{\alpha^{\vee}}^1}+\alpha_1-2n+1}(-(-q)^{-\alpha_1-n-2r})\mathlarger{\sum}_{\overline{Y} \in \overline{S}_h(\widetilde{Y})}\dfrac{\CF_{\overline{\alpha^{\vee}}^1}(\overline{Y},0)\CG(\overline{Y},\overline{B}^{\vee})}{\alpha(\overline{Y},n-1)}f_{2r}(\overline{Y}^{\wedge_{n-h}}).
			\end{array}
		\end{equation}
		Now, by \eqref{eq2.26}, we know that $f_{2r}(Y'^{\wedge_{n-h+1}})$ takes the same value for all $Y' \in S_h(\widetilde{Y})_j$, $j \in \lbrace 0,i_1,\dots,i_t\rbrace$. Therefore, we can follow the arguments in the proofs of Proposition \ref{proposition2.19} to show the following result.
		
		\begin{enumerate}
			\item If $\alpha_1 \leq -\alpha_{n}$, then we have
			\begin{equation*}\begin{array}{l}
					\mathlarger{\sum}_{j \in \lbrace 0,i_1,\dots,i_t \rbrace}\mathlarger{\sum}_{Y \in S_h(\widetilde{Y})_{j}}	\dfrac{\CF_{\alpha^{\vee}}(Y,0)\CG(Y,B^{\vee})}{\alpha(Y,n)}f_{2r}(Y^{\wedge_{n-h+1}})\\
					=q^{n\delta_{\alpha^{\vee}}-(n-1)\delta_{\overline{\alpha^{\vee}}^1}+\alpha_1-2n+1}\mathlarger{\sum}_{\overline{Y} \in \overline{S}_h(\widetilde{Y})}\dfrac{\CF_{\overline{\alpha^{\vee}}^1}(\overline{Y},0)\CG(\overline{Y},\overline{B}^{\vee})}{\alpha(\overline{Y},n-1)}f_{2r}(\overline{Y}^{\wedge_{n-h}}).
				\end{array}
			\end{equation*}
			
			\item For $n-h+2 \leq l \leq n$, assume that $-\alpha_{l}<\alpha_1 \leq -\alpha_{l-1}$. Then, we have
			\begin{equation*}\begin{array}{l}
					\mathlarger{\sum}_{j \in \lbrace 0,i_1,\dots,i_t \rbrace}\mathlarger{\sum}_{Y \in S_h(\widetilde{Y})_{j}}	\dfrac{\CF_{\alpha^{\vee}}(Y,0)\CG(Y,B^{\vee})}{\alpha(Y,n)}f_{2r}(Y^{\wedge_{n-h+1}})\\
					=q^{n\delta_{\alpha^{\vee}}-(n-1)\delta_{\overline{\alpha^{\vee}}^1}+\alpha_1-2n+1} \mathlarger{\sum}_{\overline{Y} \in \overline{S}_h(\widetilde{Y})}
					\dfrac{\CG(\overline{Y},\overline{B}^{\vee})}{\alpha(\overline{Y},n-1)}f_{2r}(\overline{Y}^{\wedge_{n-h}})
					\times \lbrace \CF_{\overline{\alpha^{\vee}}^{1}}(\overline{Y},0)((-q)^{-\alpha_1-\alpha_{n}})\\
					+\sum_{j=2}^{n-l+1} \CF_{\overline{\alpha^{\vee}}^{j}}(\overline{Y},0)((-q)^{-\alpha_1-\alpha_{n+1-j}}-(-q)^{-\alpha_1-\alpha_{n+2-j}})
					+\CF_{\overline{\alpha^{\vee}}^{n-l+2}}(\overline{Y},0)(1-(-q)^{-\alpha_1-\alpha_{l}}) \rbrace.
				\end{array}
			\end{equation*}	
		\end{enumerate}
		By \eqref{eq2.24} and \eqref{eq2.25}, we can follow the arguments in the proof of Proposition \ref{proposition2.10} to prove the proposition.
	\end{proof}

	\section{Unitary Rapoport-Zink spaces and special cycles}
	\subsection{The moduli space $\CN^{[\mathds{h}]}_{n}$}
	In this subsection, we review the definition of the Rapoport-Zink space $\CN_n^{[\mathds{h}]}$. Here $\mathds{h}:=(h_1,h_2,\dots,h_m)$ a tuple of integers such that $0\leq h_1<h_2<\dots<h_m \leq n$ and $h_1 \equiv h_2 \equiv \dots \equiv h_m (\tx{mod } 2)$. We follow the notation of \cite{RZ}, \cite{LRZ}, and \cite{Muller} whenever possible.
	\begin{definition}
		For a $\Spf O_{\breve{F}}$-scheme $S$, and $0 \leq h \leq n$, a hermitian formal $O_F$-module of signature $(1,n-1)$ and type $h$ over $S$ is a triple $(X,i_X,\lambda_X)$ where
		\begin{enumerate}
			\item $X$ is a strict formal $O_{F_0}$-module over $S$ of relative height $2n$ and dimension $n$,
			
			\item $i_X:O_F \rightarrow \End(X)$ is an action of $O_F$ on $X$ extending the action of $O_{F_0}$ and satisfying the Kottwitz condition of signature $(1,n-1)$:
			for all $a \in O_F$, $\tx{Charpol}(i_{X}(a)\mid \Lie(X))=(T-a)(T-a^*)^{n-1} \in O_S[T]$,
			
			\item $\lambda_X :X \rightarrow X^{\vee}$ is a polarization such that its Rosati involution induces the involution $^*$ on $O_F$. We also require that $\ker(\lambda_X) \subset X[\pi]$, and its order is $q^{2h}$.
		\end{enumerate}
	\end{definition}
	
	By \cite[Remark 3.32]{cho2018basic}, a hermitian formal $O_F$-module of signature $(1,n-1)$ and type $h$ over $\BF$ is unique up to quasi-isogeny. We fix one choice of $(\BX_h,i_{\BX_h},\lambda_{\BX_h})$ as the framing object.
	
	For $0 \leq h \leq h' \leq n$, $h \equiv h' (\tx{mod }2)$, there exists an isogeny
	\begin{equation*}
		\alpha_{h',h}:\BX_{h'} \rightarrow \BX_{h},
	\end{equation*}
	compatible with polarizations such that $\ker (\alpha_{h',h}) \subset \BX_{h'}[\pi]$ and its order if $q^{h'-h}$. We fix $\alpha_{h+2,h}$ for each $0 \leq h \leq n-2$ and for $h'>h$, $h \equiv h' (\tx{mod }2)$, we define
	\begin{equation*}
		\alpha_{h',h}:=\alpha_{h+2,h}\circ \dots \circ \alpha_{h',h'-2}:\BX_{h'} \rightarrow \BX_{h'-2} \rightarrow \dots \rightarrow \BX_{h}.
	\end{equation*}

	\begin{definition}\label{definition3.3.2}
		Let (Nilp) be the category of $O_{\breve{F}}$-schemes $S$ such that $\pi$ is locally nilpotent on $S$. Let $\mathds{h}=(h_1,h_2,\dots,h_m)$ a tuple of integers such that $0\leq h_1<h_2<\dots<h_m \leq n$ and $h_1 \equiv h_2 \equiv \dots \equiv h_m (\tx{mod } 2)$. Then the functor $\CN_n^{[\mathds{h}]}$ associates to a scheme $S \in$ (Nilp) the set of isomorphism classes of tuples $(X_{h_i},i_{X_{h_i}},\lambda_{X_{h_i}},\rho_{X_{h_i}})_{1 \leq i \leq m}$, where
		\begin{enumerate}
			\item $(X_{h_i},i_{X_{h_i}},\lambda_{X_{h_i}})$ is a hermitian formal $O_F$-module of signature $(1,n-1)$ and type $h_i$ over $S$,
			\item $\rho_{X_{h_i}}:X_{h_i} \times_S \overline{S} \rightarrow \BX_{h_i} \times_{\BF} \overline{S}$ is an $O_F$-isogeny of height $0$ over the reduction $\overline{S}=S \times_{O_{\breve{F}}}\BF$ such that $\rho_{X_{h_i}}^*(\lambda_{\BX_{h_i},\overline{S}})=\lambda_{X_{h_i},\overline{S}}$,
			\item For $1 \leq i <m$, $\rho_{X_{h_i}}^{-1} \circ \alpha_{h_{i+1},h_i}\circ \rho_{X_{h_{i+1}}}$ lifts to an isogeny $\tilde{\alpha}_{h_{i+1},h_i}:X_{h_{i+1}}\rightarrow X_{h_i}$,
			\begin{center}
				\begin{tikzcd}[column sep=huge, row sep=large]
					X_{h_{i+1},\overline{S}} \arrow{r}{\tilde{\alpha}_{h_{i+1},h_i,\overline{S}}}
					\arrow{d}[swap]{\rho_{X_{h_{i+1},\overline{S}}}}
					& X_{h_i,\overline{S}} \arrow{d}{\rho_{X_{h_{i},\overline{S}}}}  \\
					\BX_{h_{i+1},\overline{S}}
					\arrow{r}{\alpha_{h_{i+1},h_i,\overline{S}}}
					& \BX_{h_i,\overline{S}}
				\end{tikzcd}
			\end{center}
		\end{enumerate}
		Note that if $\tilde{\alpha}_{h_{i+1},h_i}$ exists, then it is unique (\cite[Section 3.4]{LRZ}). The functor is representable by a formal scheme over $\Spf O_{\breve{F}}$ which is locally formally of finite type and regular (\cite{RZ}, \cite{Gortz}).
		
	\end{definition}
	
	\begin{definition}\label{definition3.3.3} Let $\mathds{h}=(h_1,\dots,h_m)$ a tuple of integers such that $0\leq h_1<\dots<h_m \leq n$ and $h_1 \equiv \dots \equiv h_m (\tx{mod } 2)$.
		\begin{enumerate}
			\item We denote by $\mathds{h}^{\vee}:=(n-h_m,n-h_{m-1},\dots,n-h_1)$.
			\item We denote by $\theta:\CN_n^{[\mathds{h}]} \rightarrow \CN_n^{[\mathds{h}^{\vee}]}$ the isomorphism defined as follows (cf. \cite[Remark 5.2]{cho2018basic}): For $S \in$ (Nilp),
			\begin{equation*}
				\begin{array}{ccc}
					\CN_n^{[\mathds{h}]}(S) &\overset{\theta}{\longrightarrow} &\CN_n^{[\mathds{h}^{\vee}]}(S)\\
					(X_{h_i},i_{X_{h_i}},\lambda_{X_{h_i}},\rho_{X_{h_i}})_{1 \leq i \leq m} & \mapsto &(X_{h_i}^{\vee},\overline{i}^{\vee}_{X_{h_i}},\lambda_{X_{h_i}}',(\rho^{\vee}_{X_{h_i}})^{-1})_{1 \leq i \leq m},
				\end{array}
			\end{equation*}
			where $\lambda_{X_{h_i}}':X_{h_i}^{\vee} \rightarrow X_{h_i}$ is the unique polarization such that $\lambda_{X_{h_i}}' \circ \lambda_{X_{h_i}}=i_{X_{h_i}}(\pi)$, and for $a \in O_F$, the action $\overline{i}^{\vee}_{X_{h_i}}$ is defined as $\overline{i}^{\vee}_{X_{h_i}}(a)=i_{X_{h_i}}(a^*)^{\vee}$.
		\end{enumerate}
	\end{definition}

	\subsection{Special cycles}
	In this subsection, we define special cycles in $\CN^{[\mathds{h}]}_n$. We fix a hermitian formal $O_F$-module $(\overline{\BE},i_{\overline{\BE}},\lambda_{\overline{\BE}})$ of signature $(0,1)$ and type $0$ over $\BF$. More precisely, $(\overline{\BE},i_{\overline{\BE}},\lambda_{\overline{\BE}})$ is a triple where
	\begin{enumerate}
		\item $\overline{\BE}$ is a strict formal $O_{F_0}$-module over $\BF$ of relative height $2$ and dimension $1,$
		\item $i_{\overline{\BE}}$ is an $O_F$-action satisfying the Kottwitz signature condition $(0,1)$,
		\item $\lambda_{\BE}$ is a principal polarization such that its Rosati involution induces the involution $^*$ on $O_F$.
	\end{enumerate}
	We then choose this triple as a framing object and consider the corresponding moduli space, which we denote by $\CN^0$.

	\begin{definition}\label{definition3.3.4}
		Let $\mathds{h}=(h_1,\dots,h_m)$ a tuple of integers such that $0\leq h_1<\dots<h_m \leq n$ and $h_1 \equiv \dots \equiv h_m (\tx{mod } 2)$. We define the space of special homomorphisms to be the $F$-vector space
		\begin{equation*}
			\BV:=\Hom_{\O_F}(\overline{\BE},\BX_{h_m})\otimes_{\BZ}\BQ.
		\end{equation*}
		For $x,y \in \BV$, we define a hermitian form $h$ on $\BV$ as
		\begin{equation*}
			h(x,y)=\lambda_{\overline{\BE}}^{-1}\circ y^{\vee}\circ \lambda_{\BX_{h_m}}\circ x \in \End_{O_F}(\overline{\BE})\otimes_{\BZ}\BQ \overset{i_{\BX}^{-1}}{\simeq} F.
		\end{equation*}
	\end{definition}
	
	\begin{remark}
		The space of special homomorphisms can be defined using $\BX_{h_i}$ instead of $\BX_{h_m}$. Since $\alpha_{h_m,h_i}:\BX_{h_m} \rightarrow \BX_{h_i}$ is compatible with polarizations $\lambda_{\BX_{h_i}}$ and $\lambda_{\BX_{h_m}}$, the two definitions are isometric.
	\end{remark}
	
	\begin{definition}\label{definition3.3.6}
		Let $\mathds{h}=(h_1,\dots,h_m)$ a tuple of integers such that $0\leq h_1<\dots<h_m \leq n$ and $h_1 \equiv \dots \equiv h_m (\tx{mod } 2)$. For the framing object
		
		\[
		\begin{tikzcd}[column sep=1.7cm, row sep=1cm]
			\BX_{h_m} \arrow{r}{\alpha_{h_m,h_{m-1}}} \arrow{d}{\lambda_{\BX_{h_{m}}}}
			&
			\BX_{h_{m-1}} \arrow{r}{\alpha_{h_{m-1},h_{m-2}}} \arrow{d}{\lambda_{\BX_{h_{m-1}}}}
			&
			\cdots \arrow{r}{\alpha_{h_{3},h_{2}}}
			&
			\BX_{h_{2}} \arrow{r}{\alpha_{h_{2},h_{1}}} \arrow{d}{\lambda_{\BX_{h_{2}}}}
			&
			\BX_{h_{1}} \arrow{d}{\lambda_{\BX_{h_{1}}}}
			\\
			\BX_{h_{m}}^{\vee}
			&
			\BX_{h_{m-1}}^{\vee} \arrow{l}{\alpha_{h_{m},h_{m-1}}^{\vee}}
			&
			\cdots \arrow{l}{\alpha_{h_{m-1},h_{m-2}}^{\vee}}
			&
			\BX_{h_{2}}^{\vee} \arrow{l}{\alpha_{h_{3},h_{2}}^{\vee}}
			&
			\BX_{h_{1}}^{\vee}\arrow{l}{\alpha_{h_{2},h_{1}}^{\vee}}\tx{ },
		\end{tikzcd}
		\]
		
		and $x \in \BV=\Hom_{\O_F}(\overline{\BE},\BX_{h_m})\otimes_{\BZ}\BQ$, we define the special cycles $\CZ_{h_m}(x)$, $\CZ_{h_{m-1}}(x)$, $\dots,\CZ_{h_1}(x)$, $\CY_{h_1}(x)$, $\CY_{h_2}(x)$, $\dots$, and $\CY_{h_m}(x)$ as follows.
		\begin{enumerate}
			\item We define the special cycle $\CZ_{h_i}(x)$ to be the closed formal subscheme of $\CN^0 \times \CN^{[\mathds{h}]}_{n}$ with the following property:
			For each $O_{\breve{F}}$-scheme $S$ such that $\pi$ is locally nilpotent, $\CZ_{h_i}(x)(S)$ is the set of all points $(\overline{\CE},i_{\overline{\CE}},\lambda_{\overline{\CE}},\rho_{\overline{\CE}},(X_{h_j},i_{X_{h_j}},\lambda_{X_{h_j}},\rho_{X_{h_j}})_{1 \leq j \leq m})$ in $\CN^0 \times \CN^{[\mathds{h}]}_{n}(S)$ such that the quasi-homomorphism
			\begin{equation*}
				\rho_{X_{h_i}}^{-1} \circ \alpha_{h_m,h_i} \circ x \circ \rho_{\overline{\CE}}:\overline{\CE} \times_S \overline{S} \rightarrow X_{h_i} \times_S \overline{S}
			\end{equation*}
			extends to a homomorphism from $\overline{\CE}$ to $X_{h_i}$.
			
			\item We define the special cycle $\CY_{h_i}(x)$ in $\CN^0 \times \CN^{[\mathds{h}]}_{n}$ as follows. First, we consider the special cycle $\CZ_{n-h_i}(\lambda_{h_1}\circ \alpha_{h_m,h_1}\circ x)$ in $\CN^0 \times \CN^{[\mathds{h}^{\vee}]}_n$. We define $\CY_{h_i}(x)$ as $(id \times \theta^{-1})(\CZ_{n-h_i}(\lambda_{h_1}\circ \alpha_{h_m,h_1}\circ x))$ in $\CN^0 \times \CN^{[\mathds{h}]}_n$.
		\end{enumerate}
		
		Note that $\CN^0$ can be identified with $\Spf O_{\breve{F}}$, and hence $\CZ_{h_i}(x), \CY_{h_i}(x)$, $1 \leq i \leq m$, can be identified with closed formal subschemes of $\CN_n^{[\mathds{h}]}$.
	\end{definition}
	
	\begin{proposition}[cf. {\cite[Proposition 5.9]{cho2018basic}}] Let $\mathds{h}=(h_1,\dots,h_m)$ a tuple of integers such that $0\leq h_1<\dots<h_m \leq n$ and $h_1 \equiv \dots \equiv h_m (\tx{mod } 2)$. For $x \in \BV \setminus \lbrace 0 \rbrace$, the functors $\CZ_{h_i}(x)$, and $\CY_{h_i}(x)$, $1 \leq i \leq m$ are represented by Cartier divisors in $\CN^0 \times \CN_n^{[\mathds{h}]}$.
	\end{proposition}
	\begin{proof}
		The proof is almost identical to the proof of \cite[Proposition 5.9]{cho2018basic}.
	\end{proof}
	
	\begin{proposition}[cf. {\cite[Proposition 5.10]{cho2018basic}}] \label{proposition3.8} Let $\mathds{h}=(h_1,\dots,h_m)$ a tuple of integers such that $0\leq h_1<\dots<h_m \leq n$ and $h_1 \equiv \dots \equiv h_m (\tx{mod } 2)$. Let $x,y \in \BV \setminus \lbrace 0 \rbrace$ be special homomorphisms such that $\val(h(x,x))=0$ and $\val(h(y,y))=-1$.
		\begin{enumerate}
			\item If $h_m<n$, then $\CZ_{h_m}(x) \simeq \CN_{n-1}^{[\mathds{h}]}$.
			\item If $h_m=n$, then $\CZ_{h_m}(x)= \emptyset$.
			\item If $h_1>0$, then $\CY_{h_1}(y) \simeq \CN_{n-1}^{[\mathds{h}-1]}$, where $\mathds{h}-1=(h_1-1,h_2-1,\dots,h_m-1)$.
			\item If $h_1=0$, then $\CY_{h_1}(y)=\emptyset$.
		\end{enumerate}
	\end{proposition}
	\begin{proof} The proof is almost identical to the proof of \cite[Proposition 5.10]{cho2018basic}. Indeed, the isomorphism from $\CN_{n-1}^{[\mathds{h}]} \rightarrow \CZ_{h_m}(x)(S)$ is given by
		\begin{equation*}
			(X_{h_i},i_{X_{h_i}},\lambda_{X_{h_i}},\rho_{X_{h_i}})_{1 \leq i \leq m} \mapsto ((\overline{\CE},i_{\overline{\CE}},\lambda_{\overline{\CE}},\rho_{\overline{\CE}}),(\overline{\CE} \times X_{h_i},i_{\overline{\CE}} \times i_{X_{h_i}},\lambda_{\overline{\CE}}\times \lambda_{X_{h_i}},\rho_{\overline{\CE}}\times \rho_{X_{h_i}})_{1 \leq i \leq m}),
		\end{equation*}
		with $x=\tx{inc}_1: \overline{\BE} \rightarrow \overline{\BE}\times \BX_{h_m}$.
	\end{proof}

	\begin{proposition}[cf. {\cite[Proposition 5.11]{cho2018basic}}] \label{proposition3.9}
		Let $\mathds{h}=(h_1,\dots,h_m)$ be a tuple of integers such that $0\leq h_1<\dots<h_m \leq n$ and $h_1 \equiv \dots \equiv h_m (\tx{mod } 2)$. Assume that $x,y \in \BV \setminus \lbrace 0 \rbrace$ satisfy $\val(h(x,x))=0$ and $\val(h(y,y))=-1$. Assume further that, by rescaling by an element of $O_F^{\times}$, $x^{*}\circ x=1$ and $(\lambda_{\BX_{h_m}}\circ y)^* \circ (\lambda_{\BX_{h_m}}\circ y)=1$, where $x^*=\lambda_{\overline{\BE}}^{-1} \circ x^{\vee} \circ \lambda_{\BX_{h_m}}$ and $(\lambda_{\BX_{h_m}}\circ y)^*=\lambda_{\overline{\BE}}^{-1} \circ (\lambda_{\BX_{h_m}}\circ y)^{\vee} \circ \lambda_{\BX_{h_m}}'$. We define $e_x:=x \circ x^{*}$ and $e_y:=(\lambda_{\BX_{h_m}}\circ y)^{*}\circ (\lambda_{\BX_{h_m}}\circ y)$. Fix isomorphisms
		\begin{equation*}
			\begin{array}{l}
				\Phi:\CZ_{h_m}(x) \simeq \CN_{n-1}^{[\mathds{h}]},\\
				\Psi: \CY_{h_1}(y) \simeq \CN_{n-1}^{[\mathds{h}-1]},
			\end{array}
		\end{equation*}
		as in Proposition \ref{proposition3.8}. Then, the following statements hold.
		\begin{enumerate}
			\item For $z \in \BV$ with $h(x,z)=0$, let $z':=(1-e_x)\circ z$. Then, we have
			\begin{equation*}
				\Phi(\CZ_{h_m}(x) \cap \CZ_{h_i}(z))=\CZ_{h_i}(z') \tx{ in } \CN_{n-1}^{[\mathds{h}]},\tx{ and }h(z',z')=h(z,z).
			\end{equation*}
			\item For $w \in \BV$ with $h(x,w)=0$, let $w':=(1-e_x)\circ w$. Then, we have
			\begin{equation*}
				\Phi(\CZ_{h_m}(x) \cap \CY_{h_i}(w))=\CY_{h_i}(w') \tx{ in } \CN_{n-1}^{[\mathds{h}]},\tx{ and }h(w',w')=h(w,w).
			\end{equation*}
			\item For $z \in \BV$ with $h(y,z)=0$, let $z':=(1-e_y^{\vee})\circ z$. Then, we have
			\begin{equation*}
				\Psi(\CY_{h_1}(y) \cap \CZ_{h_i}(z))=\CZ_{h_i-1}(z') \tx{ in } \CN_{n-1}^{[\mathds{h}-1]},\tx{ and }h(z',z')=h(z,z).
			\end{equation*}
			\item For $w \in \BV$ with $h(y,w)=0$, let $w':=(1-e_y^{\vee})\circ w$. Then, we have
			\begin{equation*}
				\Psi(\CY_{h_1}(y) \cap \CY_{h_i}(w))=\CY_{h_i-1}(w') \tx{ in } \CN_{n-1}^{[\mathds{h}-1]},\tx{ and }h(w',w')=h(w,w).
			\end{equation*}
		\end{enumerate}
	\end{proposition}
	\begin{proof}
		The proof is almost identical to the proof of \cite[Proposition 5.12]{cho2018basic}.
	\end{proof}
	
	\begin{definition}\label{definition3.3.10}
		Let $0 < m \leq n$ be integers, and let $\mathds{h}=(h_1,\dots,h_m)$ a tuple of integers such that $0\leq h_1<\dots<h_m \leq n$ and $h_1 \equiv \dots \equiv h_m (\tx{mod } 2)$.
		\begin{enumerate}
			\item We define the subset $\CE_n^{\mathds{h}}$ of $\CE_n$ by
			\begin{equation*}
				\CE_n^{\mathds{h}}=\lbrace \alpha=(\alpha_1,\dots,\alpha_n) \in \CE_n \mid \forall i, \tx{ }\alpha_i=-h_j \tx{ or }h_j \tx{ for some }j \rbrace 
			\end{equation*}
			
			\item Fix $\alpha \in \CE_n^{\mathds{h}}$. Assume that $\alpha_1 \leq \alpha_2\leq \dots \alpha_{h-1} < 0 \leq \alpha_h \dots \leq \alpha_n$. Let $x_1$,$\dots$,$x_{h-1}$,$y_1$,$\dots$,$y_{n-h+1}$ $\in \BV$ be special homomorphisms. We consider special cycles
			\begin{equation*}
				\CZ_{-\alpha_1}(x_1),\CZ_{-\alpha_2}(x_2),\dots,\CZ_{-\alpha_{h-1}}(x_{h-1}),\CY_{\alpha_h}(y_1),\dots,\CY_{\alpha_n}(y_{n-h+1})
			\end{equation*}
			in $\CN:=\CN_{n}^{[\mathds{h}]}$. We denote by $\CI_{\alpha}^{\mathds{h}}(x_1,\dots,x_{h-1},y_1,\dots,y_{n-h+1})$ the arithmetic intersection number of $\CZ_{-\alpha_1}(x_1)$,$\CZ_{-\alpha_2}(x_2)$,$\dots$,$\CZ_{-\alpha_{h-1}}(x_{h-1})$,$\CY_{\alpha_h}(y_1)$,$\dots$,$\CY_{\alpha_n}(y_{n-h+1})$ in $\CN$. More precisely,
			\begin{equation*}\begin{array}{ll}
					\CI_{\alpha}^{\mathds{h}}(x_1,\dots,x_{h-1},y_1,\dots,y_{n-h+1})\\
					=\chi(O_{\CZ_{-\alpha_1}(x_1)} \otimes_{O_{\CN}}^{\BL} \dots \otimes_{O_{\CN}}^{\BL} O_{\CZ_{-\alpha_{h-1}}(x_{h-1})} \otimes_{O_{\CN}}^{\BL} O_{\CY_{\alpha_h}(y_{h})}) \otimes_{O_{\CN}}^{\BL}\dots \otimes_{O_{\CN}}^{\BL} O_{\CY_{\alpha_n}(y_{n-h+1})}).
				\end{array}
			\end{equation*}
			Here, we write $\chi$ for the Euler-Poincare characteristic and $\otimes^{\BL}_{O_{\CN}}$ for the derived tensor product of $O_{\CN}$-modules.
			
		\end{enumerate}
	\end{definition}
	
	\begin{proposition}\label{proposition3.11}
		Let $0 < m \leq n$ be integers, and let $\mathds{h}=(h_1,\dots,h_m)$ a tuple of integers such that $0\leq h_1<\dots<h_m \leq n$ and $h_1 \equiv \dots \equiv h_m (\tx{mod } 2)$. Assume that $\alpha \in \CE_n^{\mathds{h}}$. Assume that $\alpha_1 \leq \alpha_2\leq \dots \alpha_{h-1} < 0 \leq \alpha_h \dots \leq \alpha_n$. Let $x_1$,$\dots$,$x_{h-1}$,$y_1$,$\dots$,$y_{n-h+1}$ $\in \BV$ be special homomorphisms. Then, the following statements hold.
		\begin{enumerate}
			\item Assume that $\val (h(x_1,x_1))=0$, $\alpha_1=-h_m >-n$, and $x_1 \perp \lbrace x_2,\dots,x_{h-1},y_{1},\dots,y_{n-h+1}\rbrace$. Then,
			\begin{equation*}
				\CI_{\alpha}^{\mathds{h}}(x_1,\dots,x_{h-1},y_1,\dots,y_{n-h+1})=\CI_{\tx{}^1\overline{\alpha}}^{\mathds{h}}(x_2,\dots,x_{h-1},y_1,\dots,y_{n-h+1}).
			\end{equation*}
			Here, recall from Definition \ref{definition2.21} that $^1\overline{\alpha}:=(\alpha_2,\dots,\alpha_n) \in \CE_{n-1}^{\mathds{h}}$

			\item Assume that $\val (h(x_1,x_1))=0$, and $\alpha_{1}=-h_m=-n$. Then,
			\begin{equation*}
				\CI_{\alpha}^{\mathds{h}}(x_1,\dots,x_{h-1},y_1,\dots,y_{n-h+1})=0.
			\end{equation*}
			
			\item Assume that $\val (h(y_1,y_1))=-1$, $\alpha_h=-h_1 >0$, and $y_1 \perp \lbrace x_1,\dots,x_{h-1},y_{2},\dots,y_{n-h+1}\rbrace$. Then,
			\begin{equation*}
				\CI_{\alpha}^{\mathds{h}}(x_1,\dots,x_{h-1},y_1,\dots,y_{n-h+1})=\CI_{\overline{\alpha}^1}^{\mathds{h}-1}(x_1,\dots,x_{h-1},y_2,\dots,y_{n-h+1}).
			\end{equation*}
			Here, recall from Definition \ref{definition2.9} that $\overline{\alpha}^{1}:=(\alpha_1+1,\dots,\alpha_{h-1}+1,\alpha_{h+1}-1,\dots,\alpha_n-1) \in \CE_{n-1}^{\mathds{h}-1}.$

			\item Assume that $\val (h(y_1,y_1))=-1$, and $\alpha_{h}=h_1=0$. Then,
			\begin{equation*}
				\CI_{\alpha}^{\mathds{h}}(x_1,\dots,x_{h-1},y_1,\dots,y_{n-h+1})=0.
			\end{equation*}
		\end{enumerate}
	\end{proposition}
	\begin{proof}
		This follows from Proposition \ref{proposition3.8} and Proposition \ref{proposition3.9}.
	\end{proof}
	
	Recall from Definition \ref{definition3.3.3} the dual isomorphism $\theta:\CN_n^{[\mathds{h}]} \rightarrow \CN_n^{[\mathds{h}^{\vee}]}$, where $\mathds{h}=(h_1,\dots,h_m)$, and $\mathds{h}^{\vee}=(n-h_m,\dots,n-h_1)$. Under this isomorphism, we have that for $x,y \in \BV$,
	\begin{equation*}\begin{array}{l}
			\theta(\CY_{h_i}(y))=\CZ_{n-h_i}(\lambda_{h_1} \circ \alpha_{h_m,h_1}\circ y),\\
			\theta(\CZ_{h_i}(x))=\CY_{n-h_i}(\pi^{-1}\lambda_{h_1} \circ \alpha_{h_m,h_1}\circ x).
		\end{array}
	\end{equation*}
	This implies the following lemma.
	
	\begin{lemma}\label{lemma3.12}
		Let $0 < m \leq n$ be integers, and let $\mathds{h}=(h_1,\dots,h_m)$ a tuple of integers such that $0\leq h_1<\dots<h_m \leq n$ and $h_1 \equiv \dots \equiv h_m (\tx{mod } 2)$. Assume that $\alpha \in \CE_n^{\mathds{h}}$. Assume that $\alpha_1 \leq \alpha_2\leq \dots \alpha_{h-1} < 0 \leq \alpha_h \dots \leq \alpha_n$. Let $x_1$,$\dots$,$x_{h-1}$,$y_1$,$\dots$,$y_{n-h+1}$ $\in \BV$ be special homomorphisms. For a special homomorphism $z \in \BV$, we write $\check{z}$ for $\lambda_{h_1} \circ \alpha_{h_m,h_1}\circ z$. Then, we have
		\begin{equation*}\begin{array}{l}
				\CI_{\alpha}^{\mathds{h}}(x_1,\dots,x_{h-1},y_1,\dots,y_{n-h+1})
				=\CI_{\alpha^{\vee_{n-h+1}}}^{\mathds{h}^{\vee}}(\check{y}_1,\dots,\check{y}_{n-h+1},\pi^{-1}\check{x}_1,\dots,\pi^{-1}\check{x}_{h-1}).
			\end{array}
		\end{equation*}
		Here, recall from Definition \ref{definition2.20} that
		\begin{equation*}
			\alpha^{\vee_{n-h+1}}=(\alpha_h-n,\dots,\alpha_{n}-n,\alpha_1+n,\dots,\alpha_{h-1}+n) \in \CE_n^{\mathds{h}^{\vee}}.
		\end{equation*}
	\end{lemma}
	
	\begin{remark}\label{remark3.13} Let $B$ be the $n \times n$ matrix
		\begin{equation*}
			B=	\left(\begin{array}{cc} 
				h(x_i,x_j) & h(x_i,y_l)\\
				h(y_k,x_j)& h(y_k,y_l)
			\end{array} 
			\right)_{\substack{1\leq i,j \leq h-1\\ 1 \leq k,l \leq n-h+1}}.
		\end{equation*}
		Then, we have that
		\begin{equation*}
			B^{\vee_{n-h+1}}=\left(\begin{array}{cc} 
				h'(\check{y}_k,\check{y}_l) & h'(\check{y}_k,\pi^{-1}\check{x}_j)\\
				h'(\pi^{-1}\check{x}_i,\check{y}_l)& h'(\pi^{-1}\check{x}_i,\pi^{-1}\check{x}_j)
			\end{array} 
			\right)_{\substack{1 \leq k,l \leq n-h+1 \\ 1\leq i,j \leq h-1}},
		\end{equation*}
		where $h'$ is the hermitian form of special homomorphisms on the dual side. Here, we note that $	B^{\vee_{n-h+1}}$ is defined in Definition \ref{definition2.20}.
	\end{remark}

	\section{Conjectures on the arithmetic intersection numbers of special cycles}

	\subsection{Conjectures}
	In this subsection, we consider $\mathds{h}=(1,3,\dots,2n-1)$ and $\alpha=(-2n+1,-2n+3,\dots,-1,1,\dots,2n-1) \in \CE_{2n}^{\mathds{h}}$. In this setting, we formulate a conjecture for the arithmetic intersection number $\CI_{\alpha}^{\mathds{h}}(x_1,\dots,x_{n},y_1,\dots,y_{n})$, $x_i, y_i \in \BV$, of special cycles $\CZ_{2n-1}(x_1)$,$\dots$, $\CZ_{1}(x_n)$, $\CY_1(y_1)$,$\dots$, and $\CY_{2n-1}(y_n)$ in the Iwahori-level Rapoport--Zink space $\CN^{[\mathds{h}]}_{2n}$, in terms of weighted representation densities.
	
	We begin with the following remark.
	
	\begin{remark}\label{remark4.1} 	Let $0 < m \leq n$ be integers, and let $\mathds{h}=(h_1,\dots,h_m)$ a tuple of integers such that $0\leq h_1<\dots<h_m \leq n$ and $h_1 \equiv \dots \equiv h_m (\tx{mod } 2)$. Assume that $\alpha \in \CE_n^{\mathds{h}}$. Let $\lbrace x_1$,$\dots$,$x_{h-1}$,$y_1$,$\dots$,$y_{n-h+1} \rbrace$ be a basis of the space of special homomorphisms $\BV$. As in the hyperspecial level case (\cite{KR1}) and the maximal parahoric level case (\cite{Cho}), it is natural to expect that, in the general parahoric level case, the intersection number $\CI_{\alpha}^{[\mathds{h}]}(x_1$,$\dots$,$x_{h-1}$,$y_1$,$\dots$,$y_{n-h+1})$ depends only on the matrix $B$:
		\begin{equation*}
			B=	\left(\begin{array}{cc} 
				h(x_i,x_j) & h(x_i,y_l)\\
				h(y_k,x_j)& h(y_k,y_l)
			\end{array} 
			\right)_{\substack{1\leq i,j \leq h-1\\ 1 \leq k,l \leq n-h+1}}.
		\end{equation*}
		Assuming this, we write $\CI_{\alpha}^{[\mathds{h}]}(B)$ for $\CI_{\alpha}^{[\mathds{h}]}(x_1$,$\dots$,$x_{h-1}$,$y_1$,$\dots$,$y_{n-h+1})$. Then, Lemma \ref{lemma3.12} and Remark \ref{remark3.13} imply that
		\begin{equation*}
			\CI_{\alpha}^{[\mathds{h}]}(B)=\CI_{\alpha^{\vee_{n-h+1}}}^{[\mathds{h}^{\vee}]}(B^{\vee_{n-h+1}}).
		\end{equation*}
		
		Therefore, when $\mathds{h}=(1,3,\dots,2n-1)$ and $\alpha=(-2n+1,-2n+3,\dots,-1,1,\dots,2n-1) \in \CE_{2n}^{\mathds{h}}$, we expect that
		\begin{equation*}
			\CI_{\alpha}^{[\mathds{h}]}(B)=\CI_{\alpha^{\vee_n}}^{[\mathds{h}^{\vee}]}(B^{\vee_n})=\CI_{\alpha}^{[\mathds{h}]}(B^{\vee_n}).
		\end{equation*}
		This suggests that the conjectural formula for $\CI_{\alpha}^{[\mathds{h}]}(B)$ should be invariant under the change $B \leftrightarrow B^{\vee_n}$.
	\end{remark}

	\begin{theorem}\label{theorem3.1}
		For $\alpha=(-2n+1,-2n+3,\dots,-1,1,\dots,2n-1) \in \CE_{2n}^{\mathds{h}}$, we have
		\begin{equation*}
			q^{4n^2-2n}W_{\alpha}'(B)-\mathlarger{\sum}_{\beta}C_{\beta,2n-1-\beta_{2n}}W_{\beta}(B,0)=q^{4n^2-2n}W_{\alpha}'(B^{\vee_n})-\mathlarger{\sum}_{\beta}C_{\beta,2n-1-\beta_{2n}}W_{\beta}(B^{\vee_n},0),
		\end{equation*}
		where $C_{\beta,i}$ is defined by
		\begin{equation*}
			C_{\beta,i}:=\left\lbrace \begin{array}{ll}
				q^{4n^2-2n\delta_{\beta}-\sum_{i=1}^{2n} \beta_i}(-1)^k\lbrace \dfrac{1}{q^{2k-1}+1}+\dfrac{1}{q^{2k}-1} \rbrace & \tx{if }i=2k-1,\\
				q^{4n^2-2n\delta_{\beta}-\sum_{i=1}^{2n} \beta_i}(-1)^k\lbrace \dfrac{1}{q^{2k-2}-1}+\dfrac{1}{q^{2k-1}+1} \rbrace & \tx{if }i=2k-2,\\0 &\tx{if } i=0.
			\end{array}\right.
		\end{equation*}
		Here, $\beta$ runs over the subset of $\CE_{2n}$ satisfying the following conditions:
		\begin{enumerate}
			\item $\beta_i=-\beta_{2n+1-i}$ for all $1 \leq i \leq 2n$.
			\item $\alpha-\beta$ is increasing, i.e., $\alpha_i-\beta_i \leq \alpha_{i+1}-\beta_{i+1}$ for all $1 \leq i \leq 2n-1$.
		\end{enumerate}
		
	\end{theorem}
	
	\begin{proof}
		By Proposition \ref{proposition2.3}, we have that
		\begin{equation}\label{eq3.1}
			\begin{array}{ll}
				q^{4n^2-2n}W_{\alpha}'(B)-\mathlarger{\sum}_{\beta}C_{\beta,2n-1-\beta_{2n}}W_{\beta}(B,0)\\=\mathlarger{\sum}_{Y \in \CR_{2n}}\dfrac{\CG(Y,B)}{\alpha(Y,2n)}\lbrace q^{4n^2-2n}\CF_{\alpha}'(Y)-\mathlarger{\sum}_{\beta}C_{\beta,2n-1-\beta_{2n}}\CF_{\beta}(Y,0) \rbrace,
			\end{array}
		\end{equation}
		
		Also, by Proposition \ref{proposition2.3}, Proposition \ref{proposition2.22}, Proposition \ref{proposition2.23}, and Proposition \ref{proposition2.24}, we have
		\begin{equation}\label{eq3.2}
			\begin{array}{ll}
				q^{4n^2-2n}W_{\alpha}'(B^{\vee_n})-\mathlarger{\sum}_{\beta}C_{\beta,2n-1-\beta_{2n}}W_{\beta}(B^{\vee_n},0)\\=\mathlarger{\sum}_{Y \in \CR_{2n}}\dfrac{\CG(Y,B^{\vee_n})}{\alpha(Y,n)}\lbrace q^{4n^2-2n}\CF_{\alpha}'(Y)-\mathlarger{\sum}_{\beta}C_{\beta,2n-1-\beta_{2n}}\CF_{\beta}(Y,0) \rbrace\\
				=\mathlarger{\sum}_{Y \in \CR_{2n}}\dfrac{\CG(Y,B)}{\alpha(Y,2n)}\lbrace q^{4n^2-2n}\CF_{\alpha}'(Y^{\wedge_n})-\mathlarger{\sum}_{\beta}C_{\beta,2n-1-\beta_{2n}}\CF_{\beta}(Y^{\wedge_n},0) \rbrace.
			\end{array}
		\end{equation}
		
		Therefore, by \eqref{eq2.1} and \eqref{eq2.2}, it suffices to show that
		\begin{equation*}
			q^{4n^2-2n}\CF_{\alpha}'(Y)-\mathlarger{\sum}_{\beta}C_{\beta,2n-1-\beta_{2n}}\CF_{\beta}(Y,0)=q^{4n^2-2n}\CF_{\alpha}'(Y^{\wedge_n})-\mathlarger{\sum}_{\beta}C_{\beta,2n-1-\beta_{2n}}\CF_{\beta}(Y^{\wedge_n},0).
		\end{equation*}
		
		This equality will be proved in Lemma \ref{lemma3.2} below. This finishes the proof of the theorem.
	\end{proof}

	\begin{conjecture}\label{conjecture4.3} 
		Let $\mathds{h}=(1,3,\dots,2n-1)$ and $\alpha=(-2n+1,-2n+3,\dots,-1,1,\dots,2n-1) \in \CE_{2n}^{\mathds{h}}$. Let $\lbrace x_1$,$\dots$,$x_{n}$,$y_1$,$\dots$,$y_{n} \rbrace$ be a basis of the space of special homomorphisms $\BV$. Then, we have
		\begin{equation*}\begin{array}{ll}
				\CI_{\alpha}^{[\mathds{h}]}(x_1$,$\dots$,$x_{n}$,$y_1$,$\dots$,$y_{n})&=\dfrac{1}{W_{\alpha}(A_{2n,1},0)}	\lbrace W_{\alpha}'(B)-q^{-4n^2+2n}\mathlarger{\sum}_{\beta}C_{\beta,2n-1-\beta_{2n}}W_{\beta}(B,0) \rbrace\\
				&=\dfrac{q^{6n^2-n}}{(q+1)^2(q^2-1)^{n-1}}\lbrace W_{\alpha}'(B)-q^{-4n^2+2n}\mathlarger{\sum}_{\beta}C_{\beta,2n-1-\beta_{2n}}W_{\beta}(B,0) \rbrace.
			\end{array}
		\end{equation*}
		Here, $B$ is the $2n \times 2n$ matrix
		\begin{equation*}
			B=	\left(\begin{array}{cc} 
				h(x_i,x_j) & h(x_i,y_l)\\
				h(y_k,x_j)& h(y_k,y_l)
			\end{array} 
			\right)_{\substack{1\leq i,j \leq n\\ 1 \leq k,l \leq n}}.
		\end{equation*}
		We note that the matrix $A_{2n,1}$ is defined in Definition \ref{definition2.1}.
	\end{conjecture}

	\begin{remark}\label{remark4.4}
		In Conjecture \ref{conjecture4.3}, the constant $W_{\alpha}(A_{2n,1},0)$ is related to the volume of the Iwahori group associated with $\mathds{h}=(1,3,\dots,2n-1)$, with respect to our chosen measure.
	\end{remark}
	
	\begin{remark}\label{remark4.5'}
		For $\mathds{h}=(h_1,\dots,h_m)$, we have $h_1 \not\equiv \val(\det(B)) \pmod2$. Therefore, for $\beta=(\beta_1,\dots,\beta_n) \in \CE_{n}$ such that $h_1 \equiv \beta_1 \pmod2$, we have $W_{\beta}(B,0)=0$.
		
		\end{remark}
	
	\begin{remark}\label{remark4.5}
		Combining Proposition \ref{proposition3.11}, Lemma \ref{lemma3.12}, and their analytic counterparts---Proposition \ref{proposition2.19}, Proposition \ref{proposition2.25}, and Theorem \ref{theorem3.1}---Conjecture \ref{conjecture4.3} gives conjectural formulas for the following arithmetic intersection numbers:
		\begin{enumerate}
			\item $\CI_{\alpha'}^{[\mathds{h}']}(x_1$,$\dots$,$x_{n}$,$y_1$,$\dots$,$y_{n-1})$, $\alpha'=(-2n+2,-2n+4,\dots,0,\dots,2n-2)$, $\mathds{h}'=(0,2,\dots,2n-2)$;
			
			\item $\CI_{\alpha'}^{[\mathds{h}']}(x_1$,$\dots$,$x_{n-1}$,$y_1$,$\dots$,$y_{n})$, $\alpha'=(-2n+3,\dots,-1,1,\dots,2n-1)$, $\mathds{h}'=(1,3,\dots,2n-1)$;
			
			\item $\CI_{\alpha'}^{[\mathds{h}']}(x_1$,$\dots$,$x_{n-1}$,$y_1$,$\dots$,$y_{n-1})$, $\alpha'=(-2n+4,\dots,0,2,\dots,2n-2)$, $\mathds{h}'=(0,2,\dots,2n-2)$;
			
			\item $\CI_{\alpha'}^{[\mathds{h}']}(x_1$,$\dots$,$x_{n}$,$y_1$,$\dots$,$y_{n-1})$, $\alpha'=(-2n+1,\dots,-1,1,\dots,2n-3)$, $\mathds{h}'=(1,3,\dots,2n-3)$;
			
			\item $\CI_{\alpha'}^{[\mathds{h}']}(x_1$,$\dots$,$x_{n-1}$,$y_1$,$\dots$,$y_{n-1})$, $\alpha'=(-2n+2,\dots,0,2,\dots,2n-4)$, $\mathds{h}'=(0,2,\dots,2n-2)$.
		\end{enumerate}
		Therefore, Conjecture 4.3 gives conjectural formulas for the arithmetic intersection numbers of special cycles at all Iwahori levels and in arbitrary dimension, when the types of the special cycles are all distinct.

	\end{remark}
	
	\begin{remark}\label{remark4.6}
		Let $\mathds{h}=(h_1,\dots,h_m)$ a tuple of integers such that $0 \leq h_1<\dots<h_m \leq n$ and $h_1 \equiv \dots \equiv h_m \tx{ (mod 2})$. Let $\mathds{h'}=(h_1',\dots,h_l')$ a tuple of integers such that $l>m$, $0 \leq h_1'<\dots<h_l'\leq n$, $h_1' \equiv \dots \equiv h_l' \tx{ (mod 2})$. Assume further that for all $1 \leq i \leq m$, $h_i=h_j'$ for some $1 \leq j \leq l$. For example $\mathds{h}=(1,3)$ and $\mathds{h}'=(1,3,5)$. Let $\alpha \in \CE_{n}^{\mathds{h}} \subset \CE_{n}^{\mathds{h}'}$ and suppose that $\alpha_1 \leq \dots \leq \alpha_{h-1} < 0 \leq \alpha_h\leq \dots \leq \alpha_n$. Let $x_1,\dots,x_{h-1},y_1,\dots,y_{n-h+1} \in \BV$ be special homomorphisms. Then the number $\CI_{\alpha}^{[\mathds{h}]}(x_1$,$\dots$,$x_{h-1}$,$y_1$,$\dots$,$y_{n-h+1})$ (resp. $\CI_{\alpha}^{[\mathds{h}']}(x_1$,$\dots$,$x_{h-1}$,$y_1$,$\dots$,$y_{n-h+1})$) is the arithmetic intersection number of special cycles
		\begin{equation*}
			\CZ_{-\alpha_1}(x_1),\dots,\CZ_{-\alpha_{h-1}}(x_{h-1}),\CY_{\alpha_h}(y_1),\dots,\CY_{\alpha_{n}}(y_{n-h+1})
		\end{equation*}
		in $\CN_n^{[\mathds{h}]}$ (resp. $\CN_n^{[\mathds{h}']}$). Then, it is natural to ask if
		\begin{equation}\label{eq4.4.3}
			\CI_{\alpha}^{[\mathds{h}']}(x_1,\dots,x_{h-1},y_1,\dots,y_{n-h+1})\overset{?}{=}\dfrac{\tx{Vol}(P^{[\mathds{h}]})}{\tx{Vol}(P^{[\mathds{h}']})}	\CI_{\alpha}^{[\mathds{h}]}(x_1,\dots,x_{h-1},y_1,\dots,y_{n-h+1}),
		\end{equation}
		where $P^{[\mathds{h}]}$ (resp. $P^{[\mathds{h}']}$) is the parahoric group associated with $\mathds{h}$ (resp. $\mathds{h}'$).
		
		When $n=2$, \eqref{eq4.4.3} follows from the projection formula.
	\end{remark}

	\begin{remark}\label{remark4.7}
		It is natural to expect that the strategy of formulating conjectures via the duality of special cycles can also be applied in other settings. For example, when $\mathds{h}=(1,3)$ and $\alpha=(-3,-3,1,3) \in \CE_{4}^{\mathds{h}}$, by using a similar strategy, we conjecture that the arithmetic intersection number of special cycles $\CZ_3(x_1),\CZ_3(x_2),\CY_1(y_1),\CY_3(y_2)$ in $\CN_4^{[\mathds{h}]}$ is
		\begin{equation}\label{eq4.4.4}
			\CI_{\alpha}^{[\mathds{h}]}(x_1,x_2,y_1,y_2)=\dfrac{q^{22}}{(q+1)^2(q^2-1)}\lbrace W_{\alpha}'(B)+\dfrac{q^5}{q^2-1}W_{\beta}(B,0)-\dfrac{q^4}{q^2-1}W_{\gamma}(B,0)\rbrace,
		\end{equation}
		where $\beta=(-2,-2,0,2), \gamma=(-2,-2,2,2)$, and $B$ is the $4 \times 4$ matrix
		\begin{equation*}
			B=	\left(\begin{array}{cc} 
				h(x_i,x_j) & h(x_i,y_l)\\
				h(y_k,x_j)& h(y_k,y_l)
			\end{array} 
			\right)_{\substack{1\leq i,j \leq 2\\ 1 \leq k,l \leq 2}}.
		\end{equation*}
		Now, assume that $\val(h(y_1,y_1))=-1$ and that $B=Y_{\tau, \lambda}$, where $\tau(3)=3$ and $\lambda_3=-1$. Let $\overline{B}$ be the $3 \times 3$ matrix obtained by removing the $3$-rd row and the $3$-rd column from $B$. By Proposition \ref{proposition3.9}, $\CI_{\alpha}^{[\mathds{h}]}(x_1,x_2,y_1,y_2)$ is the arithmetic intersection number of special cycles $\CZ_2(x_1'),\CZ_2(x_2'),\CY_2(y_2')$ in $\CN_3^{[(0,2)]}$, where the hermitian matrix of $\lbrace x_1',x_2',y_2' \rbrace$ is $\overline{B}$. On the other hand, by Proposition \ref{proposition2.10}, the analytic side of \eqref{eq4.4.4} coincides with
		\begin{equation*}
			(q+1)\CI_{(-2,-2,2)}^{[(2)]}(x_1',x_2',y_2')=\dfrac{\tx{Vol}(P^{[(2)]})}{\tx{Vol}(P^{[(0,2)]})}	\CI_{(-2,-2,2)}^{[(2)]}(x_1',x_2',y_2').
		\end{equation*}
		Therefore, if \eqref{eq4.4.3} holds, then the conjectural formula \eqref{eq4.4.4} is consistent with the Kudla--Rapoport conjecture for maximal parahoric levels.
		
	\end{remark}

	\begin{lemma}\label{lemma3.2} For $\alpha,\beta\tx{'}s$, and $C_{\beta,i}\tx{'}s$ in Theorem \ref{theorem3.1}, we have
		\begin{equation*}
			q^{4n^2-2n}\CF_{\alpha}'(Y)-\mathlarger{\sum}_{\beta}C_{\beta,2n-1-\beta_{2n}}\CF_{\beta}(Y,0)=q^{4n^2-2n}\CF_{\alpha}'(Y^{\wedge_n})-\mathlarger{\sum}_{\beta}C_{\beta,2n-1-\beta_{2n}}\CF_{\beta}(Y^{\wedge_n},0).
		\end{equation*}
	\end{lemma}
	We will prove this lemma at the end of the subsection. We first introduce the following notation.
	
	\begin{definition}\label{definition3.3}
		\begin{enumerate}
			\item For $i \geq 0$, we define $C_{0,i}$ by
			\begin{equation*}
				C_{0,i}=\left\lbrace \begin{array}{ll}
					(-1)^{(i+1)/2}(\frac{1}{q^i+1}+\frac{1}{q^{i+1}-1}) & \tx{if }i \tx{ is odd},\\
					(-1)^{(i+2)/2}(\frac{1}{q^i-1}+\frac{1}{q^{i+1}+1}) & \tx{if }i \tx{ is even},\\
					0 & \tx{if }i=0.
				\end{array}\right.
			\end{equation*}
			
			\item We define a map $\mu:\CR_{2n} \rightarrow \lbrace 0,1 \rbrace^{2n}$ by
			\begin{equation*}
				\begin{array}{lcl}
					\mu:\CR_{2n} &\rightarrow & \lbrace 0,1\rbrace^{2n}\\
					Y=Y_{\sigma,e} & \mapsto & (\mu(Y)_i)_{1 \leq i \leq 2n},
				\end{array}
			\end{equation*}
			such that
			\begin{equation*}
				\mu(Y)_i= \left\lbrace \begin{array}{ll} 1 &\tx{if }i \in \lbrace j \mid j+\sigma(j) \geq 2n+1, e_j \leq 0\rbrace \cup \lbrace j \mid j+\sigma(j)\leq 2n, e_j \leq -1\rbrace \\
					0 & \tx{otherwise.} \end{array}\right.
			\end{equation*}
			\item For $Y \in \CR_{2n}$ we define $\kappa(Y)=\kappa(\mu(Y)):=\sum_{i=n+1}^{2n}\mu(Y)_i-\sum_{i=1}^{n}\mu(Y)_i.$
			
			\item For $\beta, \gamma \in \CE_n$ with
			\begin{equation*}
				\beta_i=-\beta_{2n+1-i}, \gamma_i=-\gamma_{2n+1-i}, \tx{ for all } i,
			\end{equation*}
			we define $\CF_{\beta,\gamma}(Y,0):=\prod_{j, \mu(Y)_j=1} (-q)^{\beta_j-\gamma_j}.$
			
		\end{enumerate}
	\end{definition}

	\begin{proposition}\label{proposition3.4}
		Let $Y=Y_{\sigma,e} \in \CR_{2n}$, $\alpha=(-(2n-1),-(2n-3),\dots,-1,1,3,\dots,2n-1)$, and let $\beta \in \CE_{2n}$ with $\beta_i=-\beta_{2n+1-i}$ for all $i$. Then, we have the following formulas.
		\begin{enumerate}
			\item $\CF'_{\alpha}(Y)=(\sum_{j=1}^{2n}\min(0,e_i))\CF_{\alpha}(Y,0).$
			
			\item $\CF'_{\alpha}(Y)-\CF'_{\alpha}(Y^{\wedge_n})=\kappa(Y)\CF_{\alpha}(Y,0)$.
			
			\item $q^{2n-2n\delta_{\beta}}\CF_{\beta}(Y,0)=\CF_{\alpha}(Y,0)\CF_{\alpha,\beta}(Y,0)$.
		\end{enumerate}

	\end{proposition}
	\begin{proof}
		\begin{enumerate}
			\item This follows directly from Proposition \ref{proposition2.6}.
			
			\item By Proposition \ref{proposition2.6}, we have $\CF_{\alpha}(Y,r)=\CF_{\alpha}(Y,0)f_{2r}(Y)$. Also, by Proposition \ref{proposition2.22} and the fact that $\alpha=\alpha^{\vee_n}$, we have
			\begin{equation*}
				\CF_{\alpha}(Y,0)=\CF_{\alpha^{\vee_n}}(Y^{\wedge_n},0)=\CF_{\alpha}(Y^{\wedge_n},0).
			\end{equation*}
			Therefore, by (1), we have
			\begin{equation*}
				\CF'_{\alpha}(Y^{\wedge_n})=\lbrace\sum_{i \in A^n(Y)} \min(0,e_i+1)+\sum_{i \in B^n(Y)} \min(0,e_i)+\sum_{i \in C^n(Y)} \min(0,e_i-1)\rbrace \CF_{\alpha}(Y,0).
			\end{equation*}
			Now, it suffices to show that
			\begin{equation}\label{eq3.3}
				\sum_{i \in A^n(Y)}\min(0,e_i)-\min(0,e_i+1) +\sum_{i \in C^n(Y)}\min(0,e_i)-\min(0,e_i-1)=\kappa(Y).
			\end{equation}
			
			If $i \in A^n(Y)$, we have that $i+\sigma(i) \leq 2n$, and hence $\mu(Y)_i=1$ if $e_i \leq -1$, and $\mu(Y)_i=0$ if $e_i \geq 0$. This shows that $-\mu(Y)_i=\min(0,e_i)-\min(0,e_i+1)$.
			
			Similarly, if $i \in C^n(Y)$, we have that $i+\sigma(i) \geq 2n$, and hence $\mu(Y)_i=1$ if $e_i \leq 0$, and $\mu(Y)_i=0$ if $e_i \geq 1$. This shows that $\mu(Y)_i=\min(0,e_i)-\min(0,e_i-1)$.
			
			Finally, we note that $\mu(Y)_i=\mu(Y)_{\sigma(i)}$ by definition. Therefore, if $i \in B^n(Y)$, then $\mu(Y)_i-\mu(Y)_{\sigma(i)}=0$. This finishes the proof of (2).
			
			\item By Proposition \ref{proposition2.6}, we have
			\begin{equation*}\begin{array}{ll}
					\CF_{\beta}(Y,0)=q^{-4n^2+2n\delta_{\beta}}f_{2n}(Y)&\prod_{\substack{1 \leq j \leq 2n\\j+\sigma(j) \leq 2n}}(-q)^{-(\beta_i+\beta_{\sigma(i)})(\min(0,e_j+1)-\min(0,e_j))/2}\\&\prod_{\substack{1 \leq j \leq 2n\\j+\sigma(j) \geq 2n+1}}(-q)^{-(\beta_i+\beta_{\sigma(i)})(\min(0,e_j)-\min(0,e_j-1))/2},
				\end{array}
			\end{equation*}
			\begin{equation*}\begin{array}{ll}
					\CF_{\alpha}(Y,0)=q^{-4n^2+2n}f_{2n}(Y)&\prod_{\substack{1 \leq j \leq 2n\\j+\sigma(j) \leq 2n}}(-q)^{-(\alpha_i+\alpha_{\sigma(i)})(\min(0,e_j+1)-\min(0,e_j))/2}\\&\prod_{\substack{1 \leq j \leq 2n\\j+\sigma(j) \geq 2n+1}}(-q)^{-(\alpha_i+\alpha_{\sigma(i)})(\min(0,e_j)-\min(0,e_j-1))/2}.
				\end{array}
			\end{equation*}
			Therefore, we have that
			\begin{equation*}\begin{array}{ll}
					q^{2n-2n\delta_{\beta}}\CF_{\beta}(Y,0)&=\CF_{\alpha}(Y,0)\prod_{\substack{1 \leq j \leq 2n\\j+\sigma(j) \leq 2n\\ e_j \leq -1}}(-q)^{(\alpha_i-\beta_{i}+\alpha_{\sigma(i)}-\beta_{\sigma(i)})/2}\prod_{\substack{1 \leq j \leq 2n\\j+\sigma(j) \leq 2n\\ e_j \leq -1}}(-q)^{(\alpha_i-\beta_{i}+\alpha_{\sigma(i)}-\beta_{\sigma(i)})/2}\\
					&=\CF_{\alpha}(Y,0)\CF_{\alpha,\beta}(Y,0).
				\end{array}
			\end{equation*}
			This finishes the proof of the proposition.
		\end{enumerate}

	\end{proof}
	
	Now, we introduce the following notation.
	
	\begin{definition}\label{definition3.5}
		\begin{enumerate}
			\item We define a subset $\CE_{2n,\alpha}$ of $\CE_{2n}$ by
			\begin{equation*}
				\CE_{2n,\alpha}:=\left\lbrace \beta \in \CE_{2n} \mid \begin{array}{lll} (1)\tx{ } \beta_i=-\beta_{2n+1-i} \text{ for all }1 \leq i \leq 2n.\\
					(2)\tx{ } \alpha_i-\beta_i \leq \alpha_{i+1}-\beta_{i+1} \text{ for all } 1\leq i \leq 2n-1. \end{array}\right\rbrace
			\end{equation*}
			\item For $0 \leq k \leq 2n-1$ and $Y \in \CR_{2n}$, we define $A_{k,2n}(Y)$ by
			\begin{equation*}
				A_{k,2n}(Y):=\sum_{\substack{\beta \in \CE_{2n,\alpha}\\\beta_{2n}=k}}\CF_{\alpha,\beta}(Y,0).
			\end{equation*}
			
			Similarly, for $\lambda \in \lbrace 0,1 \rbrace^{2n}$, we choose any $Y \in \CR_{2n}$ such that $\mu(Y)=\lambda$ and we define $A_{k,2n}(\lambda)$ by
			\begin{equation*}
				A_{k,2n}(\lambda)=A_{k,2n}(\mu(Y)):=A_{k,2n}(Y).
			\end{equation*}
			
			We also adopt the convention that $A_{-1,2n}(\lambda)=A_{-1,2n}(Y)=0$.
			
			We note that this definition does not depend on the choice of $Y$ since $\CF_{\alpha,\beta}(Y)$ depends only on $\mu(Y)$. We also note that $A_{2n-1,2n}(Y)=\CF_{\alpha,\alpha}(Y,0)=1$ for any $Y$.
			
			\item For $0 \leq k \leq 2n-1$ and $Y \in \CR_{2n}$, we define $A_{k,2n}^{\vee}(Y)$ by
			\begin{equation*}
				A_{k,2n}^{\vee}(Y):=\sum_{\substack{\beta \in \CE_{2n,\alpha}\\\beta_{2n}=k}}\CF_{\alpha,\beta^{\vee_n}}(Y,0).
			\end{equation*}
			Similarly, for $\lambda \in \lbrace 0,1 \rbrace^{2n}$, we choose any $Y \in \CR_{2n}$ such that $\mu(Y)=\lambda$ and we define $A_{k,2n}^{\vee}(\lambda)$ by
			\begin{equation*}
				A_{k,2n}^{\vee}(\lambda)=A_{k,2n}^{\vee}(\mu(Y)):=A_{k,2n}^{\vee}(Y).
			\end{equation*}
			We note that this definition does not depend on the choice of $Y$ since $\CF_{\alpha,\beta^{\vee_n}}(Y)$ depends only on $\mu(Y)$.
		\end{enumerate}
	\end{definition}
	
	To prove Lemma \ref{lemma3.2}, we need the following lemma.
	
	\begin{lemma}\label{lemma3.6}
		Assume that $0 \leq k \leq n-1$, and $\lambda \in \lbrace 0,1\rbrace^{2n}$. Then, we have the following formulas.
		\begin{enumerate}
			\item $A_{2k+1,2n}(\lambda)=(-q)^{-\kappa(\lambda)}A_{2k,2n}(\lambda)$.
			\item $A_{2k,2n}^{\vee}(\lambda)=(-q)^{-(2n-2k)\kappa(\lambda)}A_{2k,2n}(\lambda)$.
			\item $A_{2k+1,2n}^{\vee}(\lambda)=(-q)^{-(2n-2k-2)\kappa(\lambda)}A_{2k+1,2n}(\lambda)=(-q)^{-(2n-2k-1)\kappa(\lambda)}A_{2k,2n}(\lambda)$.
			\item $A_{2k+1,2n}(\lambda)=(-q)^{2\kappa(\overline{\lambda})}A_{2k+1,2n-2}(\overline{\lambda})+A_{2k-1,2n-1}(\overline{\lambda})$,
			where $\lambda=(\lambda_1,\lambda_2,\dots,\lambda_{2n}) \in \lbrace 0,1\rbrace^{2n}$, and $\overline{\lambda}=(\lambda_1,\dots,\lambda_{n-1},\lambda_{n+2},\dots,\lambda_{2n}) \in \lbrace 0,1\rbrace^{2n-2}$.
		\end{enumerate}
	\end{lemma}
	\begin{proof}
		\begin{enumerate}\label{eq4.4}
			\item For (1), we consider the bijection
			\begin{equation}\begin{array}{ccc}
					\lbrace \beta \in \CE_{2n,\alpha} \mid \beta_{2n}=2k \rbrace & \longrightarrow & \lbrace \beta \in \CE_{2n,\alpha} \mid \beta_{2n}=2k+1 \rbrace\\
					\beta =(\beta_1,\dots,\beta_{2n}) & \mapsto & \beta'=(\beta_1-1,\dots,\beta_n-1,\beta_{n+1}+1,\dots,\beta_{2n}+1).
				\end{array}
			\end{equation}
			Then, for any $Y \in \CR_{2n}$ with $\mu(Y)=\lambda$, we have
			\begin{equation*} 
				\CF_{\alpha,\beta}(Y,0)=(-q)^{-\kappa(\lambda)}	\CF_{\alpha,\beta'}(Y,0).
			\end{equation*}
			Since \eqref{eq4.4} is a bijection, we have
			\begin{equation*}
				A_{2k+1,2n}(\lambda)=(-q)^{-\kappa(\lambda)}A_{2k,2n}(\lambda).
			\end{equation*}
			This proves that (1) holds.
			
			\item For (2), we consider the set
			\begin{equation*}
				\ScS_{2n,t}:=\lbrace (\alpha_1-\beta_1,\dots,\alpha_{2n}-\beta_{2n}) \mid \beta \in \CE_{2n,\alpha},\tx{ }\beta_{2n}=t\rbrace.
			\end{equation*}
			Then, $\ScS_{2n,t}$ is the set of all $2n$-tuples of integers $(d_1,\dots,d_{2n})$ such that
			\begin{enumerate}
				\item $d_i=-d_{2n+1-i}$.
				\item $d_{2n}=2n-1-t$.
				\item $d_{n+1}=\left \lbrace \begin{array}{ll} 0 & \tx{ if }t \tx{ is odd},\\
					1 & \tx{ if }t \tx{ is even}.
				\end{array}\right.$
				\item $0 \leq d_{i+1}-d_i\leq 2$, $i \geq n+1$.
				\item $d_1 \equiv d_2\equiv \dots \equiv d_{2n}(\tx{mod }2)$.
			\end{enumerate}
			Similarly, we consider the set
			\begin{equation*}
				\ScS_{2n,t}^{\vee}:=\lbrace (\alpha_1-\beta^{\vee}_1,\dots,\alpha_{2n}-\beta^{\vee}_{2n}) \mid \beta \in \CE_{2n,\alpha},\tx{ }\beta_{2n}=t\rbrace.
			\end{equation*}
			Then, $\ScS_{2n,t}^{\vee}$ is the set of all $2n$-tuples of integers $(d_1',\dots,d_{2n}')$ such that
			\begin{enumerate}
				\item $d_i'=-d_{2n+1-i}'$.
				\item $d_{2n}'=-(2n-1-t)$.
				\item $d_{2n}'=\left \lbrace \begin{array}{ll} 0 & \tx{ if }t \tx{ is odd},\\
					-1 & \tx{ if }t \tx{ is even}.
				\end{array}\right.$
				\item $0 \leq d_{i+1}'-d_i'\leq 2$, $i \geq n+1$.
				\item $d_1' \equiv d_2\equiv \dots \equiv d_{2n}'(\tx{mod }2)$.
			\end{enumerate}
			We have a bijection between the two sets, given by
			\begin{equation*}
				\begin{array}{c}
					\ScS_{2n,t}  \longrightarrow  \ScS_{2n,t}^{\vee}\\
					(d_1,\dots,d_n,d_{n+1},\dots,d_{2n}) \mapsto\\
					\left\lbrace \begin{array}{ll}
						\scalebox{0.85}{$
							(d_1+(2n-1-t),\dots,d_n+(2n-1-t),
							d_{n+1}-(2n-1-t),\dots,d_{2n}-(2n-1-t))
							$}& \tx{if }t\tx{ is odd},\\
						\scalebox{0.85}{$
							(d_1+(2n-t),\dots,d_n+(2n-t),
							d_{n+1}-(2n-t),\dots,d_{2n}-(2n-t))
							$}& \tx{if }t\tx{ is even}.
					\end{array}\right.
				\end{array}
			\end{equation*}
			Now, (2) follows from this bijection.
			
			\item The statement in (3) follows from (1) and (2).
			
			\item For (4), we consider the subsets of $\CE_{2n,\alpha}$
			\begin{equation*}
				\CE_{2n,\alpha}=\lbrace \beta \in \CE_{2n,\alpha} \mid -\beta_{n-1}=\beta_{n+2}=0 \tx{ or } 1 \rbrace \sqcup \lbrace \beta \in \CE_{2n,\alpha} \mid -\beta_{n-1}=\beta_{n+2}\neq 0,1 \rbrace.
			\end{equation*}
			Then, for $\overline{\alpha}=(-(2n-3),\dots,-1,1,\dots,(2n-3))$, we have the following bijections:
			\begin{equation*}\begin{array}{ccl}
					\ell_1:\CE_{2n-2,\overline{\alpha}} & \rightarrow & \lbrace \beta \in \CE_{2n,\alpha} \mid -\beta_{n-1}=\beta_{n+2}=0 \tx{ or } 1 \rbrace\\
					(\overline{\beta}_1,\dots,\overline{\beta}_{2n-2}) & \mapsto & \left\lbrace \begin{array}{l}
						(\overline{\beta}_1,\dots,\overline{\beta}_{n-1},0,0,\overline{\beta}_n,\dots\overline{\beta}_{2n-2}) \tx{ if } \overline{\beta}_{2n-2} \tx{ is even},\\ (\overline{\beta}_1,\dots,\overline{\beta}_{n-1},-1,1,\overline{\beta}_n,\dots\overline{\beta}_{2n-2}) \tx{ if } \overline{\beta}_{2n-2} \tx{ is odd},
					\end{array}\right.
				\end{array}
			\end{equation*}
			
			\begin{equation*}\begin{array}{ccl}
					\ell_2:\CE_{2n-2,\overline{\alpha}} & \rightarrow & \lbrace \beta \in \CE_{2n,\alpha} \mid -\beta_{n-1}=\beta_{n+2}\neq 0,1 \rbrace\\
					(\overline{\beta}_1,\dots,\overline{\beta}_{n-2}) & \mapsto & \left\lbrace \begin{array}{l}
						(\overline{\beta}_1-2,\dots,\overline{\beta}_{n-1}-2,0,0,\overline{\beta}_n+2,\dots\overline{\beta}_{n-2}+2) \tx{ if } \overline{\beta}_{2n-2} \tx{ is even},\\ (\overline{\beta}_1-2,\dots,\overline{\beta}_{n-1}-2,-1,1,\overline{\beta}_n+2,\dots\overline{\beta}_{n-2}+2) \tx{ if } \overline{\beta}_{2n-2} \tx{ is odd}.
					\end{array}\right.
				\end{array}
			\end{equation*}
			
			Then, it is easy to see that
			\begin{equation*}
				\lbrace \beta \in \CE_{2n,\alpha} \mid \beta_{2n}=k \rbrace=\ell_1( \lbrace \beta \in \CE_{2n-2,\overline{\alpha}} \mid \overline{\beta}_{2n-2}=k \rbrace) \sqcup \ell_2( \lbrace \beta \in \CE_{2n-2,\overline{\alpha}} \mid \overline{\beta}_{2n-2}=k-2 \rbrace).
			\end{equation*}
			Therefore, the statement in (4) follows.
		\end{enumerate}
	\end{proof}

	Now, we are ready to prove Lemma \ref{lemma3.2}.
	\begin{proof}[Proof of Lemma \ref{lemma3.2}]
		By Proposition \ref{proposition3.4} (1), (2), we have that Lemma \ref{lemma3.2} is equivalent to
		\begin{equation*}
			\begin{array}{ll}
				&q^{4n^2-2n}(\CF_{\alpha}'(Y)-\CF_{\alpha}'(Y^{\wedge_n}))=\mathlarger{\sum}_{\beta}C_{\beta,2n-1-\beta_{2n}}(\CF_{\beta}(Y,0)-\CF_{\beta}(Y^{\wedge_n},0))\\
				\Longleftrightarrow & \kappa(Y)\CF_{\alpha}(Y,0)=\mathlarger{\sum}_{\beta}q^{-4n^2+2n}C_{\beta,2n-1-\beta_{2n}}((\CF_{\beta}(Y,0)-\CF_{\beta}(Y^{\wedge_n},0))).
			\end{array}
		\end{equation*}
		
		Note that
		\begin{equation*}
			q^{-4n^2+2n\delta_{\beta}}C_{\beta,2n-1-\beta_{2n}}=C_{0,2n-1-\beta_{2n}}.
		\end{equation*}
		Also, by Proposition \ref{proposition2.22} and Proposition \ref{proposition3.4}, we have
		\begin{equation*}
			\begin{array}{l}
				q^{2n-2n\delta_{\beta}}\CF_{\beta}(Y,0)=\CF_{\alpha}(Y,0)\CF_{\alpha,\beta}(Y,0),\\
				q^{2n-2n\delta_{\beta}}\CF_{\beta}(Y^{\wedge_n},0)=q^{2n-2n\delta_{\beta^{\vee_n}}}\CF_{\beta^{\vee_n}}(Y,0)=\CF_{\alpha}(Y,0)\CF_{\alpha,\beta^{\vee_n}}(Y,0).
			\end{array}
		\end{equation*}
		Therefore, Lemma \ref{lemma3.2} is equivalent to
		\begin{equation}\label{eq3.4}
			\kappa(Y)=\mathlarger{\sum}_{\beta}C_{0,2n-1-\beta_{2n}}(\CF_{\alpha,\beta}(Y,0)-\CF_{\alpha,\beta^{\vee_n}}(Y,0)).
		\end{equation}
		By Definition \ref{definition3.5}, \eqref{eq3.4} is equivalent to
		\begin{equation}\label{eq3.5}
			\begin{array}{ll}
				&\kappa(Y)=\sum_{k=0}^{2n-1} C_{0,2n-1-k}(A_{k,2n}(Y)-A_{k,2n}^{\vee}(Y)), \tx{ }\forall Y \in \CR_{2n}\\
				\Longleftrightarrow & \kappa(\lambda)=\sum_{k=0}^{2n-1} C_{0,2n-1-k}(A_{k,2n}(\lambda)-A_{k,2n}^{\vee}(\lambda)), \tx{ } \forall \lambda \in \lbrace 0,1 \rbrace^{2n}.
			\end{array}
		\end{equation}
		
		By Lemma \ref{lemma3.6} (1), (2), (3), we have that \eqref{eq3.5} is equivalent to
		\begin{equation}\label{eq3.6}
			\begin{array}{l}
				\kappa(\lambda)=C_{0,1}(1-(-q)^{-2\kappa(\lambda)})A_{2n-2,2n}(\lambda)\\+
				\mathlarger{\sum}_{k=0}^{n-2}\lbrace C_{0,2n-2k-1}(1-(-q)^{-(2n-2k)\kappa(\lambda)})+C_{0,2n-2k-2}((-q)^{-\kappa(\lambda)}-(-q)^{-(2n-2k-1)\kappa(\lambda)})\rbrace A_{2k,2n}(\lambda).
			\end{array}
		\end{equation}
		
		When $\kappa(\lambda)=0$, it is easy to see that \eqref{eq3.6} holds.
		
		Now, assume that $\kappa(\lambda) \neq 0$.
		
		Note that
		\begin{equation*}\begin{array}{l}
				C_{0,2n-2k-1}(1-(-q)^{-(2n-2k)\kappa(\lambda)})+C_{0,2n-2k-2}((-q)^{-\kappa(\lambda)}-(-q)^{-(2n-2k-1)\kappa(\lambda)})\\
				=(-1)^{n-k} (\frac{1}{q^{2n-2k-1}+1}+\frac{1}{q^{2n-2k}-1})(1-(-q)^{-(2n-2k)\kappa(\lambda)})\\
				+(-1)^{n-k}(\frac{1}{q^{2n-2k-2}-1}+\frac{1}{q^{2n-2k-1}+1})((-q)^{-\kappa(\lambda)}-(-q)^{-(2n-2k-1)\kappa(\lambda)})
				\\
				=(-1)^{n-k}\frac{1}{q^{2n-2k}-1}(1-(-q)^{-(2n-2k)\kappa(\lambda)})\\
				-(-1)^{n-k}((-q)^{-\kappa(\lambda)}+1)\frac{1}{q^{2n-2k-1}+1}(1-(-q)^{-(2n-2k-1)\kappa(\lambda)})\\
				+(-1)^{n-k}(-q)^{-\kappa(\lambda)}\frac{1}{q^{2n-2k-2}-1}(1-(-q)^{-(2n-2k-2)\kappa(\lambda)}).
			\end{array}
		\end{equation*}
		
		For $\kappa(\lambda) >0$, we define constants $d_{j,\kappa(\lambda)}(k)$, $0 \leq k \leq n-1$, $0 \leq j \leq \kappa(\lambda)-1$ by
		\begin{equation*}\begin{array}{l}
				d_{j,\kappa(\lambda)}(k):
				=(-1)^{n-k}(-q)^{-(2n-2k)(j+1)}\lbrace 1-((-q)^{-\kappa(\lambda)}+1)(-q)^{j+1}+(-q)^{-\kappa(\lambda)}(-q)^{2(j+1)}\rbrace.
			\end{array}
		\end{equation*}
		
		For $\kappa(\lambda) <0$, we define constants $d_{j,\kappa(\lambda)}(k)$, $0 \leq k \leq n-1$, $0 \leq j \leq -\kappa(\lambda)-1$ by
		\begin{equation*}\begin{array}{l}
				d_{j,\kappa(\lambda)}(k):=(-1)^{n-k}(-q)^{(2n-2k)j}\lbrace-1+((-q)^{-\kappa(\lambda)}+1)(-q)^{-j}-(-q)^{-\kappa(\lambda)}(-q)^{-2j}\rbrace.
			\end{array}
		\end{equation*}
		
		Then, for $0 \leq k \leq n-2$, we have
		\begin{equation}\label{eq3.7}
			\begin{array}{l}
				C_{0,2n-2k-1}(1-(-q)^{-(2n-2k)\kappa(\lambda)})+C_{0,2n-2k-2}((-q)^{-\kappa(\lambda)}-(-q)^{-(2n-2k-1)\kappa(\lambda)})\\
				=\left\lbrace \begin{array}{ll} \sum_{j=0}^{\kappa(\lambda)-1} d_{j,\kappa(\lambda)}(k), & \tx{ if } \kappa(\lambda) >0,\\
					\sum_{j=0}^{-\kappa(\lambda)-1} d_{j,\kappa(\lambda)}(k),& \tx{ if } \kappa(\lambda) <0.
				\end{array}\right.
			\end{array}
		\end{equation}	
		
		For $k=n-1$, we have
		\begin{equation}\label{eq3.8}
			C_{0,1}(1-(-q)^{-2\kappa(\lambda)})-\kappa(\lambda)(-q)^{-\kappa(\lambda)}=\left\lbrace
			\begin{array}{ll} \sum_{j=0}^{\kappa(\lambda)-1} d_{j,\kappa(\lambda)}(n-1), & \tx{ if } \kappa(\lambda) >0,\\
				\sum_{j=0}^{-\kappa(\lambda)-1} d_{j,\kappa(\lambda)}(n-1),& \tx{ if } \kappa(\lambda) <0.
			\end{array}\right.
		\end{equation}
		
		By \eqref{eq3.7} and \eqref{eq3.8}, we have that \eqref{eq3.6} is equivalent to
		\begin{equation*}
			0=\left\lbrace \begin{array}{ll}
				\sum_{i=0}^{n-1}\sum_{j=0}^{\kappa(\lambda)-1} d_{j,\kappa(\lambda)}(i)A_{2i,2n}(\lambda) & \tx{ if }\kappa(\lambda)>0,\\
				\sum_{i=0}^{n-1}\sum_{j=0}^{-\kappa(\lambda)-1} d_{j,\kappa(\lambda)}(i)A_{2i,2n}(\lambda) & \tx{ if }\kappa(\lambda)<0.
			\end{array}\right.
		\end{equation*}
		
		To prove this, it suffices to show that
		\begin{equation}\label{eq3.9}
			0=\sum_{i=0}^{n-1} d_{j,\kappa(\lambda)}(i)A_{2i,2n}(\lambda)
		\end{equation}
		when  $\kappa(\lambda)>0 \tx{, } 0 \leq j \leq \kappa(\lambda)-1$, or $\kappa(\lambda)<0 \tx{,}-\kappa(\lambda)+1 \leq j \leq 0$.
		
		For $\kappa(\lambda)>0$, we have that
		\begin{equation*}
			d_{\kappa(\lambda)-1,\kappa(\lambda)}(i)=(-1)^{n-i}(-q)^{-(2n-2i)(\kappa(\lambda)}\lbrace 1-((-q)^{-\kappa(\lambda)}+1)(-q)^{\kappa(\lambda)}+(-q)^{-\kappa(\lambda)}(-q)^{2\kappa(\lambda)}\rbrace=0.
		\end{equation*}
		Similarly, for $\kappa(\lambda)<0$, we have that $d_{0,\kappa(\lambda)}(i)=0.$ Therefore, \eqref{eq3.9} holds when $\kappa(\lambda)>0, j=\kappa(\lambda)-1$ or $\kappa(\lambda)<0, j=0$. Now, it remains to prove \eqref{eq3.9} when $\kappa(\lambda) \geq 2, 0 \leq j \leq \kappa(\lambda)-2$ or $\kappa(\lambda) \leq -2, 1 \leq j \leq -\kappa(\lambda)-1$.
		
		By the definition of $d_{j,\kappa(\lambda)}(i)$ and Lemma \ref{lemma3.6} (1), it suffices to prove that
		\begin{equation}\label{eq3.10}
			0=\sum_{i=0}^{n-1}(-1)^{n-i}(-q)^{-(2n-2i)j}A_{2i,2n}(\lambda)=\sum_{i=0}^{n-1}(-1)^{n-i}(-q)^{-(2n-2i)j}A_{2i+1,2n}(\lambda)(-q)^{\kappa(\lambda)},
		\end{equation}
		when $\kappa(\lambda)>0, 1 \leq j \leq \kappa(\lambda)-1$, or $\kappa(\lambda)<0, \kappa(\lambda)+1 \leq j \leq -1$.
		
		We prove \eqref{eq3.10} by induction on $n$.
		
		Since $-n \leq \kappa(\lambda) \leq n$, there is nothing to prove when $n=1$.
		
		When $n=2$, we only need to prove \eqref{eq3.10} in the following two cases: $(\kappa(\lambda),j)=(2,1)$, $(-2,-1)$.
		
		By Lemma \ref{lemma3.6} (4) and Definition \ref{definition3.5} (2), we have
		\begin{equation*}\begin{array}{l}
				(1) A_{1,4}(\lambda)=(-q)^{2\kappa(\overline{\lambda})}A_{1,2}(\overline{\lambda})=\left\lbrace \begin{array}{ll}
					(-q)^2 &\tx{if }\kappa(\lambda)=2,j=1,\\
					(-q)^{-2} &\tx{if }\kappa(\lambda)=-2, j=-1,
				\end{array}\right.\\
				(2) A_{3,4}(\lambda)=1,
			\end{array}
		\end{equation*}
		where $\lambda=(\lambda_1,\lambda_2,\lambda_3,\lambda_4)$, and $\overline{\lambda}=(\lambda_1,\lambda_4)$.
		Therefore, we have
		\begin{equation*}
			\sum_{i=0}^{n-1}(-1)^{n-i}(-q)^{-(2n-2i)j}A_{2i+1,2n}(\lambda)=(-q)^{-4j}A_{1,4}(\lambda)-(-q)^{-2j}A_{3,4}(\lambda)=0.
		\end{equation*}
		This shows that \eqref{eq3.10} holds for $n=2$.
		
		Now, we assume that \eqref{eq3.10} holds for $n=N-1$. By Lemma \ref{lemma3.6} (4), we have
		\begin{equation*}
			\begin{array}{l}
				\sum_{i=0}^{N-1}(-1)^{N-i}(-q)^{-(2N-2i)j}A_{2i+1,2N}(\lambda)\\
				=\sum_{i=0}^{N-1}(-1)^{N-i}(-q)^{-(2N-2i)j+2\kappa(\overline{\lambda})}A_{2i+1,2N-2}(\overline{\lambda})+\sum_{i=0}^{N-1}(-1)^{N-i}(-q)^{-(2N-2i)j}A_{2i-1,2N-2}(\overline{\lambda}),
			\end{array}
		\end{equation*}
		where $\lambda=(\lambda_1,\dots,\lambda_{2N})$, and $\overline{\lambda}=(\lambda_1,\dots,\lambda_{N-1},\lambda_{n+2},\dots\lambda_{2N})$.
		
		Since $A_{2N-1,2N-2}(\overline{\lambda})=A_{-1,2N-2}(\overline{\lambda})=0$, we have
		\begin{equation}\label{eq3.11}
			\begin{array}{l}
				\sum_{i=0}^{N-1}(-1)^{N-i}(-q)^{-(2N-2i)j}A_{2i+1,2N}(\lambda)\\
				=\sum_{i=0}^{N-2}(-1)^{N-i}(-q)^{-(2N-2i)j+2\kappa(\overline{\lambda})}A_{2i+1,2N-2}(\overline{\lambda})+\sum_{i=0}^{N-2}(-1)^{N-i-1}(-q)^{-(2N-2i-2)j}A_{2i+1,2N-2}(\overline{\lambda})\\
				=\sum_{i=0}^{N-2}(-1)^{N-1-i}(-q)^{-(2(N-1)-2i)j}A_{2i+1,2N-2}(\overline{\lambda})\lbrace 1-(-q)^{2\kappa(\overline{\lambda})-2j}\rbrace.
			\end{array}
		\end{equation}
		
		Note that $\kappa(\overline{\lambda})=\kappa(\lambda)-1, \kappa(\lambda),$ or $\kappa(\lambda)+1$.
		
		We first assume that $\kappa(\lambda)\geq 2$, $1 \leq j \leq \kappa(\lambda)-1$.
		
		If $\kappa(\overline{\lambda}) \geq \kappa(\lambda)$, then $1 \leq j \leq \kappa(\lambda)-1 \leq \kappa(\overline{\lambda})-1$, and hence
		\begin{equation*}
			\sum_{i=0}^{n-2}(-1)^{n-1-i}(-q)^{-(2(n-1)-2i)j}A_{2i+1,2n-2}(\overline{\lambda})=0,
		\end{equation*}
		by our inductive hypothesis.
		
		If $\kappa(\overline{\lambda})=\kappa(\lambda)-1$, then for $1 \leq j \leq \kappa(\lambda)-2=\kappa(\overline{\lambda})-1$, we have that
		\begin{equation*}
			\sum_{i=0}^{n-2}(-1)^{n-1-i}(-q)^{-(2(n-1)-2i)j}A_{2i+1,2n-2}(\overline{\lambda})=0,
		\end{equation*}
		by our inductive hypothesis.
		
		For $j=\kappa(\lambda)-1=\kappa(\overline{\lambda})$, we have that $\lbrace 1-(-q)^{2\kappa(\overline{\lambda})-2j}\rbrace=0$.
		
		Therefore, by \eqref{eq3.11}, we have that \eqref{eq3.10} holds when $\kappa(\lambda)>0$, $1 \leq j \leq \kappa(\lambda)-1$.
		
		Now, assume that $\kappa(\lambda) \leq -2$, $\kappa(\lambda)+1 \leq j \leq -1$.
		
		If $\kappa(\overline{\lambda}) \leq \kappa(\lambda)$, then $\kappa(\overline{\lambda})+1 \leq \kappa(\lambda)+1 \leq j \leq -1$, and hence
		\begin{equation*}
			\sum_{i=0}^{n-2}(-1)^{n-1-i}(-q)^{-(2(n-1)-2i)j}A_{2i+1,2n-2}(\overline{\lambda})=0,
		\end{equation*}
		by our inductive hypothesis.
		
		If $\kappa(\overline{\lambda})=\kappa(\lambda)+1$, then for $\kappa(\lambda)+2=\kappa(\overline{\lambda})+1 \leq j \leq -1$, we have
		\begin{equation*}
			\sum_{i=0}^{n-2}(-1)^{n-1-i}(-q)^{-(2(n-1)-2i)j}A_{2i+1,2n-2}(\overline{\lambda})=0,
		\end{equation*}
		by our inductive hypothesis.
		
		For $j=\kappa(\lambda)+1=\kappa(\overline{\lambda})$, we have that  $\lbrace 1-(-q)^{2\kappa(\overline{\lambda})-2j}\rbrace=0$.
		
		Therefore, by \eqref{eq3.11}, we have that \eqref{eq3.10} holds when $\kappa(\lambda)<0$, $\kappa(\lambda)+1 \leq j \leq -1$.
		
		This finishes the proof of the lemma.
	\end{proof}

	\section{Proof of $n=2$ case}\label{section5}
		In this section, we always assume $n=2$ and $\mathds{h}=(0,2)$. We will drop the subscript $n$ in $\cN_n^{[\mathds{h}]}$ for simplicity.

\subsection{Pure weight case}
In this subsection, we consider the case  $\alpha=(0,0)$ or $(2,2)$.

    The following is a special case of  \cite[Proposition 3.5]{LiMihatsch}.
	\begin{proposition}\label{proposition5.5.1}
		The natural forgetful morphism 
		\begin{align*}
			\pi_0: \cN^{[\mathds{h}]}\to \cN^{[0]}, \quad \text{ and } \quad  \pi_2: \cN^{[\mathds{h}]}\to \cN^{[2]}
		\end{align*}
		are finite and flat of degree $q+1$.
	\end{proposition}
	\begin{proof}
		Since $\cN^{[\mathds{h}]}$, $\cN^{[0]}$ and  $\cN^{[2]}$  are regular,  the flatness follows from the miracle flatness theorem. Moreover, since it is well-known that $\pi$ is proper, one can also conclude that $\pi$ is finite once we know it is quasi-finite. Note that there is natural identification between $\cN^{[0]}$ and $\cM_2$ via the Serre construction, where $\cM_2$ is the Lubin-Tate space parameterizes formal $p$-divisible groups of dimension 1 and height 2. Moreover, $\cN^{[\mathds{h}]}$ is also naturally identified with $\cM_{\Gamma_0(p)}$. Then the flatness and finiteness properties of $\pi_0$ and $\pi_2$ are proved in \cite[Proposition 3.5]{LiMihatsch}.  Then the fact that $\pi_0$ and $\pi_2$ have degree $q+1$ can be checked easily on the generic fiber.
	\end{proof}
	Using Proposition \ref{proposition5.5.1}, we can verify a variant of Conjecture \ref{conjecture4.3}.
	\begin{theorem}Assume $B$ is a non-degenerate hermitian matrix of rank $2$ and $\alpha=(0,0)$. We have
	\begin{align*}
\CI_{\alpha}^{[\mathds{h}]}(B)= (q+1)   
\partial \operatorname{Den}_{2,0}(B).\end{align*}
Here $\partial\operatorname{Den}_{2,0}(B)$ is the normalized derived local density as in \cite{chz}.
	\end{theorem}
    \begin{proof}
      Using Proposition \ref{proposition5.5.1} and the projection formula, we know that $\CI_{\alpha}^{[\mathds{h}]}(B)=(q+1)\CI_{\alpha}^{[0]}(B)$. Then the theorem reduces to the Kudla-Rapoport conjecture proved in \cite{KR1,LZ}.
    \end{proof}

\subsection{Mixed weight case}

     In this subsection,   we assume that $\alpha=(0,2) \in \CE_{2}^{\mathds{h}}$. Let $\lbrace x ,y \rbrace$ be a basis of the space of special homomorphisms $(\BV,h)$ associated to the Rapoport-Zink space $\CN^{[\mathds{h}]}_2$.	Let $B=\left(\begin{array}{cc} h(x,x)& h(x,y)\\h(y,x)&h(y,y) \end{array}\right).$ We note that the valuation of $\det B$ is odd. Also, by Definition \ref{definition2.5}, $\CR_2$ form a complete set of representatives of $\Gamma_2 \backslash X_2$. Therefore, $B$ is $\Gamma_2$-equivalent to $\left(\begin{array}{cc} \pi^{\lambda_1}& \\&\pi^{\lambda_2} \end{array}\right)$ for some $\lambda_1, \lambda_2 \in \BZ$.

     Combining Conjecture \ref{conjecture4.3} and Remark \ref{remark4.5}, we have the following conjecture.
	\begin{conjecture}\label{conjecture4.12}
		Assume that $n=2$, $\mathds{h}=(0,2)$, and $\alpha=(0,2) \in \CE_{2}^{\mathds{h}}$. Then, we have
		\begin{equation*}
			\CI_{\alpha}^{[\mathds{h}]}(x,y)=\dfrac{q^7}{q^2-1} W_{\alpha}'(B).
		\end{equation*}
	\end{conjecture}
    \begin{proof}
        This follows from Proposition \ref{proposition4.13} and Proposition \ref{prop: div intersection} below.
    \end{proof}
    First, we explicitly compute the analytic side.
	\begin{proposition}\label{proposition4.13}
		Assume that $n=2$, $\mathds{h}=(0,2)$, and $\alpha=(0,2) \in \CE_{2}^{\mathds{h}}$. Then, we have
		\begin{equation*}\begin{array}{l}
				\dfrac{q^7}{q^2-1} W_{\alpha}'(B)
				=\left\lbrace \begin{array}{cc}
					\dfrac{q+1}{2}\mathlarger{\sum}_{l=0}^{\lambda_2+1} q^l(\lambda_1+\lambda_2-2l+1)+\mathlarger{\sum}_{l=0}^{\lambda_2+1}q^l & \tx{if }B \overset{\Gamma_2}{\simeq} \diag(\pi^{\lambda_1},\pi^{\lambda_2}), \lambda_1 >\lambda_2,\\
					\dfrac{q+1}{2}\mathlarger{\sum}_{l=0}^{\lambda_1} q^l(\lambda_1+\lambda_2-2l+1)+\mathlarger{\sum}_{l=0}^{\lambda_1}q^l & \tx{if }B \overset{\Gamma_2}{\simeq} \diag(\pi^{\lambda_1},\pi^{\lambda_2}), \lambda_1 <\lambda_2.
				\end{array}\right.
			\end{array}
		\end{equation*}
	\end{proposition}
	\begin{proof}
		For $\beta=(0,0) \in \CE_{2}^{\mathds{h}}$, $W_{\beta}'(B)$ is a usual representation density with respect to our chosen measure. Therefore, by \cite[Corollary 9.5]{KR1}, we have
		\begin{equation*}
			W_{\beta}'(B)=\left\lbrace \begin{array}{ll}\dfrac{1}{q^4}(1+\dfrac{1}{q})(1-\dfrac{1}{q^2})\dfrac{1}{2}\mathlarger{\sum}_{l=0}^{\min(\lambda_1,\lambda_2)}q^{l}(\lambda_1+\lambda_2-2l+1) & \tx{if }\lambda_1,\lambda_2 \geq 0,\\
				0 & \tx{otherwise}.
			\end{array}\right.
		\end{equation*}
		By Proposition \ref{proposition2.3}, we have
		\begin{equation*}
			W_{\alpha}'(B)-W_{\beta}'(B)=\mathlarger{\sum}_{Y \in \CR_2}\dfrac{\CG(Y,B)}{\alpha(Y;2)}\lbrace \CF_{\alpha}'(Y)-\CF_{\beta}'(Y)\rbrace.
		\end{equation*}
		By Proposition \ref{proposition2.6}, we have
		\begin{equation*}
			\CF_{\alpha}'(Y)-\CF_{\beta}'(Y)=\left\lbrace \begin{array}{ll}
				e_1 (-q)^{2e_1}(\dfrac{1}{q^2}-\dfrac{1}{q^4}) & \tx{if }Y=\diag(\pi^{e_1},\pi^{e_2}), e_1 \leq -1, e_2 \geq 1,\\
				0 & \tx{otherwise}.
			\end{array}\right.
		\end{equation*}
		Therefore, we only need to consider $Y=\diag(\pi^{e_1},\pi^{e_2}) \in \CR_2$, $e_1 \leq -1, e_2 \geq 1$. For such $Y$, we have
		\begin{equation*}
			\begin{array}{l}
				\CG(Y,B)=(1-\dfrac{1}{q^2})^2q^{-2}I_1^*(e_1+\lambda_1)(-q)^{\min(0,e_1+\lambda_2+2)},\\
				\alpha(Y;2)=(q+1)^2q^{-2+3e_1+e_2}.
			\end{array}
		\end{equation*}
		Combining these, we prove the proposition.
	\end{proof}

	To compute $\CI_{\alpha}^{[\mathds{h}]}(x,y)$, we first review some relevant background about quasi-canonical liftings following \cite{KR1}.
	Let $G$ be a strict formal $O_{F_0}$-module   of relative height 2 and dimension 1 over $\mathbb{F}$.   Note that $\operatorname{End}(G)=O_D$, the maximal order in the quaternion division algebra $D$ over $F_0$. If we fix an embedding $i: O_F \hookrightarrow O_D,$ then  $G$ becomes a strict formal $O_F$-$O_{F_0}$-module of height 2. There is a unique lifting $F_0$ of this formal $O_F$-$O_{F_0}$-module to $O_{\breve{F}_0}$ which we call the canonical lift. Let
	$$
	O_{F, s}=O_{F_0}+\pi^s O_F .
	$$
	be the order of conductor $s$ in $O_F$  For an integer $s \geq 1$. To simplify the notation, we will use $O$ to denote $O_{\breve{F}_0}$ from now on.
	
	  Let $F_s$ denote the  quasi-canonical lift  of level $s$ of $G$. More precisely, it is a lifting of $G$ to some finite extension $A$ of $O_{F_0}$ with endomorphism ring equal to $O_{F, s}$, such that the induced action of $O_{F, s}$ on Lie $F_s$ is through the embedding $O_{F, s} \subset O_{\breve{F_0}} \subset A$ and such that the image of $O_{F, s}$ in $\operatorname{End}(G)$ is contained in the image of the fixed embedding of $O_F$ in $O_D$. Any quasi-canonical lift of level $s$ is defined over $O_s$, which is the ring of integers in the ring class field $M_s$ of $O_{F, s}^{\times}$. More precisely, $M_s$ is the finite abelian extension of $M=\breve{F}_0$ corresponding to the subgroup $O_{F, s}^{\times}$ of $O_F^{\times}$ under the reciprocity isomorphism. Note that $M_s$ is a totally ramified extension of $M$ of degree
	$$
	e_s=q^{s-1}(q+1).
	$$

    Now let $X^{(s)}=O_F \otimes F_s$, which is a strict formal $O_F$-$O_{F_0}$-module   over $O_s$. Moreover, Lie $X^{(s)}=O_F \otimes_{O_{F_0}}$ Lie $F_s$, hence $X^{(s)}$ satisfies the signature condition $(1,1)$. In fact, there exists $
\left(X^{(s)}, \iota, \lambda^{(s)}\right)$, which can be regarded as a $O_s$-point of $\cN_2^{[0]}$ by  \cite[Lemma 6.5]{KR1}.  
	We use $\cZ_{s}^{[0]}$ to denote this $W_s$-point of $\cN_2^{[0]}$. Via the duality isomorphism $\cN_2^{[0]}\cong \cN_2^{[2]}$, we also obtain a corresponding $W_s$-point on  $\cN_2^{[2]}$ which we denote as $\cZ_{s}^{[2]}$.  Let $\cZ_{s}^+$ denote   $\cZ_{s}^{[0]}\stackrel{\alpha_{0,2}^+}{\to} \cZ_{s+1}^{[2]}$ which we regard as an $O_{s+1}$-point of $\cN_2^{[\mathds{h}]}$. Similarly, let $\cZ_{s}^{-}$ denote $\cZ_{s}^{[0]}\stackrel{\alpha_{0,2}^-}{\to}  \cZ_{s-1}^{[2]}$ which we regard as an $O_{s}$-point of $\cN_2^{[\mathds{h}]}$ again.

	\begin{proposition}\label{prop: pullback along pi_0}
		We have
		\begin{align*}
			\pi_0^*(\cZ_{s}^{[0]})&=\cZ_{s}^{+}+\cZ_{s}^{-},\\
			\pi_2^*(\cZ_{s}^{[2]})&=\cZ_{s-1}^{+}+\cZ_{s+1}^{-}.
		\end{align*}
	\end{proposition}
	\begin{proof}
		It is easy to check that $\cZ_{s}^{+}$ and $\cZ_{s}^{-}$ are contained in $\pi_0^*(\cZ_{s}^{[0]})$. Then one can conclude the identity by the fact that $\pi_0$ is of degree $q+1$, $\cZ_{s}^{[0]}\stackrel{\alpha_{0,2}^+}{\to} \cZ_{s+1}^{[2]}$ is  an $O_{s+1}$-point of $\cN_2^{[\mathds{h}]}$ , $\cZ_{s}^{[0]}\stackrel{\alpha_{0,2}^-}{\to}  \cZ_{s-1}^{[2]}$is an $O_{s}$-point of $\cN_2^{[\mathds{h}]}$, and $e_{s+1}/e_s=q$. 
		
		Similalry, we know $\cZ_{s-1}^{+}$ and  $\cZ_{s+1}^{-}$ are contained in $	\pi_2^*(\cZ_{s})$ and we can conclude the second statement by comparing degrees again.
	\end{proof}
Assume $\val(x)=s$. To simplify the notations, we let $\cZ_i(x)$ denote $\cZ_i^{[\mathds{h}]}$. By \cite{KR1}[Proposition 8.1],  we have 
\begin{align}\label{eq: decomposition}
\cZ^{[0]}(x)-\cZ^{[0]}(\pi^{-1}x)=\cZ_s^{[0]}.
\end{align} 
Then $$\cZ_0(x)-\cZ_0(\pi^{-1}x)=\pi_0^*(\cZ^{[0]}(x)-\cZ^{[0]}(\pi^{-1}x))=\cZ_{s}^{+}+\cZ_{s}^{-}$$
by Proposition \ref{prop: pullback along pi_0}.

Note that  $\cZ_0(\pi^{-1}x)\subset \cY_2(\pi^{-1}x) \subset \cZ_0(x)$. And it is natural to ask what  are $   \cZ_0(x) -\cY_2(\pi^{-1}x)$ and 
$\cY_2(\pi^{-1}x)-\cZ_0(\pi^{-1}x)$. The following proposition answers this.
\begin{proposition}\label{prop: diff div as quasicanonical lifting}
	Assume $\val(x)=s\ge 1$. We have 
	\begin{align*}
	 \cZ_0(x) -\cY_2(\pi^{-1}x)=\cZ_s^+ ,\\
	 \cY_2(\pi^{-1}x)-\cZ_0(\pi^{-1}x)=\cZ_s^-.
	\end{align*}
\end{proposition}
\begin{proof}
Assume $X_2\to X_0$ lies on $\cY_2(\pi^{-1}x)-\cZ_0(\pi^{-1}x)\subset \cY_2(\pi^{-1}x)-\cY_2(\pi^{-2}x) $. In particular, $X_2 \in \cY^{[2]}(\pi^{-1}x)-\cY^{[2]}(\pi^{-2}x) \subset \cN_2^{[2]}$. Note that $\cY^{[2]}(\pi^{-1}x)-\cY^{[2]}(\pi^{-2}x)\cong  \cZ_{s-1}^{[2]}$ by duality and \eqref{eq: decomposition}. Therefore, we know $\cZ_s^{-}\subset \cY_2(\pi^{-1}x) -\cZ_0(\pi^{-1}x) $ and $\cZ_s^{+}\not \subset \cY_2(\pi^{-1}x) -\cZ_0(\pi^{-1}x) $. On the other hand, we know $( \cZ_0(x) -\cY_2(\pi^{-1}x))+
(\cY_2(\pi^{-1}x)-\cZ_0(\pi^{-1}x))=\cZ_s^+ +\cZ_s^-.$ Therefore, $ \cY_2(\pi^{-1}x)-\cZ_0(\pi^{-1}x)=\cZ_s^-$ and $ \cZ_0(x) -\cY_2(\pi^{-1}x)=\cZ_s^+ $.
\end{proof}

To compute the intersection between special cycles, we recall the following lifting problem following \cite[\S 7]{KR1}.  Let $A$ be a finite extension of $O$ with uniformizer $\lambda$, and set $A_m=A / \lambda^{m+1}$. Let $e$ denote the ramification index of $A$ over $O$ and denote by $\val_A$ the discrete valuation on $A$ with $\val_A(\lambda)=1 / e$.   Let
$$
\mu:\left(O_F \otimes F_0\right) \otimes_W \mathbb{F} \longrightarrow\left(O_F \otimes F_s\right) \otimes_{W_s} \mathbb{F}
$$
be a homomorphism. Let $m_s(\mu)$ be the maximum $m$ such that $\mu$ lifts to a homomorphism from $\left(O_F \otimes F_0\right) \otimes_W A_m$ to $\left(O_F \otimes F_s\right) \otimes_{W_s} A_m$. Also, as in the first part of (iii) of Lemma 6.1, write
$$
\mu=\left(\begin{array}{ll}
	\mu_1 & \mu_2 \\
	\mu_3 & \mu_4
\end{array}\right),
$$
with $\mu_i \in \operatorname{Hom}\left(F_0 \otimes_W \mathbb{F}, F_s \otimes_{W_s} \mathbb{F}\right)=O_D$.
The following is \cite[Theorem 7.1]{KR1}.
\begin{theorem}\label{thm: lifting}
 Write $\mu$ as in (7.1), and suppose that
$$
\mu_i \in\left(\Pi^s O_F+\Pi^{l_i} O_D\right) \backslash\left(\Pi^s O_F+\Pi^{l_i+1} O_D\right)
$$
for integers $l_i \geq 0$. Let $l=\min _i\left\{l_i\right\}$. Then
$$
m_s(\mu)=e / e_s \cdot\left\{\begin{array}{ll}
	\frac{q^{l+1}-1}{q-1} & \text { if } l<s, \\
	\frac{q^s-1}{q-1}+\frac{1}{2}(l+1-s) e_s & \text { if } l \geq s .
\end{array}\right.
$$
\end{theorem}

To simplify notations, from now on, we use $\cdot $ to denote the intersection number we defined previously using the Euler characteristic of derived tensor product.
\begin{proposition}\label{prop: div intersection}
	Assume $\val(x)=a$ and $\val(y)=b$ such that  $a>b$. Then we have
	\begin{align*}
	 (\cY_2(x)-\cZ_0(x))\cdot  (\cY_2(y)-\cY_2(\pi^{-1}y))=(q+1)q^{b+1},\\
	 (\cZ_0(x)-\cY_2(\pi^{-1}x))\cdot  (\cY_2(y)-\cY_2(\pi^{-1}y))=(q+1)q^{b}.
	\end{align*}
\end{proposition}
\begin{proof}
	By Proposition \ref{prop: diff div as quasicanonical lifting}, we have 
   \begin{align*}
	 \cZ_0(x) -\cY_2(\pi^{-1}x)&=\cZ_a^+ ,\\
	 \cY_2(x)-\cZ_0(x)&=\cZ_{a+2}^-.
	\end{align*}
By the duality isomorphism, 
\begin{align*}
     (\cY_2(x)-\cZ_0(x))\cdot  (\cY_2(y)-\cY_2(\pi^{-1}y))
     &= \cZ_{a+2}^- \cdot  (\cY_2(y)-\cY_2(\pi^{-1}y))\\
     &= \cZ_{a+1}^+ \cdot  (\cZ_0(y')-\cZ_0(\pi^{-1}y')),
\end{align*}
where $\val(y')=b+1$. 
Now we can combine \cite[Theorem 7.1]{KR1} with the same proof of \cite[Proposition 8.4]{KR1} to compute $\cZ_{a+1}^+ \cdot  (\cZ_0(y')-\cZ_0(\pi^{-1}y'))$. The only difference is when we apply  \cite[Theorem 7.1]{KR1} here, we need to take the $A$ in  \cite[Theorem 7.1]{KR1} to be $O_{a+2}$ instead of $O_{a+1}$ since $\cZ_{a+1}^+$ is defined over $O_{a+2}$. Therefore, we get
    \begin{align*}
        \cZ_{a+1}^+ \cdot  (\cZ_0(y')-\cZ_0(\pi^{-1}y'))&= 
e_{a+2}/e_{a+1}(\frac{q^{b+2}-1}{q-1}-\frac{q^{b}-1}{q-1})=(q^{b+1}+q^b)\\
&=e_{a+2}/e_{a+1}(q+1)q^b\\&=(q+1)q^{b+1}.
    \end{align*}

Similarly,   
\begin{align*}
     (\cZ_0(x)-\cY_2(\pi^{-1}x))\cdot  (\cY_2(y)-\cY_2(\pi^{-1}y))&=\cZ_a^+\cdot (\cY_2(y)-\cY_2(\pi^{-1}y))\\
     &= \cZ_{a+1}^-\cdot (\cZ_0(y')-\cZ_0(\pi^{-1}y')).
\end{align*}
Note that $\cZ_{a+1}^-$ is defined over $O_{a+1}$. By the same argument as before, we know 
    \begin{align*}
        \cZ_{a+1}^- \cdot  (\cZ_0(y')-\cZ_0(\pi^{-1}y'
        ))&= 
e_{a+1}/e_{a+1}(\frac{q^{b+2}-1}{q-1}-\frac{q^{b}-1}{q-1})=(q^{b+1}+q^b)\\
 &=(q+1)q^{b}.
    \end{align*}
\end{proof}

	\bibliographystyle{alpha}
	\bibliography{reference}
\end{document}